\documentclass{amsart}
\usepackage{amsmath,amssymb, amsbsy}
\usepackage[dvips]{graphicx}
\newcommand{\R}{{\mathbb R}}
\newcommand{\N}{{\mathbb N}}
\newcommand{\Z}{{\mathbb Z}}
\newcommand{\C}{{\mathbb C}}
\newcommand{\SN}{{\mathbb S}^{N-1}}
\newcommand{\weakly}{\rightharpoonup}
\newcommand{\e }{\varepsilon}
\newcommand{\Di}{{\mathcal D}^{1,2}(\R^N,\C)}
\newcommand{\Distar}{{\mathcal D}^{1,2}_{*}(\R^N,\C)}
\newcommand{\DiA}{{\mathcal D}^{1,2}_{\A}(\R^N)}
\newcommand{\DiAa}{{\mathcal D}^{1,2}_{\A,a}(\R^N)}

\newcommand{\alchi}{\raisebox{1.7pt}{$\chi$}}
\newcommand{\dive }{\mathop{\rm div}}

\newcommand{\A}{{\mathbf A}}

\renewcommand{\geq }{\geqslant}
\renewcommand{\le }{\leqslant}
\renewcommand{\leq }{\leqslant}

\newenvironment{pf}{\noindent{\sc Proof}.\enspace}{\hfill\qed\medskip}
\newenvironment{pfn}[1]{\noindent{\bf Proof of
    {#1}.\enspace}}{\hfill\qed\medskip}
\newtheorem{Theorem}{Theorem}[section]
\newtheorem{Corollary}[Theorem]{Corollary}
\newtheorem{Lemma}[Theorem]{Lemma}
\newtheorem{Proposition}[Theorem]{Proposition}
\theoremstyle{definition} 

\newtheorem{remark}[Theorem]{Remark}
\begin{document}

\title[Schr\"odinger equations with singular electromagnetic
potentials]{Asymptotic behavior of solutions to Schr\"odinger
  equations near an isolated singularity of the electromagnetic potential}

\author[Veronica Felli \and Alberto Ferrero \and Susanna
Terracini]{Veronica Felli \and Alberto Ferrero \and Susanna Terracini}
\address{\hbox{\parbox{5.7in}{\medskip\noindent{Universit\`a di Milano
        Bicocca,\\
        Dipartimento di Ma\-t\-ema\-ti\-ca e Applicazioni, \\
        Via Cozzi
        53, 20125 Milano, Italy. \\[3pt]
        \em{E-mail addresses: }{\tt veronica.felli@unimib.it,
          alberto.ferrero@unimib.it, susanna.terracini@unimib.it}.}}}}

\date{Revised version, April 7, 2009}

\thanks{First and third author supported by Italy MIUR, national
  project ``Variational Methods and Nonlinear Differential
  Equations''.
  }

 \begin{abstract}
   \noindent Asymptotics of solutions to Schr\"odinger equations with
   singular magnetic and electric potentials is investigated.  By
   using a Almgren type monotonicity formula, separation of variables,
   and an iterative Brezis-Kato type procedure, we describe the
   exact behavior near the singularity of solutions to linear and
   semilinear (critical and subcritical) elliptic equations with an
   inverse square electric potential and a singular magnetic potential with a
   homogeneity  of order $-1$.
 \end{abstract}

\maketitle

% \tableofcontents
%  \newpage

\bigskip
\section{Introduction}\label{intro}

In quantum mechanics, the hamiltonian of a non-relativistic charged
particle in an electromagnetic field has the form
$(-i\nabla+{\boldsymbol{\mathcal A}})^2+V$, where $V:\R^N\to\R$ is the
electric potential and ${\boldsymbol{\mathcal A}}:\R^N\to\R^N$ is a
magnetic potential associated to the magnetic field $B=\mathop{\rm
  curl} {\boldsymbol{\mathcal A}}$. For $N=2,3$, ``$\mathop{\rm
  curl}$'' denotes the usual
$\mathop{\rm curl}$ operator, whereas for $N>3$ by $B=\mathop{\rm curl}
{\boldsymbol{\mathcal A}}$ we mean the 2-form $(B_{jk})$ with
$B_{jk}:=\partial_j{\mathcal A}_k-\partial_k{\mathcal A}_j$, where
${\boldsymbol{\mathcal A}}=({\mathcal A}_j)_{j=1,\dots,N}$.
Linear and nonlinear elliptic equations associated to electromagnetic
hamiltonians have been the object of a wide recent mathematical
research; we quote, among others,
\cite{as,cs,cingolani,cingolanisecchi,el,Kurata}.

In this paper we are concerned with singular homogeneous
electromagnetic potentials $({\boldsymbol{\mathcal A}},V)$ which make
the operator invariant by scaling, namely of the form
\begin{equation*}
{\boldsymbol{\mathcal A}}(x)=\frac{\A\big(\frac{x}{|x|}\big)}
    {|x|}\quad\text{and}\quad
V(x)=-
\dfrac{a\big(\frac{x}{|x|}\big)}{|x|^2}
\end{equation*}
in $\R^N$, where $N\geq 2$, $\A\in C^1({\mathbb S}^{N-1},\R^N)$, and
$a\in L^{\infty}({\mathbb S}^{N-1},\R)$.
A prototype in dimension $2$ is given by potentials associated to thin
solenoids: if the radius of the solenoid tends to zero while the flux
through it remains constant, then the particle is subject to a
$\delta$-type magnetic field, which is called \emph{Aharonov-Bohm}
field. A vector potential associated to the Aharonov-Bohm
magnetic field in $\R^2$ has the form
\begin{equation}\label{eq:31}
{\boldsymbol{\mathcal
    A}}(x_1,x_2)=\alpha\bigg(-\frac{x_2}{|x|^2},\frac{x_1}{|x|^2}\bigg),\quad
(x_1,x_2)\in\R^2,
\end{equation}
with $\alpha\in\R$ representing the circulation of
${\boldsymbol{\mathcal A}}$ around the solenoid. We notice that the
potential in~(\ref{eq:31}) is singular at~$0$, homogeneous of degree
$-1$ and satisfies the following transversality condition
\begin{equation*}
\A(\theta)\cdot\theta=0\quad\text{for all }\theta\in{\mathbb S}^{N-1}.
\end{equation*}
We refer to \cite{balinsky,hhlt,mor} for properties of Aharonov-Bohm
magnetic potentials and related Hardy inequalities. In the present
paper, we consider, for $N\geq 2$, a larger class of singular vector
potentials, characterized by the presence of a homogeneous isolated
singularity of order $-1$ and by the transversality (or Poincar\'e)
condition (we address the reader to \cite{jackson} and \cite[\S
8.4.2]{thaller} for details about the transversal or Poincar\'e
gauge). Such a class includes, for $N=2$, the Aharonov-Bohm magnetic
potential \eqref{eq:31}. The Aharonov-Bohm potential in dimension
$N=3$ is singular on a straight line and is not covered by the
analysis performed here, which only allows treating isolated
singularities. In a forthcoming paper, we will extend the present
results to potentials with cylindrical singularity including the
$3$-dimensional Aharonov-Bohm case.

Singular homogeneous electric potentials which scale as the laplacian
arise in nonrelativistic molecular physics, where the interaction between an
electric charge and the dipole moment ${\mathbf D}\in\R^N$ of a
molecule is described by an inverse square potential with an
anisotropic coupling strength of the form
\[
V(x)=-\dfrac{\lambda\,(x\cdot{\mathbf d})}{|x|^3}
\]
in $\R^N$, where $\lambda>0$ is proportional to the magnitude of the
dipole moment ${\mathbf D}$ and ${\mathbf d}={\mathbf D}/{|{\mathbf
    D}|}$ denotes the orientation of ${\mathbf D}$, see
\cite{FMT2,FMT3,leblond}.  We notice that the above electric potential
is singular at~$0$ and homogeneous of degree $-2$.

We aim to describe the asymptotic behavior near the singularity of
solutions to equations associated to the following class of
Schr\"odinger operators with singular homogeneous electromagnetic
potentials:
\begin{equation*}
  {\mathcal L}_{\A,a}:=
  \left(-i\,\nabla+\frac{\A\big(\frac{x}{|x|}\big)}
    {|x|}\right)^{\!\!2}-\dfrac{a\big(\frac{x}{|x|}\big)}{|x|^2}.
\end{equation*}
We study both linear and nonlinear equations obtained as perturbations
of the operator ${\mathcal L}_{\A,a}$ in a domain $\Omega\subset \R^N$
containing either the origin or a neighborhood of $\infty$. More
precisely, we deal with linear equations of the type
\begin{equation} \label{u}
\mathcal L_{\A,a} u=h(x)\, u, \qquad {\rm in} \ \Omega,
\end{equation}
where $h\in L^\infty_{{\rm loc}}(\Omega\setminus \{0\})$ is negligible
with respect to the inverse square potential $|x|^{-2}$ near the
singularity, and semilinear equations
\begin{equation} \label{nonlin0}
\mathcal L_{\A,a}u(x)=f(x,u(x))
\end{equation}
with $f$ having at most critical growth.

Regularity properties of solutions to Schr\"odinger equations with
less singular magnetic and electric potentials have been studied by
several authors. In particular, in \cite{cs}, boundedness and decay at
$\infty$ of solutions are proved in dimensions $N\geq 3$ for $L^2_{\rm
  loc}$ magnetic potentials and electric potentials with $L^{N/2}$
negative part. It is also worth quoting \cite{Kurata2} and
\cite{Kurata}, where, in dimensions $N\geq 3$, local boundedness and,
respectively, a unique continuation property are established under the
assumption that the electric potential and the square of the magnetic
one belong to the Kato class. In \cite{Kurata2} continuity of
solutions is also obtained under restricted assumptions on the
potentials.

Due to the presence of a
more strong singularity which keeps  potentials in $  {\mathcal
  L}_{\A,a}$ out of the Kato
class, it is natural to expect that solutions to equations (\ref{u})
and (\ref{nonlin0}) behave singularly at the origin: our purpose is to
describe the rate and the shape of the singularity of
solutions, by relating them to the eigenvalues and the eigenfunctions of
a Schr\"odinger operator on the sphere ${\mathbb S}^{N-1}$
 corresponding to the angular part of ${\mathcal L}_{\A,a}$.

 As remarked in \cite{FMT1, FMT3} for the case $\A=0$ (i.e. no
 magnetic vector potential), the estimate of the behavior of solutions
 to elliptic equations with singular potentials near the singularities
 has several important applications to the study of spectral
 properties of the associated Schr\"odinger operator, such as
 essential self-adjointness, positivity, etc. In \cite{FMT2}, the
 exact asymptotic behavior near the singularity of solutions to
 Schr\"odinger equations with singular dipole-type electric potentials
 is established, using separation of variables combined with a
 comparison method.  Comparison and maximum principles play a crucial
 role also in \cite{pinchover94}, where the existence of the limit at
 the singularity of any quotient of two positive solutions to Fuchsian
 type elliptic equations is proved.  In the presence of a singular
 magnetic potential, comparison methods are no more available,
 preventing us from a direct extension of the results of
 \cite{FMT2,pinchover94}. This difficulty is overcome by a Almgren
 type monotonicity formula (see \cite{almgren,GL}) and blow-up methods
 which allow avoiding the use of comparison methods.

\subsection{Assumptions and functional setting}\label{assumptions}

As already mentioned, we shall deal with electromagnetic potentials
$({\boldsymbol{\mathcal A}},V)$ in $\R^N$, $N\geq 2$, satisfying the
following assumptions:
\begin{align}
  &{\boldsymbol{\mathcal
      A}}(x)=\dfrac{\A\big(\frac{x}{|x|}\big)}
  {|x|}\quad\text{and}\quad V(x)=-
  \dfrac{a\big(\frac{x}{|x|}\big)}{|x|^2}\tag{\bf A.1}\label{eq:36}
  &\text{\bf\hfil(homogeneity)}\\[5pt]
  &\A\in C^1({\mathbb S}^{N-1},\R^N)\quad\text{and}\quad
  a\in L^{\infty}({\mathbb S}^{N-1},\R)\tag{\bf A.2}\label{eq:reg}
  &\text{\bf (regularity of angular coefficients)}\\[5pt]
&\A(\theta)\cdot\theta=0\quad\text{for all }\theta\in{\mathbb S}^{N-1}.
\tag{\bf A.3}\label{eq:transversality}
  &\text{\bf (transversality)}
\end{align}
Under the transversality assumption (\ref{eq:transversality}),
the operator ${\mathcal L}_{\A,a}$ acts on functions $u:\R^N\to\C$ as
\begin{displaymath}
  {\mathcal L}_{\A,a}u= -\Delta
  u-\frac{a\big(\frac{x}{|x|}\big)-|\A\big(\frac{x}{|x|}\big)|^2
    +i\dive\nolimits_{{\mathbb S}^{N-1}}\A\big(\frac{x}{|x|}\big)
  }{|x|^2}\,u-2i\,\frac{\A\big(\frac{x}{|x|}\big)}{|x|}\cdot\nabla u,
\end{displaymath}
 where $\dive\nolimits_{{\mathbb S}^{N-1}}\A$ denotes the Riemannian
 divergence of $\A$ on the unit sphere ${\mathbb S}^{N-1}$ endowed
 with the standard metric.

 The positivity properties of the Schr\"odinger operator ${\mathcal
   L}_{\A,a}$ are strongly related to the first eigenvalue of the
 angular component of the operator on the sphere $\mathbb S^{N-1}$.
 More precisely, the positivity  of the quadratic form associated to ${\mathcal
   L}_{\A,a}$ is ensured under the assumption
\begin{align} \label{positivity}
\quad\!\!\mu_1(\A,a)>-\bigg(\frac{N-2}2\bigg)^{\!\!2}\tag{\bf A.4},
\hskip5.5cm\text{\bf (positive definiteness)},
\end{align}
see Lemma \ref{l:pos}, where $\mu_1(\A,a)$ is the first eigenvalue of
the angular component of the operator on the sphere $\mathbb S^{N-1}$,
i.e. of the operator
\[
L_{\A,a}:=\big(-i\,\nabla_{\mathbb S^{N-1}}+\A\big)^2-a.
\]
When dealing with the nonlinear problem \eqref{nonlin0} we introduce the stronger
condition

\begin{align} \label{pos-mu10a}
\quad\!\!\mu_1(0,a)>-\bigg(\frac{N-2}2\bigg)^{\!\!2}\tag{\bf A.5}.
\hskip5.5cm\text{}
\end{align}

From the diamagnetic inequality it follows that $\mu_1(0,a)\leq \mu_1(\A,a)$
with equality holding if and only if curl$\frac{\A}{|x|}=0$ in the sense
of distributions, see Lemma \ref{irrotational}. In particular the assumption
\eqref{pos-mu10a} is in general stronger than \eqref{positivity}.

The spectrum of the angular operator  $L_{\A,a}$ is
discrete and consists in  a nondecreasing sequence of  eigenvalues
\[
\mu_1(\A,a)\leq\mu_2(\A,a)\leq\cdots\leq\mu_k(\A,a)\leq\cdots
\]
diverging to $+\infty$, see Lemma \ref{l:spe} in the Appendix.
Condition (\ref{positivity}) is fundamental to introduce a proper
functional setting in which to frame our analysis. Let us define
$\Distar$ as the completion of $C^\infty_{\rm
  c}(\R^N\setminus\{0\},\C)$ with respect to the norm
\begin{equation}\label{eq:47}
  \|u\|_{\Distar}:=\bigg(\int_{\R^N}\bigg(\big|\nabla u(x)\big|^2+
\frac{|u(x)|^2}{|x|^2}\bigg)
\,dx\bigg)^{\!\!1/2}.
\end{equation}
It is easy to verify that
\[
\Distar=\left\{u\in L^1_{\rm
    loc}(\R^N\setminus\{0\},\C): \frac{u}{|x|}\in L^2(\R^N,\C)\text{
    and } \nabla u\in L^2(\R^N,\C^N)\right\}.
\]
The following lemma ensures that, under assumption \eqref{positivity},
the space $\Distar$ coincides with the Hilbert space originated by the
quadratic form $Q_{\A,a}$ associated to the operator ${\mathcal L}_{\A,a}$
\begin{equation}\label{eq:48}
Q_{\A,a}:\Distar\to\R,\quad
Q_{\A,a}(u):=\int_{\R^N} \bigg[ \bigg|\bigg(\nabla+i\,\frac
{\A\big({x}/{|x|}\big)}{|x|}\bigg)u(x)\bigg|^2
-\frac{a\big({x}/{|x|}\big)}{|x|^2}|u(x)|^2\bigg]\,dx.
\end{equation}

\begin{Lemma}\label{l:posN}
  Assume that $N\geq 2$ and \eqref{eq:reg}, \eqref{eq:transversality},
  \eqref{positivity} hold. Then
\begin{align*}
{\rm i)}\quad& \inf_{u\in \Distar\setminus\{0\}}\frac{Q_{\A,a}(u)}
{\int_{\R^N}|x|^{-2}|u(x)|^2\,dx}>0\\[5pt]
{\rm ii)}\quad& \text{$Q_{\A,a}$ is positive definite in $\Distar$, i.e. }
\inf_{u\in \Distar\setminus\{0\}}\frac{Q_{\A,a}(u)}
{\|u\|^2_{\Distar}}>0 \\[5pt]
{\rm iii)}\quad& \text{$\Distar=\DiAa$, where $\DiAa$ is the
  completion of $C^\infty_{\rm c}(\R^N\setminus\{0\},\C)$} \\
  \quad & \text{with respect to the norm} \\[-7pt]
\noalign{$$\|u\|_{\DiAa}:=\big(Q_{\A,a}(u)\big)^{1/2}.$$}\\[-12pt]
& \text{Moreover the norms $\|\cdot\|_{\Distar}$
  and $\|\cdot\|_{\DiAa}$ are equivalent.}
\end{align*}
\end{Lemma}
In any open bounded domain $\Omega\subset \R^N$ containing $0$,
we introduce the functional space $H^1_{*}(\Omega,\C)$ as the completion of
\[
\{u\in
H^1(\Omega,\C)\cap C^\infty(\Omega,\C):u\text{ vanishes in a
  neighborhood of }0\}
\]
 with respect to the norm
 $$
 \|u\|_{H^1_{*}(\Omega,\C)}=\left(\left\|\nabla u\right\|^2
   _{L^2(\Omega,\C^N)} +\|u\|^2_{L^2(\Omega,\C)}+\Big\|\frac{u}{|x|}
   \Big\|^2_{L^2(\Omega,\C)}\right)^{\!\!1/2}.
$$
It is easy to verify that
\[
H^1_{*}(\Omega,\C)=\left\{u\in H^1(\Omega,\C):\frac{u}{|x|}\in
  L^2(\Omega,\C)\right\}.
\]
If $N\geq 3$, $H^1_{*}(\Omega,\C)=H^1(\Omega,\C)$ and their norms are
equivalent, as one can easily deduce from the Hardy type inequality
with boundary terms due to \cite{wz} (see (\ref{eq:32})) and
continuity of Sobolev trace imbeddings. On
the other hand, if $N=2$,
$H^1_{*}(\Omega,\C)$ is strictly smaller than $H^1(\Omega,\C)$.

For any $h$ satisfying
\begin{equation} \label{hph}
h\in L^{\infty}_{\rm loc}(\Omega\setminus\{0\},\C),
\qquad |h(x)|=O(|x|^{-2+\e})  \quad \mbox{as } |x|\rightarrow 0
\quad\text{for some $\e>0$,}
\end{equation}
we introduce the notion of \emph{weak solution} to (\ref{u}): we say
that a function $u\in H^1_{*}(\Omega,\C)$ is a $H^1_{*}(\Omega,\C)$-weak
solution to (\ref{u}) if, for all $w\in H^1_0(\Omega,\C)$ such that
$\frac w{|x|}\in L^2(\Omega,\C)$,
\begin{equation*}
{\mathcal Q}_{\A,a}^{\Omega}(u,w)
=\int_{\Omega} h(x)
  u(x)\overline{w(x)}\,dx,
\end{equation*}
where ${\mathcal Q}_{\A,a}^{\Omega}:H^1_{*}(\Omega,\C)\times
H^1_{*}(\Omega,\C)\to\C$ is defined by
\[ {\mathcal Q}_{\A,a}^{\Omega}(u,w)\!:=\!  \int_{\Omega}\!
\bigg(\!\nabla u(x)+i\,\frac
{\A\big(\frac{x}{|x|}\big)}{|x|}u(x)\!\bigg)\!\cdot\!
\bigg(\!\overline{\nabla w(x)+i\,\frac
  {\A\big(\frac{x}{|x|}\big)}{|x|}w(x)}\!\bigg)\,dx\!
-\!\int_{\Omega} \!\frac{a(x/|x|)}{|x|^2}\,u(x)\overline{w(x)}\,dx.
\]
In an analogous way, we define the notion of weak solutions to
(\ref{nonlin0}) in a bounded domain for  every
Carath\'eodory function $f:\Omega\times \C\rightarrow \C$ satisfying
the growth restriction
\begin{equation}\label{subcrit_0}
  \left| \frac{f(x,z)}{z} \right|
  \leq
\begin{cases}
C_f(1+|z|^{2^*-2}),  &\text{ if $N\geq 3$},\\
C_f(1+|z|^{p-2})  \quad\text{for some $p>2$ },&\text{ if $N=2$ },
\end{cases}
\end{equation}
for a.e. $x\in \Omega$ and for all $z\in \C\setminus \{0\}$, where
$2^*=\frac{2N}{N-2}$ is the critical Sobolev exponent and the constant
$C_f>0$ is independent of $x\in \Omega$ and $z\in \C\setminus \{0\}$:
we say that a function $u\in
H^1_{*}(\Omega,\C)$ is a $H^1_{*}(\Omega,\C)$-weak solution to (\ref{nonlin0})
if, for all $w\in H^1_0(\Omega,\C)$ such that
$\frac w{|x|}\in L^2(\Omega,\C)$,
\begin{equation*}
{\mathcal Q}_{\A,a}^{\Omega}(u,w)
=\int_{\Omega} f(x,u(x))\overline{w(x)}\,dx.
\end{equation*}
Regularity of solutions either to (\ref{u}) or to (\ref{nonlin0})
outside the singularity follows from classical elliptic regularity
theory, as described in the following remark.

\begin{remark} \label{remA0} If $\A\in C^1(\SN,\R^N)$, $a\in
  L^\infty(\SN,\R)$, and $h\in L^\infty_{{\rm loc}}(\Omega\setminus
  \{0\})$, then, from standard regularity theory and bootstrap
  arguments, it follows that any $H^1_{*}(\Omega,\C)$-weak solution $u$ of
  (\ref{u}) satisfies $u\in W^{2,p}_{{\rm loc}}(\Omega\setminus
  \{0\})$ for any $1\leq p<\infty$ and in particular $u\in
  C^{1,\tau}_{{\rm loc}}(\Omega\setminus\{0\},\C)$ for any $\tau\in
  (0,1)$. The Brezis-Kato technique introduced in \cite{BrezisKato},
  standard regularity theory, and bootstrap arguments, lead to the
  same conclusion also for $H^1_{*}(\Omega,\C)$-weak solutions to
  (\ref{nonlin0}) with $f$ as in (\ref{subcrit_0}).
\end{remark}

\subsection{Statement of the main results}\label{sec:stat-main-results}

The following theorem provides a classification of the behavior of
any solution $u$ to (\ref{u}) near the singularity based on the limit as
$r\to 0^+$ of the \emph{Almgren's frequency} function (see \cite{GL})
\begin{equation}\label{eq:37}
  {\mathcal N}_{u,h}(r)=\frac{r\int_{B_r} \big[\big|\nabla u(x)+i
    \frac{A(x/|x|)}{|x|}u(x)\big|^2-\frac{a(x/|x|)}{|x|^2}|u(x)|^2 -(\Re
    h(x))|u(x)|^2 \big] \, dx}{\int_{\partial B_r}|u(x)|^2 \, dS},
  \end{equation}
where, for any $r>0$, $B_r$ denotes the ball $\{x\in\R^N: |x|<r\}$.

\begin{Theorem} \label{t:asym} Let $\Omega\subset\R^N$, $N\geq 2$, be
  a bounded open set containing $0$,
\eqref{eq:36}, \eqref{eq:reg}, \eqref{eq:transversality}, \eqref{positivity}
hold,
and $u$ be a weak $H^1_{*}(\Omega,\C)$-solution to (\ref{u}),
  $u\not\equiv 0$, with $h$ satisfying (\ref{hph}).  Then, letting
  ${\mathcal N}_{u,h}(r)$ as in (\ref{eq:37}), there exists $k_0\in
  \N$, $k_0\geq 1$, such that
  \begin{align}\label{eq:20}
    \lim_{r\to 0^+}{\mathcal
      N}_{u,h}(r)=-\frac{N-2}2+\sqrt{\bigg(\frac{N-2}{2}
      \bigg)^{\!\!2}+\mu_{k_0}(\A,a)}.
  \end{align}
  Furthermore, if $\gamma$ denotes the limit in (\ref{eq:20}), $m\geq
  1$ is the multiplicity of the eigenvalue $\mu_{k_0}(\A,a)$, and
  $\{\psi_i:\,j_0\leq i\leq j_0+m-1\}$ ($j_0\leq k_0\leq j_0+m-1$) is
  an $L^2({\mathbb S}^{N-1},\C)$-orthonormal basis for the eigenspace
  of the operator $L_{\A,a}$ associated to $\mu_{k_0}(\A,a)$, then
\begin{equation} \label{estu}
\lambda^{-\gamma} u(\lambda\theta)\longrightarrow\sum_{i=j_0}^{j_0+m-1}
  \beta_i\psi_{i}(\theta)\quad \text{in }C^{1,\tau}({\mathbb S}^{N-1},\C)
\quad\text{as }\lambda\to 0^+,
\end{equation}
 and
\begin{equation} \label{estgradu}
\lambda^{1-\gamma}\nabla u(\lambda\theta)\longrightarrow
\sum_{i=j_0}^{j_0+m-1} \beta_i
\big(\gamma \psi_{i}(\theta)\theta+\nabla_{{\mathbb S}^{N-1}}
\psi_{i}(\theta)\big)\quad \text{in }C^{0,\tau}({\mathbb S}^{N-1},\C^N)
\quad\text{as }\lambda\to 0^+,
\end{equation}
for any $\tau\in(0,1)$, where
\begin{align}\label{eq:38}
  \beta_i= \int_{{\mathbb S}^{N-1}}\bigg[ R^{-\gamma}u(R\theta)+
  \int_{0}^R\frac{ h(s\,\theta)u(s\,\theta)}{2\gamma+N-2}
  \bigg(s^{1-\gamma}-\frac{s^{\gamma+N-1}}{R^{2\gamma+N-2}}\bigg)ds
  \bigg]\overline{\psi_{i}(\theta)}\,dS(\theta),
 \end{align}
 for all $R>0$ such that $\overline{B_{R}}=
\{x\in\R^N:|x|\leq R\}\subset\Omega$
and $(\beta_{j_0},\beta_{j_0+1},\dots,\beta_{j_0+m-1})\neq(0,0,\dots,0)$.
\end{Theorem}
We notice that (\ref{eq:38}) is actually a {\it Cauchy's integral type
  formula} for $u$ which allows retracing the behavior of $u$ at the
singularity from the values of $u$ along any circle centered at $0$,
up to some term depending on the perturbation $h$.

An application of Theorem \ref{t:asym} to the special case of
Aharonov-Bohm magnetic fields in $\R^2$ of the form
(\ref{eq:31}) is described in section \ref{sec:ahar-bohm-magn}.

Theorem \ref{t:asym} implies a {\em{strong unique continuation
  property}} as the following corollary states. Moreover, if $\gamma>0$
(as e.g. it happens under assumption (\ref{positivity}) in dimension
$N=2$) then the solutions to \eqref{u} are H\"older continuous for
$0<\gamma<1$ and Lipschitz continuous for $\gamma\geq1$.

\begin{Corollary}\label{cor:holdercont}
  Suppose that all the assumptions of Theorem \ref{t:asym} hold true.
  Let $\gamma$ denote the limit in (\ref{eq:20}) and  $u$ be a weak
  $H^1_{*}(\Omega,\C)$-solution to (\ref{u}).
\begin{itemize}
\item[\rm (i)] If $u(x)=O(|x|^k)$ as $|x|\to 0$ for all $k\in
  \N$, then $u\equiv 0$ in $\Omega$.
\item[\rm (ii)] If $0<\gamma<1$ then
  $u\in C^{0,\gamma}_{\rm loc}(\Omega,\C)$.
\item[\rm (iii)] If $\gamma\geq 1$ then $u$ is
  locally Lipschitz continuous in $\Omega$.
\end{itemize}
\end{Corollary}

We notice that the unique continuation property proved in
\cite{Kurata} for electromagnetic potentials in the Kato class does
not contain the result stated in part (i) of Corollary
\ref{cor:holdercont} for singular homogeneous magnetic potentials. We
also remark that the monotonicity argument used to prove Theorem
\ref{t:asym} (see sections \ref{sec:monot-prop} and
\ref{sec:proofs-theor-reft}) actually applies when perturbing the
magnetic homogeneous potential with a non singular term, namely with a
magnetic potential of the form
\begin{equation}\label{eq:44}
 {\boldsymbol{\mathcal
      A}}(x)=\frac{\A\big(\frac{x}{|x|}\big)} {|x|}+{\mathbf b}(x)
\end{equation}
where ${\mathbf b}\in C^1(\Omega\setminus\{0\},\C^N)$ satisfies
$|{\mathbf b}(x)|=O(|x|^{-1+\e})$ and $|\nabla {\mathbf
  b}(x)|=O(|x|^{-2+\e})$ as $|x|\to 0$ for some $\e>0$ as $|x|\to 0$.
For the sake of simplicity, we omit the details of case (\ref{eq:44}),
which can be treated following closely the strategy developed in
sections \ref{sec:monot-prop} and \ref{sec:proofs-theor-reft}.

Due to the homogeneity of the potentials, Schr\"odinger operators
${\mathcal L}_{\A,a}$ are invariant by the Kelvin transform,
\begin{equation*}
\tilde u(x)=|x|^{-(N-2)}u\bigg(\frac{x}{|x|^2}\bigg),
\end{equation*}
which is an isomorphism of $\Distar$. Indeed, if $u\in H^1_{*}(\Omega,\C)$
weakly solves (\ref{u}) in a bounded open set $\Omega$ containing $0$,
then its Kelvin's transform $\tilde u$ weakly solves (\ref{u}) with
$h$ replaced by $|x|^{-4}h(\frac x{|x|^2})$ in the external domain
$\widetilde \Omega=\big\{x\in\R^N: x/{|x|^2}\in\Omega\big\}$.
Weak solution $u$ of problem (\ref{u}) with $h$ satisfying
\begin{equation} \label{hph_00}
h\in L^{\infty}_{\rm loc}(\Omega,\C),
\qquad h(x)=O(|x|^{-2-\e})  \quad \mbox{as } |x|\rightarrow +\infty
\quad\text{for some $\e>0$,}
\end{equation}
in an external domain $\Omega$ (i.e. a domain  $\Omega$  such that
$\R^N\setminus B_{R_0}\subset\Omega\subset\R^N\setminus B_{R_1}$
 for some $R_0>R_1>0$), we mean a function $u$ such
that $\frac{u}{|x|}\in L^{2}(\Omega,\C)$, $\nabla u\in L^2(\Omega,\C^N)$, and
\begin{equation*}
{\mathcal Q}_{\A,a}^{\Omega}(u,w)
=\int_{\Omega} h(x)
  u(x)\overline{w(x)}\,dx,
\end{equation*}
for any $w\in \mathcal D^{1,2}_{*}(\Omega,\C)$, where
$\mathcal D^{1,2}_{*}(\Omega,\C)$ is the completion of $C^\infty_{\rm
  c}(\Omega,\C)$ with respect to the norm $\|u\|_{\mathcal
  D^{1,2}_{*}(\Omega)}:= \big(\big\|\nabla u\big\|^2 _{L^2(\Omega,\C^N)}
  +\big\|\frac{u}{|x|}\big\|^2_{L^2(\Omega,\C)}\big)^{1/2}$.

Theorem \ref{t:asym} and invariance by the Kelvin
transform provide the following description of the behavior of
solutions to (\ref{u}) as $|x|\to\infty$.
The Almgren's frequency type function in exterior domains has the form
\begin{equation}\label{eq:37ext}
  \widetilde{\mathcal N}_{u,h}(r)=
\frac{r\int_{\R^N\setminus B_{r}}
     \big[\big|\nabla u(x)+i
     \frac{A(x/|x|)}{|x|}u(x)\big|^2-\frac{a(x/|x|)}{|x|^2}|u(x)|^2 -(\Re
     h(x))|u(x)|^2 \big] \, dx}{\int_{\partial B_{r}}|u(x)|^2 \, dS}.
 \end{equation}

 \begin{Theorem} \label{asy_infinity} Let $\Omega\subset\R^N$, $N\geq
   2$, be an open set such that $\R^N\setminus
   B_{R_0}\subset\Omega\subset\R^N\setminus B_{R_1}$ for some
   $R_0>R_1>0$, \eqref{eq:36}, \eqref{eq:reg},
   \eqref{eq:transversality}, \eqref{positivity} hold, and $u$ be a
   weak solution to (\ref{u}), $u\not\equiv 0$, with $h$ satisfying
   (\ref{hph_00}). Then, letting $\widetilde{\mathcal N}_{u,h}$ as in
   (\ref{eq:37ext}), there exists $k_0\in \N$, $k_0\geq 1$, such that
 \begin{align}\label{eq:20_00}
   \lim_{r\to +\infty} \widetilde{\mathcal N}_{u,h}(r)=\frac{N-2}2+
   \sqrt{\bigg(\frac{N-2}{2}\bigg)^{\!\!2}+\mu_{k_0}(\A,a)}.
 \end{align}
 Moreover,if $\tilde\gamma$ denotes the limit in
 (\ref{eq:20_00}), $m\geq 1$ is the multiplicity of the eigenvalue
 $\mu_{k_0}(\A,a)$, and $\{\psi_i:\,j_0\leq i\leq j_0+m-1\}$ ($j_0\leq
 k_0\leq j_0+m-1$) is an $L^2({\mathbb S}^{N-1},\C)$-orthonormal basis
 for the eigenspace of the operator $L_{\A,a}$ associated to
 $\mu_{k_0}(\A,a)$, then
\[
 \lambda^{\tilde \gamma}u(\lambda\theta)\longrightarrow
 \sum_{i=j_0}^{j_0+m-1} \widetilde\beta_i\psi_{i}(\theta)\quad\text{in
 }C^{1,\tau}({\mathbb S}^{N-1},\C)
\quad \text{as }\lambda\to+\infty
\]
and
\[
 \lambda^{\tilde \gamma+1}\nabla u(\lambda\theta)\longrightarrow
 \sum_{i=j_0}^{j_0+m-1} \widetilde\beta_i\big(
-\tilde \gamma
\psi_{i}(\theta)\theta+\nabla_{{\mathbb S}^{N-1}}\psi_i(\theta)\big)
\quad\text{in
 }C^{0,\tau}({\mathbb S}^{N-1},\C^N)\quad \text{as }\lambda\to+\infty
\]
for every $\tau\in(0,1)$, where
\begin{align*}
  \widetilde \beta_i= \int_{{\mathbb S}^{N-1}}\bigg[
  R^{\tilde\gamma}u(R\theta)
  +\int_{R}^{+\infty}\frac{h(s\,\theta)u(s\,\theta)}{2\tilde\gamma-N+2}
  \bigg(s^{\tilde\gamma+1}-R^{2\tilde\gamma-N+2}s^{-\tilde\gamma+N-1}\bigg)
  ds \bigg]\overline{\psi_{i}(\theta)}\,dS(\theta),
\end{align*}
for all $R>0$ such that $\R^N\setminus B_R\subset\Omega$ and
$(\widetilde\beta_{j_0},\widetilde\beta_{j_0+1},\dots,
\widetilde\beta_{j_0+m-1})\neq(0,0,\dots,0)$.
\end{Theorem}

A Brezis-Kato type iteration, see \cite{BrezisKato}, allows us to
obtain asymptotics of solutions also for semilinear problems with at
most critical growth. In order to
start such an iterative procedure, we require assumption
(\ref{pos-mu10a}) which allows transforming equation
(\ref{nonlin0}) into a degenerate elliptic equation without singular
potentials on which the Brezis-Kato method applies successfully, see
Lemmas \ref{l:bk1} and \ref{l:bk12d}. The iteration scheme developed in
sections~\ref{sec:bk} and \ref{sec:brezis-kato-type} provides an upper
bound for solutions and then reduces the semilinear problem to a
linear one with enough control on the perturbing potential at the
singularity to apply Theorem \ref{t:asym} and to recover the exact
asymptotic behavior, as stated in the following theorem.

\begin{Theorem} \label{t:asym-nonlin} Let $\Omega\subset\R^N$, $N\geq
  2$, be a bounded open set containing $0$, \eqref{eq:36},
  \eqref{eq:reg}, \eqref{eq:transversality}, \eqref{pos-mu10a} hold,
  and $u$ be a weak $H^1_{*}(\Omega,\C)$-solution to (\ref{nonlin0}),
  $u\not\equiv 0$, with $f$ being a Carath\'eodory function satisfying
  (\ref{subcrit_0}).  Then, there exists $k_0\in \N$, $k_0\geq 1$,
  such that
  \begin{align}\label{eq:20_bis}
    \lim_{r\to 0^+} {\mathcal N}_{u,{f(\cdot,u)}/{u}}(r)
    =-\frac{N-2}2+\sqrt{\bigg(\frac{N-2}{2}\bigg)^{\!\!2}+\mu_{k_0}(\A,a)}.
  \end{align}
  Furthermore, if $\gamma$ denotes the limit in
  (\ref{eq:20_bis}), $m\geq 1$ is the multiplicity of the eigenvalue
  $\mu_{k_0}(\A,a)$, and $\{\psi_i:\,j_0\leq i\leq j_0+m-1\}$ ($j_0\leq
  k_0\leq j_0+m-1$) is an $L^2({\mathbb S}^{N-1},\C)$-orthonormal basis
  for the eigenspace of the operator $L_{\A,a}$ associated to
  $\mu_{k_0}(\A,a)$, then

\[
\lambda^{-\gamma} u(\lambda\theta)\longrightarrow\sum_{i=j_0}^{j_0+m-1}
  \beta_i\psi_{i}(\theta)\quad \text{in }C^{1,\tau}({\mathbb S}^{N-1},\C)
\quad\text{as }\lambda\to 0^+,
\]
 and
\[
\lambda^{1-\gamma}\nabla u(\lambda\theta)\longrightarrow
\sum_{i=j_0}^{j_0+m-1} \beta_i
\big(\gamma \psi_{i}(\theta)\theta+\nabla_{{\mathbb S}^{N-1}}
\psi_{i}(\theta)\big)\quad \text{in }C^{0,\tau}({\mathbb S}^{N-1},\C^N)
\quad\text{as }\lambda\to 0^+,
\]
for any $\tau\in (0,1)$, where
$$
\beta_i=
\int_{{\mathbb S}^{N-1}}\!\bigg[
R^{-\gamma}u(R\theta)+
\int_{0}^R\frac{f(s\,\theta,u(s\,\theta)) }{2\gamma+N-2}
\bigg(s^{1-\gamma}-\frac{s^{\gamma+N-1}}{R^{2\gamma+N-2}}\bigg)
ds \bigg]\overline{\psi_{i}(\theta)}\,dS(\theta)
$$
for all $R>0$ such that $\overline{B_{R}}\subset\Omega$
and $(\beta_{j_0},\beta_{j_0+1},\dots,\beta_{j_0+m-1})\neq(0,0,\dots,0)$.
\end{Theorem}

Similar conclusions as those in Corollary \ref{cor:holdercont} can be
deduced from the above theorem for solutions to semilinear equations of type
(\ref{nonlin0}): under the same
assumption as in Theorem \ref{t:asym-nonlin}, if $\gamma>0$ then the
solutions to (\ref{nonlin0}) are $\gamma$-H\"older continuous for
$0<\gamma<1$ and Lipschitz continuous for $\gamma\geq1$.

 The following result is the counterpart of Theorem \ref{t:asym-nonlin}
in exterior domains.
\begin{Theorem} \label{asy_infinitynonlin} Let $\Omega\subset\R^N$,
  $N\geq 2$, be an open set such that $\R^N\setminus
  B_{R_0}\subset\Omega\subset\R^N\setminus B_{R_1}$ for some
  $R_0>R_1>0$, \eqref{eq:36},
  \eqref{eq:reg}, \eqref{eq:transversality}, \eqref{pos-mu10a}  hold,
 and $u$ be a weak
  solution to (\ref{nonlin0}) in $\Omega$, $u\not\equiv 0$, with $f$
  satisfying, for some $\widetilde C_f>0$,
\begin{equation*}
  \left| \frac{f(x,z)}{z} \right|
  \leq
\begin{cases}
  \widetilde C_f(|x|^{-4}+|z|^{2^*-2}),  &\text{ if $N\geq 3$},\\
  \widetilde C_f|x|^{-4}(1+|z|^{p-2}) \quad\text{for some $p>2$
  },&\text{ if $N=2$},
\end{cases}
\end{equation*}
for a.e. $x\in \Omega$ and  for all $z\in \C\setminus \{0\}$.
Then there exists $k_0\in \N$, $k_0\geq 1$, such that
 \begin{align}\label{eq:39}
   \lim_{r\to +\infty} \widetilde{\mathcal N}_{u,f(\cdot,u)/u}(r)=\frac{N-2}2+
   \sqrt{\bigg(\frac{N-2}{2}\bigg)^{\!\!2}+\mu_{k_0}(\A,a)}.
 \end{align}
 Moreover, if $\tilde\gamma$ denotes the limit in (\ref{eq:39}),
$m\geq 1$ is the multiplicity of the eigenvalue
 $\mu_{k_0}(\A,a)$, and $\{\psi_i:\,j_0\leq i\leq j_0+m-1\}$ ($j_0\leq
 k_0\leq j_0+m-1$) is an $L^2({\mathbb S}^{N-1},\C)$-orthonormal basis
 for the eigenspace of the operator $L_{\A,a}$ associated to
 $\mu_{k_0}(\A,a)$, then
\[
 \lambda^{\tilde \gamma}u(\lambda\theta)\longrightarrow
 \sum_{i=j_0}^{j_0+m-1} \widetilde\beta_i\psi_{i}(\theta)\quad\text{in
 }C^{1,\tau}({\mathbb S}^{N-1},\C)
\quad \text{as }\lambda\to+\infty
\]
and
\[
 \lambda^{\tilde \gamma+1}\nabla u(\lambda\theta)\longrightarrow
 \sum_{i=j_0}^{j_0+m-1} \widetilde\beta_i\big(
-\tilde \gamma
\psi_{i}(\theta)\theta+\nabla_{{\mathbb S}^{N-1}}\psi_i(\theta)\big)
\quad\text{in
 }C^{0,\tau}({\mathbb S}^{N-1},\C^N)\quad \text{as }\lambda\to+\infty
\]
for every $\tau\in(0,1)$, where
$$
\widetilde\beta_i
=
\int_{{\mathbb S}^{N-1}}\bigg[
  R^{\tilde\gamma}u(R\theta)
+\int_{R}^{+\infty}\frac{f(s\,\theta,u(s\,\theta)) }{2\tilde\gamma-N+2}
  \bigg(s^{\tilde\gamma+1}-R^{2\tilde\gamma-N+2}s^{-\tilde\gamma+N-1}\bigg)
  ds \bigg]\overline{\psi_{i}(\theta)}\,dS(\theta)
$$
for all $R>0$ such that $\R^N\setminus B_R\subset\Omega$ and
$(\widetilde\beta_{j_0},\widetilde\beta_{j_0+1},\dots,
\widetilde\beta_{j_0+m-1})\neq(0,0,\dots,0)$.
\end{Theorem}

The paper is organized as follows. In section
\ref{sec:posit-quadr-form} we prove Lemma \ref{l:posN} and discuss the
relation between the positivity of the quadratic form associated to
${\mathcal L}_{\A,a}$ and the first eigenvalue of the angular operator
on the sphere ${\mathbb S}^{N-1}$.  In section \ref{sec:hardyboundary}
we prove a Hardy type inequality with boundary terms and singular
electromagnetic potential, while in section \ref{sec:pohoz-type-ident}
we derive a Pohozaev-type identity for solutions to~(\ref{u}). Section
\ref{sec:monot-prop} contains an Almgren type monotonicity formula,
which is used in section~\ref{sec:proofs-theor-reft} together with a
blow-up method to prove Theorems \ref{t:asym} and \ref{asy_infinity}.
Section \ref{sec:ahar-bohm-magn} contains an application of
Theorem~\ref{t:asym} to Aharonov-Bohm magnetic potentials.  In section
\ref{sec:magn-hardy-sobol} we prove a Hardy-Sobolev inequality with
magnetic potentials which is needed in section \ref{sec:bk} to start a
Brezis-Kato iteration procedure in order to obtain a-priori pointwise
bounds for solutions to the nonlinear equation and to prove Theorems
~\ref{t:asym-nonlin} and \ref{asy_infinitynonlin} in dimension
$N\geq3$. The proof of Theorems ~\ref{t:asym-nonlin} and
\ref{asy_infinitynonlin} in dimension $N=2$ can be found in section
\ref{sec:brezis-kato-type}.  In a final appendix, we recall well-known
results such as the diamagnetic inequality, Hardy's inequality with
boundary terms, and  the description the spectrum of angular operator $L_{\A,a}$.

\medskip
\noindent
{\bf Notation. } We list below some notation used throughout the
paper.\par
\begin{itemize}
\item[-] For all $r>0$, $B_r$ denotes the ball $\{x\in\R^N: |x|<r\}$
  in $\R^N$ with center at $0$ and radius $r$.
\item[-] For all $r>0$, $\overline{B_{r}}= \{x\in\R^N:|x|\leq
r\}$ denotes the closure of $B_r$.
\item[-] $dS$ denotes the volume element on the spheres $\partial B_r$, $r>0$.
\item[-] For every complex number $z\in\C$, $\Re z$ denotes its real
  part and $\Im z$ its imaginary part.
\item[-] For every complex number $z\in\C$, $\overline{z}$ denotes its
  complex conjugate.
\end{itemize}

\medskip\noindent\textbf{Acknowledgements.} The authors wish to express their
gratitude to the unknown referee for his/her helpful remarks, which
stimulated them to revise and improve the paper, both in the results and in the
exposition.

\section{Positivity of the quadratic form}\label{sec:posit-quadr-form}

In this section, we study the quadratic form associated to the
Schr\"odinger operator ${\mathcal L}_{\A,a}$ and defined in
\eqref{eq:48}. To study the sign of $Q_{\A,a}$, we define the first
eigenvalue of $Q_{\A,a}$ with respect to the Hardy singular weight as
\begin{displaymath}
  \lambda_1(\A,a):=\inf_{u\in\Distar\setminus\{0\}}
  \frac{Q_{\A,a}(u)}{\int_{\R^N}\frac{|u(x)|^2}{|x|^2}dx}
\end{displaymath}
and discuss the relation between $\lambda_1(\A,a)$ and the first
eigenvalue of the angular component of the operator on the sphere
$\mathbb S^{N-1}$, i.e. of the operator
\begin{align*}
  L_{\A,a}=\big(-i\,\nabla_{\mathbb S^{N-1}}+\A\big)^2-a=-\Delta_{\mathbb
    S^{N-1}}-\big(a(\theta)-|\A|^2+i\,\dive\nolimits_{{\mathbb
      S}^{N-1}}\A\big)-2i\,\A\cdot\nabla_{\mathbb S^{N-1}}.
\end{align*}
We notice that, by \eqref{eq:reg}, $\lambda_1(\A,a)$ is well defined
and finite.  Let us introduce the Sobolev space
\begin{equation}\label{eq:1magn}
  H^1_{\A}(\mathbb S^{N-1}):=\big\{\psi\in L^2(\mathbb S^{N-1},\C):\,
  \nabla_{\mathbb S^{N-1}}\psi+i\A(\theta)\psi\in L^2(\mathbb S^{N-1},\C^N)\big\},
\end{equation}
  endowed with the norm
\begin{equation}\label{eq:2magn}
  \|\psi\|_{H^1_{\A}(\mathbb S^{N-1})}:= \bigg( \int_{\mathbb S^{N-1}}
  \Big[\big|\big(\nabla_{\mathbb
    S^{N-1}}+i\A(\theta)\big)\psi(\theta)\big|^2 +
  |\psi(\theta)|^2\Big]\,dS(\theta) \bigg)^{\!\!1/2},
\end{equation}
$dS$ denoting the volume element on the sphere ${\mathbb S}^{N-1}$.
We observe that, if $\A\in C^1({\mathbb S}^{N-1},\R^N)$, then
$H^1_{\A}(\mathbb S^{N-1})$ is equal to the classical Sobolev space
$H^1(\mathbb S^{N-1},\C)$ and its norm is equivalent to the
$H^1(\mathbb S^{N-1},\C)$-norm, see Lemma \ref{l:h1sph} in the appendix.

Under assumption \eqref{eq:reg}, the operator $L_{\A,a}$ on $\mathbb
S^{N-1}$ admits a diverging sequence of real eigenvalues
$\mu_1(\A,a)\leq\mu_2(\A,a)\leq\cdots\leq\mu_k(\A,a)\leq\cdots$ the
first of which can be characterized as
\begin{equation}\label{firsteig}
  \mu_1(\A,a)=\min_{\psi\in H^1_{\A}(\mathbb
    S^{N-1})\setminus\{0\}}\frac{\int_{\mathbb S^{N-1}}
    \big[\big|\big(\nabla_{\mathbb S^{N-1}}+i\A(\theta)\big)\psi(\theta)\big|^2
    -a(\theta)|\psi(\theta)|^2\big]\,dS(\theta)}{\int_{\mathbb S^{N-1}}
    |\psi(\theta)|^2\,dS(\theta)},
\end{equation}
see Lemma \ref{l:spe} in the appendix.
 The relation between $\lambda_1(\A,a)$ and $\mu_1(\A,a)$ is
clarified in the following lemma.
\begin{Lemma}\label{l:la1mu1}
 If $N\geq 2$, \eqref{eq:reg} and \eqref{eq:transversality} hold, then
  \begin{displaymath}
    \lambda_1(\A,a)=\mu_1(\A,a)+\bigg(\frac{N-2}{2}\bigg)^{\!\!2}.
  \end{displaymath}
\end{Lemma}
\begin{pf}
  Let $\psi\in H^1_{\A}(\mathbb S^{N-1})$, $\psi\not\equiv0$, attaining
  $\mu_1(\A,a)$ and let $\varphi\in C^\infty_{\rm c}\big((0,+\infty),\R\big)$
so that  $\widetilde\varphi:\,x\mapsto\varphi(|x|)\in
C^\infty_{\rm c}(\R^N\setminus\{0\},\R)$.
If $u(x)=\varphi(|x|)\psi\big(\frac{x}{|x|}\big)$, there holds
\begin{displaymath}
  \left(\nabla+i\,{\textstyle{\frac{\A(\frac{x}{|x|})} {|x|}}}\right)u(x)
  ={\textstyle{\varphi'(|x|)\psi\big(\frac{x}{|x|}\big)\frac{x}{|x|}+
      \frac1{|x|}\varphi(|x|)\nabla_{\mathbb S^{N-1}}\psi\big(\frac{x}{|x|}\big)
      +\frac{i}{|x|}\A\big(\frac{x}{|x|}\big)\varphi(|x|)
      \psi\big(\frac{x}{|x|}\big)}}
\end{displaymath}
and, by assumption (\ref{eq:transversality}),
\begin{displaymath}
  \left|\left(\nabla+i\,{\textstyle{
          \frac{\A(\frac{x}{|x|})} {|x|}}}\right)u(x)\right|^2
  ={\textstyle{|\varphi'(|x|)|^2\big|\psi\big(\frac{x}{|x|}\big)\big|^2+
      \frac{|\varphi(|x|)|^2}{|x|^2}\Big|
      \nabla_{\mathbb S^{N-1}}\psi\big(\frac{x}{|x|}\big)
      +i\,\A\big(\frac{x}{|x|}\big)\psi\big(\frac{x}{|x|}\big)\Big|^2}}.
\end{displaymath}
Therefore, from the definition of  $\lambda_1(\A,a)$ it follows
\begin{align*}
  \lambda_1(\A,a)&\bigg(\int_0^{+\infty}r^{N-1}\frac{|\varphi(r)|^2}{r^2}\,dr\bigg)
  \bigg(\int_{{\mathbb S}^{N-1}}|\psi(\theta)|^2\,dS(\theta)\bigg)\\
  &\leq
  \int_{\R^N}\left[\left|\left(\nabla+i\,{\textstyle{\frac{\A(\frac{x}{|x|})}
            {|x|}}}\right)u(x)\right|^2-a\big({\textstyle{\frac{x}{|x|}}}
    \big)\frac{|u(x)|^2}{|x|^2} \right]\,dx\\
  &=\bigg(\int_0^{+\infty}r^{N-1}|\varphi'(r)|^2\,dr\bigg)
  \bigg(\int_{{\mathbb S}^{N-1}}|\psi(\theta)|^2\,dS(\theta)\bigg)\\
  &\quad
  +\bigg(\int_0^{+\infty}r^{N-1}\frac{|\varphi(r)|^2}{r^2}\,dr\bigg)
  \bigg(\int_{{\mathbb S}^{N-1}} \big[\big|\big(\nabla_{\mathbb
    S^{N-1}}\psi(\theta)+i\A(\theta)\psi(\theta)\big|^2
  -a(\theta)|\psi(\theta)|^2\big]\,dS(\theta)\bigg)\\
  &=\bigg(\int_{{\mathbb S}^{N-1}}|\psi(\theta)|^2\,dS(\theta)\bigg)
  \left[\int_0^{+\infty}r^{N-1}|\varphi'(r)|^2\,dr+\mu_1(\A,a)
    \int_0^{+\infty}r^{N-1}\frac{|\varphi(r)|^2}{r^2}\,dr\right].
\end{align*}
Hence
\begin{align*}
  \lambda_1(\A,a)-\mu_1(\A,a)\leq\frac{\int_0^{+\infty}r^{N-1}
    |\varphi'(r)|^2\,dr}{\int_0^{+\infty}r^{N-3}|\varphi(r)|^2\,dr}
  =\frac{\int_{\R^N}|\nabla\widetilde\varphi(x)|^2\,dx}
  {\int_{\R^N}\frac{|\widetilde\varphi(x)|^2}{|x|^2}\,dx}
\end{align*}
for every radial function $\widetilde\varphi\in C^\infty_{\rm
  c}(\R^N\setminus\{0\},\R)$.  Hence by Schwarz symmetrization
\begin{align*}
  \lambda_1(\A,a)-\mu_1(\A,a)&\leq\inf_{{\substack{\widetilde\varphi\in
C^\infty_{\rm
  c}(\R^N\setminus\{0\},\R)\setminus\{0\}\\\widetilde\varphi\text{
          radial}}}}
  \frac{\int_{\R^N}|\nabla\widetilde\varphi(x)|^2\,dx}
  {\int_{\R^N}\frac{|\widetilde\varphi(x)|^2}{|x|^2}\,dx}
\\&=\inf_{v\in C^\infty_{\rm
  c}(\R^N\setminus\{0\},\R)\setminus\{0\}}
  \frac{\int_{\R^N}|\nabla v(x)|^2\,dx}
  {\int_{\R^N}\frac{|v(x)|^2}{|x|^2}\,dx}=\bigg(\frac{N-2}2\bigg)^{\!\!2},
\end{align*}
where the last identity is due to the optimality of the classical best
Hardy constant for $N\geq 3$ and to direct calculations for $N=2$.
To prove the reverse inequality, let $u\in C^\infty_{\rm
  c}(\R^N\setminus\{0\},\C)$.  The magnetic gradient of $u$ can be
written in polar coordinates as
\begin{displaymath}
  \nabla u(x)+i\,\frac{\A\big(\frac{x}{|x|}\big)}{|x|}\,u(x)
  =\big(\partial_ru(r,\theta)\big)\theta+
  \frac1r\nabla_{{\mathbb S}^{N-1}}u(r,\theta)
  +i\,\frac{u(r,\theta)}{r}\,\A(\theta), \quad r=|x|,\quad \theta=\frac{x}{|x|}.
\end{displaymath}
  By assumption (\ref{eq:transversality}), there holds
  \begin{equation}\label{eq:46}
     \left|\nabla u(x)+i\,\frac{\A\big(\frac{x}{|x|}\big)}{|x|}\,u(x)\right|^2
=\big|\partial_ru(r,\theta)\big|^2+
\frac1{r^2}\big|\nabla_{{\mathbb S}^{N-1}}u(r,\theta)
+i\,\A(\theta)u(r,\theta)\big|^2,
\end{equation}
hence
\begin{multline}\label{eq:1}
  Q_{\A,a}(u)= \int_{{\mathbb S}^{N-1}}\bigg(
  \int_0^{+\infty}r^{N-1}|\partial_ru(r,\theta)|^2\,dr\bigg)\,dS(\theta)\\
  + \int_0^{+\infty}\frac{r^{N-1}}{r^2}\bigg(\int_{{\mathbb
      S}^{N-1}}\left[|\nabla_{{\mathbb S}^{N-1}}u(r,\theta)
    +i\,\A(\theta)u(r,\theta)|^2-a(\theta)
    |u(r,\theta)|^2\right]\,dS(\theta)\bigg)\,dr.
\end{multline}
For all $\theta\in{\mathbb S}^{N-1}$, let $\varphi_{\theta}\in
C^\infty_{\rm c}((0,+\infty),\C)$ be defined by
$\varphi_{\theta}(r)=u(r,\theta)$, and $\widetilde\varphi_{\theta}\in
C^\infty_{\rm c}(\R^N\setminus\{0\},\C)$ be the radially symmetric
function given by
$\widetilde\varphi_{\theta}(x)=\varphi_{\theta}(|x|)$. If $N\geq 3$, Hardy's
inequality yields
\begin{align}\label{eq:2}
  \int_{{\mathbb S}^{N-1}}&\bigg(
  \int_0^{+\infty}r^{N-1}|\partial_ru(r,\theta)|^2\,dr\bigg)\,dS(\theta)=
  \int_{{\mathbb S}^{N-1}}\bigg(
  \int_0^{+\infty}r^{N-1}|\varphi_{\theta}'(r)|^2\,dr\bigg)\,dS(\theta)\\
  \notag& =\frac1{\omega_{N-1}}  \int_{{\mathbb S}^{N-1}}\bigg(
\int_{\R^N}|\nabla\widetilde \varphi_{\theta}(x)|^2\,dx\bigg)\,dS(\theta)
\\&\notag\geq
\frac1{\omega_{N-1}}  \bigg(\frac{N-2}2\bigg)^{\!\!2}
  \int_{{\mathbb S}^{N-1}}\bigg(\int_{\R^N}\frac{|\widetilde
    \varphi_{\theta}(x)|^2}{|x|^2}\,dx\bigg)\,dS(\theta)\\
\notag&=\bigg(\frac{N-2}2\bigg)^{\!\!2}
  \int_{{\mathbb S}^{N-1}}\bigg(
  \int_0^{+\infty}\frac{r^{N-1}}{r^2}|u(r,\theta)|^2\,dr\bigg)\,dS(\theta)
=\bigg(\frac{N-2}2\bigg)^{\!\!2}\int_{\R^N}\frac{|u(x)|^2}{|x|^2}\,dx,
\end{align}
where $\omega_{N-1}$ denotes the volume of the unit sphere ${\mathbb
  S}^{N-1}$, i.e. $\omega_{N-1}=\int_{{\mathbb S}^{N-1}}dS(\theta)$. For $N=2$
(\ref{eq:2}) trivially holds.
On the other hand, from the definition of $\mu_1(\A,a)$ it follows
that
\begin{equation} \label{eq:3}
  \int_{{\mathbb
      S}^{N-1}}\!\!\left[|\nabla_{{\mathbb S}^{N-1}}u(r,\theta)
    +i\,\A(\theta)u(r,\theta)|^2\!-a(\theta)|u(r,\theta)|^2\right]dS(\theta)
  \geq \mu_1(\A,a)\int_{{\mathbb S}^{N-1}}\!|u(r,\theta)|^2dS(\theta).
\end{equation}
From (\ref{eq:1}), (\ref{eq:2}), and (\ref{eq:3}), we deduce that
\begin{displaymath}
  Q_{\A,a}(u)\geq \left[\bigg(\frac{N-2}2\bigg)^{\!\!2}+\mu_1(\A,a)\right]
  \int_{\R^N}\frac{|u(x)|^2}{|x|^2}\,dx
\quad \text{for all }u\in C^\infty_{\rm
  c}(\R^N\setminus\{0\},\C),
\end{displaymath}
which, by density of $C^\infty_{\rm c}(\R^N\setminus\{0\},\C)$ in
$\Distar$, implies
\begin{displaymath}
  \lambda_1(\A,a)\geq \bigg(\frac{N-2}2\bigg)^{\!\!2}+\mu_1(\A,a),
\end{displaymath}
thus completing the proof.
\end{pf}

The relation between positivity of $Q_{\A,a}$ and the values
$\mu_1(\A,a),\lambda_1(\A,a)$ is described  in the following lemma.
\begin{Lemma}\label{l:pos}
 If $N\geq 2$, \eqref{eq:reg} and \eqref{eq:transversality} hold, then
  the following
  conditions are equivalent:
\begin{align*} {\rm i)}\quad&Q_{\A,a} \text{ is positive definite in $\Distar$,
    i.e. } \inf_{u\in\Distar\setminus\{0\}}\frac{Q_{\A,a}(u)}
  {\|u\|^2_{\Distar}}>0; \\
  {\rm ii)} \quad&\lambda_1(\A,a)>0;\\
  {\rm iii)}\quad&
  \mu_1(\A,a)>-\big({\textstyle{\frac{N-2}2}}\big)^{\!2}.
\end{align*}
\end{Lemma}
\begin{pf}
  The equivalence between (ii) and (iii) is an immediate consequence
  of Lemma \ref{l:la1mu1}.  The fact that (i) implies (ii) follows
  easily from \eqref{eq:47}.  It remains to prove that (ii) implies
  (i). One can proceed as in the proof of \cite[Proposition 1.3]{Ter}.
  For completeness we give here the details.  Assume (ii) and suppose
  by contradiction that (i) is not true.  Then for any $\e>0$ there
  exists $u_\e\in \Distar$ such that
\begin{align*}
  Q_{\A,a}(u_\e)&<\e\|u_\e\|^2_{\Distar} \\
  &\leq 2(\|\A\|_{L^{\infty}({\mathbb
      S}^{N-1},\R^N)}^2+1)\,\e\bigg(\int_{\R^N}\bigg|\bigg(\nabla+i\,\frac
  {\A\big({x}/{|x|}\big)}{|x|}\bigg)u(x)\bigg|^2+
  \int_{\R^N}\frac{|u(x)|^2}{|x|^2}\,dx\bigg)
\end{align*}
and hence, for $\e$ small,
\[
\lambda_1\Big(\A,\frac1{1-2\e(\|\A\|_{L^{\infty}({\mathbb
      S}^{N-1},\R^N)}^2+1)}a\Big)<\frac{2\e(\|\A\|_{L^{\infty}({\mathbb
      S}^{N-1},\R^N)}^2+1)}{1-2\e(\|\A\|_{L^{\infty}({\mathbb
      S}^{N-1},\R^N)}^2+1)}.
\]
  On the other hand, from
  the characterization of $\lambda_1(\A,a)$ given in Lemma
  \ref{l:la1mu1}, we have that the map $a\mapsto \lambda_1(\A,a)$ is
  continuous with respect to the $L^\infty(\SN)$-norm and hence,
  letting $\e\rightarrow 0$, we obtain $\lambda_1(\A,a)\leq 0$, a
  contradiction.
\end{pf}

The previous lemma allows relating $\Distar$ with the Hilbert space
$\DiAa$ generated by the quadratic form $Q_{\A,a}$, thus proving
Lemma \ref{l:posN}.

\medskip\noindent \begin{pfn}{Lemma \ref{l:posN}} i) follows from
  Lemma \ref{l:la1mu1} and assumption \eqref{positivity}.  ii) is a
  direct consequence of Lemma \ref{l:pos} and \eqref{positivity}.
  From ii) we deduce that $(Q_{\A,a}(\cdot))^{1/2}$ defines a norm in
  $C^\infty_{\rm c}(\R^N\setminus\{0\},\C)$ which is equivalent to
  $\|\cdot\|_{\Distar}$. Hence completion of $C^\infty_{\rm
    c}(\R^N\setminus\{0\},\C)$ with respect to the norms
  $(Q_{\A,a}(\cdot))^{1/2}$ and $\|\cdot\|_{\Distar}$ yields two
  coinciding spaces with equivalent norms.
\end{pfn}

By Hardy type inequalities, it is possible to compare the functional
space $\Distar$ with the classical Sobolev space $\Di$ defined as the
completion of $C^\infty_{\rm c}(\R^N,\C)$ with respect to the norm
\begin{displaymath}
  \|u\|_{\Di}:=\bigg(\int_{\R^N}\big|\nabla u(x)\big|^2\,dx\bigg)^{\!\!1/2}
\end{displaymath}
and with the space $\DiA$ given by the completion of $C^\infty_{\rm
  c}(\R^N\setminus\{0\},\C)$ with respect to the magnetic Dirichlet
norm
\begin{displaymath}
  \|u\|_{\DiA}:=\bigg(\int_{\R^N}\bigg|\nabla u(x)+i\,
  \frac{\A\big({x}/{|x|}\big)} {|x|}\,u(x)\bigg|^2\,dx\bigg)^{\!\!1/2}.
\end{displaymath}
The presence of a vector potential satisfying a suitable non-degeneracy
condition, allows recovering a Hardy's inequality even for
$N=2$. Indeed, if $N=2$, (\ref{eq:transversality}) holds, and
\begin{equation}\label{eq:circuit}
  \Phi_\A:=\frac1{2\pi}\int_0^{2\pi}\alpha(t)\,dt
  \not\in\Z,\quad \text{where }\alpha(t):=\A(\cos t,\sin t)\cdot(-\sin t,\cos
  t),
\end{equation}
then functions in ${\mathcal D}^{1,2}_\A(\R^2)$ satisfy the following Hardy
inequality
\begin{equation}\label{eq:hardyN2}
\Big(\min_{k\in\Z}|k-\Phi_\A|\Big)^2\int_{\R^2}\frac{|u(x)|^2}{|x|^2}\,dx
\leq \int_{\R^2}\bigg|\nabla u(x)+i\,
  \frac{\A\big({x}/{|x|}\big)} {|x|}\,u(x)\bigg|^2\,dx
\end{equation}
being $\Big(\min_{k\in\Z}|k-\Phi_\A|\Big)^2$
the best constant, as proved in \cite{lw}.
 It is easy to verify
  that, for $N=2$,
\[
\mu_1(\A,0)=\min_{\substack{\psi\in
    H^1((0,2\pi),\C)\\\psi(0)=\psi(2\pi)}}
\frac{\int_0^{2\pi}|\psi'(t)+i\alpha(t)
  \psi(t)|^2\,dt}{\int_0^{2\pi}|\psi(t)|^2\,dt},
\]
where $\alpha(t):=\A(\cos t,\sin t)\cdot(-\sin t,\cos
t)$. Furthermore, $\mu_1(\A,0)>0$ if and only if (\ref{eq:circuit})
holds. Combining Lemma \ref{l:la1mu1} (in the case $N=2$ and $a\equiv
0$) with \cite{lw}, we conclude that, for $N=2$,
\begin{equation}\label{eq:49}
 \mu_1(\A,0)=\Big(\min_{k\in\Z}|k-\Phi_\A|\Big)^2.
\end{equation}

\begin{Lemma}\label{l:spaces}

\noindent\begin{itemize}
\item[\rm (i)] If $N\geq 3$ then $\Distar=\Di$ and the norms $\|\cdot\|_{\Distar}$
and $\|\cdot\|_{\Di}$ are equivalent.
\item[\rm (ii)] If $\A\in C^1({\mathbb S}^{N-1},\R^N)$ and
either $N\geq 3$ or $N=2$ and
  (\ref{eq:transversality}), (\ref{eq:circuit}) hold, then
$\Distar=\DiA$ with equivalent norms.
\end{itemize}
\end{Lemma}
\begin{pf}
  By classical Hardy's inequality, for $N\geq 3$ the norms
  $\|\cdot\|_{\Di}$ and $\|\cdot\|_{\Distar}$ are equivalent over the
  space $C^\infty_c(\R^N\setminus\{0\},\C)$.  The proof of i) then
  follows by completion after observing that, for $N\geq 3$,
  $C^\infty_c(\R^N\setminus\{0\},\C)$ is dense in $\Di$.

  In order to prove ii), let us consider $u\in C^\infty_{\rm
    c}(\R^N\setminus\{0\},\C)$. Then
\begin{align*}
  \|u\|&_{\DiA} =\left\|\nabla u+i \frac{{\bf
        A}(x/|x|)}{|x|}u\right\|_{L^2(\R^N,\C^N)} \leq \|\nabla
  u\|_{L^2(\R^N,\C^N)}+\left\|
    \frac{{\bf A}(x/|x|)}{|x|}u\right\|_{L^2(\R^N)} \\
  & \leq \|\nabla u\|_{L^2(\R^N,\C^N)}+\sup_{\mathbb S^{N-1}} |{\bf
    A}| \left(\int_{\R^N}\frac{|u|^2}{|x|^2} \, dx\right)^{\!\!1/2}
  \leq  {\rm const}\|u\|_{\Distar}.
\end{align*}
On the other hand, by the diamagnetic inequality in Lemma \ref{l:diam},
  classical Hardy's inequality for $N\geq 3$, and \eqref{eq:hardyN2}
  for $N=2$, we have
\begin{align*}
\|u\|&_{\Di}=\|\nabla u\|_{L^2(\R^N,\C^N)}
\leq
\left\|\nabla u+i\frac{{\bf A}(x/|x|)}{|x|}u\right\|_{L^2(\R^N,\C^N)}
+\left\|\frac{{\bf A}(x/|x|)}{|x|}u\right\|_{L^2(\R^N,\C^N)}\\
&
\leq\|u\|_{\DiA}+\sup_{\mathbb S^{N-1}} |{\bf A}|
\left(\int_{\R^N}\frac{|u|^2}{|x|^2} \, dx\right)^{\!\!1/2}
\leq
\|u\|_{\DiA}+{\rm const}\,
\|\nabla |u|\|_{L^2(\R^N,\C^N)} \\
&
\leq
\big(1+{\rm const}\big)
\|u\|_{\DiA}.
\end{align*}
The above inequalities show that
$\|\cdot\|_{\Distar}$
and $\|\cdot\|_{\DiA}$ are equivalent norms over
the space $C^\infty_c(\R^N\setminus\{0\},\C)$.
The proof of the lemma then follows immediately
from the definition of the spaces $\Distar$ and $\DiA$.
\end{pf}

\section{A Hardy type inequality with boundary
  terms}\label{sec:hardyboundary}

We extend to singular electromagnetic potentials the Hardy type
inequality with boundary terms proved by Wang and Zhu in \cite{wz}
(see Lemma \ref{l:wz} in the Appendix).

\begin{Lemma}\label{l:hardyboundary}
 If $N\geq 2$, \eqref{eq:reg} and \eqref{eq:transversality} hold, then
  \begin{align}\label{eq:10}
    \int_{B_r} \bigg[ \bigg|\bigg(\nabla+i\,\frac
    {\A\big({x}/{|x|}\big)}{|x|}&\bigg)u\bigg|^2
    -\frac{a\big({x}/{|x|}\big)}{|x|^2}|u|^2\bigg]\,dx+
    \frac{N-2}{2r}\int_{\partial B_r}|u(x)|^2\,dS\\
    \notag&\geq \left(\mu_1(\A,a)+\bigg(\frac{N-2}{2}\bigg)^{\!\!2}\right)
      \int_{B_r}\frac{|u(x)|^2}{|x|^2}\,dx
  \end{align}
for all $r>0$ and $u\in H^1_{*}(B_r,\C)$.
  \end{Lemma}
\begin{pf}
  By scaling, it is enough to prove the inequality for $r=1$.  Let
  $u\in C^\infty(B_1,\C)\cap H^1_{*}(B_1,\C)$ with $0\not\in
  {\mathop{\rm supp}}u$.  Passing to polar coordinates and using
  \eqref{eq:46}, we have that
\begin{align}\label{eq:4}
  \int_{B_1} \bigg[ \bigg|&\bigg(\nabla+i\,\frac
  {\A\big({x}/{|x|}\big)}{|x|}\bigg)u\bigg|^2
  -\frac{a\big({x}/{|x|}\big)}{|x|^2}|u|^2\bigg]\,dx+
  \frac{N-2}{2}\int_{\partial B_1}|u(x)|^2\,dS\\
  \notag=& \int_{{\mathbb S}^{N-1}}\bigg(
  \int_0^{1}r^{N-1}|\partial_ru(r,\theta)|^2\,dr\bigg)\,dS(\theta)
  +\frac{N-2}2\int_{{\mathbb S}^{N-1}}|u(1,\theta)|^2\,dS(\theta)
  \\
  \notag& + \int_0^{1}\frac{r^{N-1}}{r^2}\bigg(\int_{{\mathbb
      S}^{N-1}}\left[|\nabla_{{\mathbb S}^{N-1}}u(r,\theta)
    +i\,\A(\theta)u(r,\theta)|^2-a(\theta)
    |u(r,\theta)|^2\right]\,dS(\theta)\bigg)\,dr.
\end{align}
For all $\theta\in{\mathbb S}^{N-1}$, let $\varphi_{\theta}\in
C^\infty_{\rm c}((0,+\infty),\C)$ be defined by
$\varphi_{\theta}(r)=u(r,\theta)$, and $\widetilde\varphi_{\theta}\in
C^\infty_{\rm c}(\R^N\setminus\{0\},\C)$ be the radially symmetric
function given by
$\widetilde\varphi_{\theta}(x)=\varphi_{\theta}(|x|)$. The Hardy
inequality with boundary term proved in \cite{wz} (see Lemma
\ref{l:wz} in the appendix) yields, for $N\geq 3$,
\begin{align}\label{eq:5}
  \int_{{\mathbb S}^{N-1}}&\bigg(
  \int_0^{1}r^{N-1}|\partial_ru(r,\theta)|^2\,dr
+\frac{N-2}2|u(1,\theta)|^2\bigg)\,dS(\theta)\\
\notag&=
  \int_{{\mathbb S}^{N-1}}\bigg(
  \int_0^{1}r^{N-1}|\varphi_{\theta}'(r)|^2\,dr
+\frac{N-2}2|\varphi_{\theta}(1)|^2\bigg)\,dS(\theta)\\
  \notag& =\frac1{\omega_{N-1}}  \int_{{\mathbb S}^{N-1}}\bigg(
\int_{B_1}|\nabla\widetilde \varphi_{\theta}(x)|^2\,dx
+\frac{N-2}2\int_{\partial B_1}|\widetilde \varphi_{\theta}(x)|^2\,dS
\bigg)\,dS(\theta)
\\&\notag\geq
\frac1{\omega_{N-1}}  \bigg(\frac{N-2}2\bigg)^{\!\!2}
  \int_{{\mathbb S}^{N-1}}\bigg(\int_{B_1}\frac{|\widetilde
    \varphi_{\theta}(x)|^2}{|x|^2}\,dx\bigg)\,dS(\theta)\\
\notag&=\bigg(\frac{N-2}2\bigg)^{\!\!2}
  \int_{{\mathbb S}^{N-1}}\bigg(
  \int_0^{1}\frac{r^{N-1}}{r^2}|u(r,\theta)|^2\,dr\bigg)\,dS(\theta)
=\bigg(\frac{N-2}2\bigg)^{\!\!2}\int_{B_1}\frac{|u(x)|^2}{|x|^2}\,dx.
\end{align}
On the other hand, (\ref{eq:5}) trivially holds also for $N=2$.
From (\ref{eq:4}), (\ref{eq:5}), and (\ref{eq:3}), we deduce that
\begin{align*}
  &\int_{B_1} \bigg[ \bigg|\bigg(\nabla+i\,\frac
  {\A\big({x}/{|x|}\big)}{|x|}\bigg)u\bigg|^2
  -\frac{a\big({x}/{|x|}\big)}{|x|^2}|u|^2\bigg]\,dx+
  \frac{N-2}{2}\int_{\partial B_1}|u(x)|^2\,dS \\
  &\quad\geq \left[\bigg(\frac{N-2}2\bigg)^{\!\!2}+\mu_1(\A,a)\right]
  \int_{B_1}\frac{|u(x)|^2}{|x|^2}\,dx \quad \text{for all } u\in
  C^\infty(B_1,\C)\cap H^1_{*}(B_1,\C)\text{ with }0\not\in {\mathop{\rm
      supp}}u,
\end{align*}
which, by density, yields the stated inequality for all
$H^1_{*}(B_r,\C)$-functions for $r=1$.
\end{pf}
\begin{remark}\label{r:hardyN2}
In view of (\ref{eq:49}), Lemma \ref{l:hardyboundary} for $N=2$ and $a\equiv 0$
yields
\begin{equation*}
  \int_{B_r} \bigg[ \bigg|\bigg(\nabla+i\,\frac
  {\A\big({x}/{|x|}\big)}{|x|}\bigg)u\bigg|^2 dx\geq
  \Big(\min_{k\in\Z}|k-\Phi_\A|\Big)^2\int_{B_r}\frac{|u(x)|^2}{|x|^2}\,dx
\end{equation*}
for all $r>0$ and $u\in H^1_{*}(B_r,\C)$.
 \end{remark}

\section{A Pohozaev-type identity}\label{sec:pohoz-type-ident}

Solutions to (\ref{u}) satisfy the following Pohozaev-type identity.

\begin{Theorem} \label{t:pohozaev} Let $\Omega\subset\R^N$, $N\geq 2$,
  be a bounded open set such that  $0\in\Omega$.
Let  $a,\A$ satisfy (\ref{eq:reg}), and $u$ be a
  weak $H^1_{*}(\Omega,\C)$-solution to (\ref{u}) in $\Omega$, with $h$
  satisfying (\ref{hph}).  Then
\begin{multline}\label{eq:26}
-\frac{N-2}2\int_{B_r}\bigg[ \bigg|\bigg(\nabla+i\,\frac
  {\A\big({x}/{|x|}\big)}{|x|}\bigg)u\bigg|^2
  -\frac{a\big({x}/{|x|}\big)}{|x|^2}|u|^2\bigg]\,dx\\
+\frac{r}{2}\int_{\partial B_r}\bigg[ \bigg|\bigg(\nabla+i\,\frac
  {\A\big({x}/{|x|}\big)}{|x|}\bigg)u\bigg|^2
  -\frac{a\big({x}/{|x|}\big)}{|x|^2}|u|^2\bigg]\,dS\\
=r\int_{\partial B_r}\bigg|\frac{\partial u}{\partial \nu}\bigg|^2\,dS+
\int_{B_r}\Re\big(h(x)u(x)\,(x\cdot
\overline{\nabla u(x)})\big)\,dx
\end{multline}
for all $r>0$ such that $\overline{B_{r}}=
\{x\in\R^N:|x|\leq r\}\subset\Omega$, where $\nu=\nu(x)$
is the unit outer normal vector $\nu(x)=\frac{x}{|x|}$.
\end{Theorem}
\begin{pf}
  Let $r>0$ such that $\overline{B_{r}}\subset\Omega$.  Since
\begin{multline*}
  \int_0^r\bigg[\int_{\partial B_s} \bigg[\bigg|\bigg(\nabla+i\,\frac
  {\A\big({x}/{|x|}\big)}{|x|}\bigg)u\bigg|^2 +\frac{|u|^2}{|x|^2}
  +\bigg|\frac{\partial u}{\partial
    \nu}\bigg|^2\bigg]\,dS\bigg]ds\\
  = \int_{B_r}\bigg[\bigg|\bigg(\nabla+i\,\frac
  {\A\big({x}/{|x|}\big)}{|x|}\bigg)u\bigg|^2
  +\frac{|u|^2}{|x|^2}+\bigg|\frac{\partial u}{\partial
    \nu}\bigg|^2\bigg]\,dx<+\infty
\end{multline*}
there exists a sequence $\{\delta_n\}_{n\in\N}\subset(0,r)$ such that
$\lim_{n\to+\infty}\delta_n=0$ and
\begin{equation}\label{eq:33}
\delta_n \int_{\partial B_{\delta_n}} \bigg[\bigg|\bigg(\nabla+i\,\frac
  {\A\big({x}/{|x|}\big)}{|x|}\bigg)u\bigg|^2 +\frac{|u|^2}{|x|^2}
  +\bigg|\frac{\partial u}{\partial
    \nu}\bigg|^2\bigg]\,dS\longrightarrow 0\quad \text{as }n\to+\infty.
\end{equation}
From classical regularity theory for elliptic equations, $u\in
W^{2,p}_{\rm loc}(\Omega\setminus\{0\})$ for all $p\in[1,\infty)$ and
$u\in C^{1,\tau}_{\rm loc}(\Omega\setminus\{0\},\C)$ for any
$\tau\in(0,1)$ (see Remark \ref{remA0}), hence we can multiply
equation (\ref{u}) by $x\cdot\overline{\nabla u(x)}$, integrate over
$B_r\setminus B_{\delta_n}$, and take the real part, thus obtaining
  \begin{gather}\label{eq:27}
    \int_{B_r\setminus B_{\delta_n}}\Re\big(\nabla
    u(x)\cdot\nabla(x\cdot\overline{\nabla u(x)})\big)\,dx+
    \int_{B_r\setminus
      B_{\delta_n}}\frac{|\A(\frac{x}{|x|})|^2-a(\frac{x}{|x|})}{|x|^2}
    \Re\big(u(x)\,(x\cdot\overline{\nabla u(x)})\big)\,dx\\
    \nonumber+ \int_{B_r\setminus B_{\delta_n}}
    \frac{\A(x/|x|)}{|x|}\cdot\Im\left(\overline u(x)\nabla(\nabla
      u(x)\cdot x) \right)\,dx +\int_{B_r\setminus B_{\delta_n}}
    \frac{\A(x/|x|)}{|x|}\cdot\Im\left((\overline {\nabla u(x)} \cdot x)
      \nabla u(x)
    \right) \, dx\\
    \nonumber=r\int_{\partial B_r}\bigg|\frac{\partial u}{\partial
      \nu}\bigg|^2\,dS -{\delta_n}\int_{\partial
      B_{\delta_n}}\bigg|\frac{\partial u}{\partial \nu}\bigg|^2\,dS
    +\int_{B_r\setminus B_{\delta_n}}\Re\big(h(x)\,
    u(x)(x\cdot\overline{\nabla u(x)})\big)\,dx.
  \end{gather}
Integration by parts yields
\begin{gather}\label{eq:29}
  \int_{B_r\setminus B_{\delta_n}}\nabla
  u(x)\cdot\nabla(x\cdot\overline{\nabla u(x)})\,dx
  =-(N-1)\int_{B_r\setminus B_{\delta_n}}|\nabla u(x)|^2\,dx+r\int_{\partial B_r} |\nabla
  u(x)|^2\,dS\\
  \nonumber-{\delta_n}\int_{\partial B_{\delta_n}} |\nabla
  u(x)|^2\,dS-\sum_{i,j=1}^N\int_{B_r\setminus B_{\delta_n}} x_j\overline{\frac{\partial
      u}{\partial x_i}}\frac{\partial^2 u}{\partial x_i\partial
    x_j}\,dx.
\end{gather}
A further integration by parts leads to
\begin{gather*}
  \sum_{i,j=1}^N\int_{B_r\setminus B_{\delta_n}} x_j\overline{\frac{\partial u}{\partial
      x_i}}\frac{\partial^2 u}{\partial x_i\partial
    x_j}\,dx=-N\int_{B_r\setminus B_{\delta_n}}|\nabla u(x)|^2\,dx +r \int_{\partial
    B_r}|\nabla u(x)|^2\,dS\\
  -{\delta_n} \int_{\partial B_{\delta_n}}|\nabla u(x)|^2\,dS -
  \sum_{i,j=1}^N\int_{B_r\setminus B_{\delta_n}} x_j\frac{\partial u}{\partial
    x_i}\overline{\frac{\partial^2 u}{\partial x_i\partial x_j}}\,dx
\end{gather*}
and hence
\begin{align}\label{eq:28}
  \sum_{i,j=1}^N\int_{B_r\setminus B_{\delta_n}} \Re\bigg(x_j\overline{\frac{\partial
      u}{\partial x_i}}\frac{\partial^2 u}{\partial x_i\partial
    x_j}\bigg)\,dx&=-\frac N2\int_{B_r\setminus B_{\delta_n}}|\nabla u(x)|^2\,dx \\
  \nonumber&\quad+\frac r2 \int_{\partial B_r}|\nabla u(x)|^2\,dS-\frac
  {\delta_n}2 \int_{\partial B_{\delta_n}}|\nabla u(x)|^2\,dS.
\end{align}
Collecting (\ref{eq:29}) and (\ref{eq:28}) we obtain
\begin{multline}\label{eq:30}
  \int_{B_r\setminus B_{\delta_n}}\Re\big(\nabla
  u(x)\cdot\nabla(x\cdot\overline{\nabla u(x)})\big)\,dx\\
  =-\frac{N-2}2\int_{B_r\setminus B_{\delta_n}}|\nabla u(x)|^2\,dx
  +\frac r2\int_{\partial B_r} |\nabla u(x)|^2\,dS
  -\frac{\delta_n}2\int_{\partial B_{\delta_n}} |\nabla u(x)|^2\,dS.
\end{multline}
Letting $f(\theta)=|\A(\theta)|^2-a(\theta)$, we have that $f\in
L^{\infty}({\mathbb S}^{N-1},\R)$ and, passing to polar coordinates
$r=|x|$, $\theta=\frac x{|x|}$, and observing
that $\partial_ru(r,\theta)=\nabla u(r\theta)\cdot \theta$,
\begin{align*}
  \int_{B_r\setminus B_{\delta_n}}&\frac{f(\frac x{|x|})}{|x|^2}
  u(x)\,(x\cdot\overline{\nabla u(x)})\,dx= \int_{{\mathbb
      S}^{N-1}}f(\theta)\bigg[\int_{\delta_n}^rs^{N-2}u(s\theta)
  \overline{\partial_s u(s\theta)}\,ds\bigg]dS(\theta)\\
  &=\int_{{\mathbb S}^{N-1}}f(\theta)\bigg[
  r^{N-2}|u(r\theta)|^2-\delta_n^{N-2}|u(\delta_n\theta)|^2\\
  &\hskip3cm- (N-2)\int_{\delta_n}^rs^{N-3}|u(s\theta)|^2\,ds
  -\int_{\delta_n}^rs^{N-2}\overline{u(s\theta)}
  \partial_s u(s\theta)\,ds\bigg]dS(\theta)\\
  &=r\int_{\partial B_r}\frac{f(\frac x{|x|})}{|x|^2} |
  u(x)|^2\,dS-\delta_n\int_{\partial B_{\delta_n}}\frac{f(\frac
    x{|x|})}{|x|^2}
  | u(x)|^2\,dS\\
  &\hskip2cm-(N-2)\int_{B_r\setminus B_{\delta_n}}\frac{f(\frac
    x{|x|})}{|x|^2} | u(x)|^2\,dx- \int_{B_r\setminus
    B_{\delta_n}}\frac{f(\frac x{|x|})}{|x|^2}
  \overline{u(x)}\,(x\cdot\nabla u(x))\,dx,
\end{align*}
thus leading to
\begin{multline}\label{eq:34}
  \int_{B_r\setminus B_{\delta_n}}\frac{|\A(\frac x{|x|})|^2-a(\frac
    x{|x|})}{|x|^2} \Re\big(u(x)\,(x\cdot\overline{\nabla
    u(x)})\big)\,dx\\
\quad= -\frac{N-2}2\int_{B_r\setminus
    B_{\delta_n}}\frac{|\A(\frac x{|x|})|^2-a(\frac x{|x|})}{|x|^2}
  | u(x)|^2\,dx+\frac r2\int_{\partial B_r}\frac{|\A(\frac
    x{|x|})|^2-a(\frac x{|x|})}{|x|^2} |
  u(x)|^2\,dS\\
   -\frac{\delta_n}2\int_{\partial
    B_{\delta_n}}\frac{|\A(\frac x{|x|})|^2-a(\frac x{|x|})}{|x|^2} |
  u(x)|^2\,dS.
\end{multline}
From integration by parts it follows
\begin{gather*}
  \int_{B_r\setminus B_{\delta_n}} \overline
  {u(x)}\,\frac{\A(x/|x|)}{|x|}\cdot\nabla(\nabla u(x)\cdot x) \,dx
  =-(N-2)\int_{B_r\setminus B_{\delta_n}}\overline{u(x)}\,
  \frac{\A(x/|x|)}{|x|}\cdot\nabla u(x) \,dx\\
  +r\int_{\partial B_r}\overline{u(x)}\,
  \frac{\A(x/|x|)}{|x|}\cdot\nabla u(x) \,dS- \delta_n\int_{\partial
    B_{\delta_n}}\overline{u(x)}\,
  \frac{\A(x/|x|)}{|x|}\cdot\nabla u(x) \,dS\\
  -\int_{B_r\setminus B_{\delta_n}} \frac{\A(x/|x|)}{|x|}\cdot\nabla
  u(x)\big(x\cdot \overline{\nabla u(x)}\big) \,dx
\end{gather*}
and therefore
\begin{multline}\label{eq:35}
  \int_{B_r\setminus B_{\delta_n}}
  \frac{\A(x/|x|)}{|x|}\cdot\Im\left(\overline u(x)\nabla(\nabla
    u(x)\cdot x) \right)\,dx +\int_{B_r\setminus B_{\delta_n}}
  \frac{\A(x/|x|)}{|x|}\cdot\Im\left((\overline {\nabla u(x)} \cdot x)
    \nabla u(x) \right) \, dx\\
=-(N-2)\int_{B_r\setminus B_{\delta_n}}
 \Im\bigg(\frac{\A(x/|x|)}{|x|}\cdot\nabla u(x)\,\overline {u(x)}\bigg) \,dx\\
+r\int_{\partial B_r}
\Im\bigg( \frac{\A(x/|x|)}{|x|}\cdot\nabla u(x) \,\overline {u(x)}\bigg)\,dS-
\delta_n\int_{\partial B_{\delta_n}}\Im\bigg(
 \frac{\A(x/|x|)}{|x|}\cdot\nabla u(x)\,\overline{u(x)}\bigg) \,dS.
\end{multline}
Putting together (\ref{eq:27}), (\ref{eq:30}), (\ref{eq:34}), and
(\ref{eq:35}) and taking into account that
\[
\bigg|\bigg(\nabla+i\,\frac
{\A\big({x}/{|x|}\big)}{|x|}\bigg)u\bigg|^2 =|\nabla
u|^2+2\,\frac{\A\big({x}/{|x|}\big)}{|x|}\cdot\Im(\overline{u}\nabla
u)+\frac{|\A\big({x}/{|x|}\big)|^2}{|x|^2}|u|^2,
\]
we obtain
\begin{multline*}
  -\frac{N-2}2\int_{B_r\setminus B_{\delta_n}}\bigg[
  \bigg|\bigg(\nabla+i\,\frac
  {\A\big({x}/{|x|}\big)}{|x|}\bigg)u\bigg|^2
  -\frac{a\big({x}/{|x|}\big)}{|x|^2}|u|^2\bigg]\,dx\\
  +\frac{r}{2}\int_{\partial B_r}\bigg[ \bigg|\bigg(\nabla+i\,\frac
  {\A\big({x}/{|x|}\big)}{|x|}\bigg)u\bigg|^2
  -\frac{a\big({x}/{|x|}\big)}{|x|^2}|u|^2\bigg]\,dS\\
  -\frac{\delta_n}{2}\int_{\partial B_{\delta_n}}\bigg[
  \bigg|\bigg(\nabla+i\,\frac
  {\A\big({x}/{|x|}\big)}{|x|}\bigg)u\bigg|^2
  -\frac{a\big({x}/{|x|}\big)}{|x|^2}|u|^2\bigg]\,dS
  \\
  =r\int_{\partial B_r}\bigg|\frac{\partial u}{\partial
    \nu}\bigg|^2\,dS -\delta_n\int_{\partial
    B_{\delta_n}}\bigg|\frac{\partial u}{\partial \nu}\bigg|^2\,dS+
  \int_{B_r\setminus B_{\delta_n}}\Re\big(h(x)u(x)\,(x\cdot\overline{\nabla
  u(x)})\big)\,dx.
\end{multline*}
Letting $n\to+\infty$ in the above identity and using (\ref{eq:33}) we obtain
(\ref{eq:26}).
\end{pf}

\section{The Almgren type frequency
  function}\label{sec:monot-prop}

Let  $u$ be a weak
$H^1_{*}(\Omega,\C)$-solution to equation (\ref{u}) in a bounded domain
$\Omega\subset \R^N$ containing the origin with $h$ satisfying
(\ref{hph}).  Let $\overline{R}>0$ be such that
$\overline{B_{\overline{R}}}= \{x\in\R^N:|x|\leq \overline
R\}\subseteq\Omega$.  Thus, the following functions are well defined
for every $r\in (0,\overline{R}]$:
\begin{equation} \label{D(r)} D(r)=\frac{1}{r^{N-2}} \int_{B_r}
  \left[\left|\nabla u(x)+i
      \frac{A(x/|x|)}{|x|}u(x)\right|^2-\frac{a(x/|x|)}{|x|^2}|u(x)|^2
    -(\Re h(x))|u(x)|^2 \right] \, dx,
\end{equation}
and
\begin{equation} \label{H(r)}
H(r)=\frac{1}{r^{N-1}}\int_{\partial B_r}|u|^2 \, dS.
\end{equation}
We are going to study regularity of functions $D$ and $H$. We first
differentiate  $H$.
\begin{Lemma} \label{l:hprime}
Let $\Omega\subset\R^N$, $N\geq 2$,
  be a bounded open set such that  $0\in\Omega$.
Let  $a,\A$ satisfy (\ref{eq:reg}), and $u$ be a
  weak $H^1_{*}(\Omega,\C)$-solution to (\ref{u}) in $\Omega$, with $h$
  satisfying (\ref{hph}). If $H$
  is the function defined in (\ref{H(r)}), then $H\in C^1(0,\overline
  R)$ and
\begin{equation}\label{H'}
  H'(r)=\frac{2}{r^{N-1}} \int_{\partial
    B_r} \Re\left(u\,\frac{\partial\, \overline u}{\partial \nu}\right)
  dS  \qquad \text{for every } r\in (0,\overline R).
\end{equation}
\end{Lemma}

\begin{pf}
Fix $r_0\in (0,\overline R)$ and consider the limit
\begin{equation} \label{limite} \lim_{r\rightarrow r_0}
  \frac{H(r)-H(r_0)}{r-r_0} = \lim_{r\rightarrow r_0} \int_{\partial
    B_1} \frac{|u(r\theta)|^2-|u(r_0\theta)|^2}{r-r_0}dS(\theta).
\end{equation}
Since $u\in C^1(B_{\overline R}\setminus \{0\},\C)$ (see Remark \ref{remA0})
then, for every $\theta\in{\partial B_1}$,
\begin{equation} \label{55} \lim_{r\rightarrow r_0}
  \frac{|u(r\theta)|^2-|u(r_0\theta)|^2}{r-r_0}
  =2\Re\left(\frac{\partial\, \overline u}{\partial
      \nu}(r_0\theta)u(r_0\theta)\right).
\end{equation}
On the other hand, for any $r\in (r_0/2,\overline R)$ and
$\theta\in{\partial B_1}$ we have
$$
\left| \frac{|u(r\theta)|^2-|u(r_0\theta)|^2}{r-r_0} \right| \leq 2
\sup_{B_{\overline R}\setminus B_{\frac{r_0}{2}}} |u| \cdot \sup_{B_{\overline R}\setminus
  B_{\frac{r_0}{2}}} |\nabla u|
$$
and hence, by (\ref{limite}), (\ref{55}), and the Dominated Convergence
Theorem, we obtain that
$$
H'(r_0)=\int_{\partial B_1} 2\Re\left(\frac{\partial\, \overline
    u}{\partial \nu}(r_0\theta)u(r_0\theta)\right) dS(\theta) =
\frac{2}{r_0^{N-1}} \int_{\partial B_{r_0}} \Re\left(u\frac{\partial
    \,\overline u}{\partial \nu}\right) dS.
$$
The continuity of $H'$ on the interval $(0,\overline R)$ follows by
the representation of $H'$ given above, the fact that $u\in
C^1(B_{\overline R}\setminus \{0\},\C)$, and the Dominated Convergence Theorem.
\end{pf}

\noindent In the lemma below, we study the regularity of the function
$D$.

\begin{Lemma}\label{l:dprime}
  Let $\Omega\subset\R^N$, $N\geq 2$, be a bounded open set such that
  $0\in\Omega$.  Let $a,\A$ satisfy (\ref{eq:reg}), and $u$ be a weak
  $H^1_{*}(\Omega,\C)$-solution to (\ref{u}) in $\Omega$, with $h$
  satisfying (\ref{hph}).  If $D$ is the function defined in
  (\ref{D(r)}), then $D\in W^{1,1}_{{\rm loc}}(0,\overline R)$.
  Moreover
\begin{align}\label{D'F}
  D'(r)=\frac{2}{r^{N-1}} \bigg[ r\int_{\partial B_r}
    \left|\frac{\partial u}{\partial \nu}\right|^2 dS &+\int_{B_r}
    \Re \big(h(x)\,\overline {u(x)}\, (x\cdot \nabla u(x))\big) \, dx\\
    \notag & + \frac{N-2}2
    \int_{B_r} (\Re
    h(x))|u(x)|^2\,dx-\frac r2\int_{\partial B_r} (\Re h(x))|u(x)|^2 \, dS \bigg]
\end{align}
in a distributional sense and for a.e. $r\in (0,\overline R)$.
\end{Lemma}

\begin{pf}
For any $r\in (0,\overline R)$ let
\begin{align} \label{I(r)} I(r)&= \int_{B_r} \left[\left|\nabla u(x)+i
      \frac{A(x/|x|)}{|x|}u(x)\right|^2-\frac{a(x/|x|)}{|x|^2}|u(x)|^2
    -(\Re h(x))|u(x)|^2 \right] \, dx  \\
  \notag &= \int_0^r \left( \int_{\partial B_\rho} \left[\left|\nabla
        u(x)+i \frac{A(x/|x|)}{|x|}u(x)\right|^2-\frac{a(x/|x|)}{|x|^2}|u(x)|^2
      -(\Re h(x))|u(x)|^2 \right] dS \right) d\rho.
\end{align}
From  the fact that $u\in H^1_{*}(B_{\overline R},\C)$,  we deduce that $I\in
W^{1,1}(0,\overline R)$ and
\begin{equation} \label{I'(r)} I'(r) = \int_{\partial B_r}
  \left[\left|\nabla u(x)+i
      \frac{A(x/|x|)}{|x|}u(x)\right|^2-\frac{a(x/|x|)}{|x|^2}|u(x)|^2 -(\Re
    h(x))|u(x)|^2 \right] dS
\end{equation}
 for a.e. $r\in (0,\overline R)$ and in the distributional sense.
Therefore by (\ref{eq:26}), (\ref{I(r)}), and (\ref{I'(r)}), we deduce
that $D\in W^{1,1}_{{\rm loc}}(0,\overline R)$ and
\begin{align}
  D'(r)&=r^{1-N}[-(N-2)I(r)+rI'(r)] \\
  \notag & =r^{1-N} \bigg[ 2r\int_{\partial B_r} \left|\frac{\partial
      u}{\partial \nu}\right|^2 dS
  +2\int_{B_r} \Re \big(h(x)\overline{ u(x)} (x\cdot \nabla u(x))\big) \, dx \\
  \notag & \quad+(N-2) \int_{B_r} (\Re h(x))|u(x)|^2\,dx -r\int_{\partial
    B_r} (\Re h(x))|u(x)|^2 \, dS \bigg],
\end{align}
 for a.e. $r\in (0,\overline R)$ and in the distributional sense.
This completes the proof of the lemma.
\end{pf}

\noindent We now show that $H(r)$ does not vanish for every $r>0$
sufficiently close to zero.

\begin{Lemma} \label{welld} Let $\Omega\subset\R^N$, $N\geq 2$, be a
  bounded open set such that $0\in\Omega$, $a,\A$ satisfy
  (\ref{eq:reg}), (\ref{eq:transversality}), (\ref{positivity}), and
  $u\not\equiv 0$ be a weak $H^1_{*}(\Omega,\C)$-solution to (\ref{u}) in
  $\Omega$, with $h$ satisfying (\ref{hph}). Let $H=H(r)$ be the
  function defined in (\ref{H(r)}). Then there exists $\overline r>0$
  such that $H(r)>0$ for any $r\in (0,\overline r)$.
\end{Lemma}

\begin{pf} Suppose by contradiction that there exists a sequence
  $r_n\to 0^+$ such that $H(r_n)=0$. Then for any $n$, $u\equiv
  0$ on $\partial B_{r_n}$.  Multiplying both sides of (\ref{u}) by
  $\overline u$ and integrating by parts over $B_{r_n}$ we obtain
\begin{align*}
\int_{B_{r_n}} & \left|\nabla u(x)+i\frac{\A(x/|x|)}{|x|}u(x)\right|^2 dx
-\int_{B_{r_n}} \frac{a(x/|x|)}{|x|^2}|u(x)|^2 dx \\
\notag &
=
\int_{B_{r_n}} h(x)|u(x)|^2 dx+
\int_{\partial B_{r_n}} \frac{\partial u}{\partial \nu}\,\overline u \, dS
=  \int_{B_{r_n}} h(x)|u(x)|^2 dx.
\end{align*}
Taking the real part on both sides it follows
\begin{equation*}
  \int_{B_{r_n}}  \left|\nabla u(x)+i\frac{\A(x/|x|)}{|x|}u(x)\right|^2 dx
  -\int_{B_{r_n}} \frac{a(x/|x|)}{|x|^2}|u(x)|^2 dx=
  \int_{B_{r_n}} \Re(h(x))|u(x)|^2 dx.
\end{equation*}
Since $u\equiv 0$ on $\partial B_{r_n}$, Lemma \ref{l:hardyboundary}
and (\ref{hph}) yield, for some positive constant $c_h>0$ depending
only on $h$,
\begin{align} \label{sdo} 0 & \geq \int_{B_{r_n}} \left|\nabla
    u(x)+i\frac{\A(x/|x|)}{|x|}u(x)\right|^2 dx -\int_{B_{r_n}}
  \frac{a(x/|x|)}{|x|^2}|u(x)|^2 dx -c_hr_n^\e
  \int_{B_{r_n}} \frac{|u(x)|^2}{|x|^2} dx  \\
  \notag & \geq \left(\mu_1(\A,a)+\bigg(\frac{N-2}{2}\bigg)^{\!\!2}-c_h
    r_n^\e \right) \int_{B_{r_n}}\frac{|u(x)|^2}{|x|^2}\,dx.
\end{align}
Since $\mu_1(\A,a)+\big(\frac{N-2}{2}\big)^2>0$ and $r_n\to 0^+$, we
conclude that $u\equiv 0$ in $B_{r_n}$ for $n$ sufficiently large.
Since $u\equiv 0$ in a neighborhood of the origin, we may apply, away
from the origin, a unique continuation principle for second order
elliptic equations with locally bounded coefficients (see e.g.
\cite{wolff}) to conclude that $u\equiv 0$ in $\Omega$, a
contradiction.
\end{pf}

\noindent
By virtue of Lemma \ref{welld}, the \emph{Almgren type frequency
  function}
\begin{equation} \label{N(r)}
{\mathcal N}(r)={\mathcal N}_{u,h}(r)=\frac{D(r)}{H(r)}
\end{equation}
is well defined in a suitably small interval $(0,\bar r)$.
Collecting Lemmas \ref{l:hprime} and \ref{l:dprime}, we
 compute the derivative of ${\mathcal N}$.

 \begin{Lemma} \label{mono} Let $\Omega\subset\R^N$, $N\geq 2$, be a
   bounded open set such that $0\in\Omega$, $a,\A$ satisfy
   (\ref{eq:reg}), (\ref{eq:transversality}), (\ref{positivity}), and
   $u\not\equiv 0$ be a weak $H^1_{*}(\Omega,\C)$-solution to (\ref{u})
   in $\Omega$, with $h$ satisfying (\ref{hph}). Then, letting
   ${\mathcal N}$ as in (\ref{N(r)}), there holds ${\mathcal N}\in
   W^{1,1}_{{\rm loc}}(0,\overline r)$ and
\begin{gather} \label{formulona} {\mathcal N}'(r)= \frac{2r\left[
      \left(\int_{\partial B_r} \left|\frac{\partial
            u}{\partial\nu}\right|^2 dS\right) \cdot
      \left(\int_{\partial B_r} |u|^2 dS\right)-\left(\int_{\partial
          B_r} \Re\left(u\frac{\partial \overline u}{\partial
            \nu}\right) dS\right)^{\!2} \right]}
  {\left(\int_{\partial B_r} |u|^2 dS\right)^2} \\
  \notag +\frac{2 \left[ \int_{B_r} \Re (h(x)\, \overline{u(x)}\, (x\cdot
      \nabla u(x))) \, dx + \frac{N-2}2 \int_{B_r} (\Re
      h(x))|u(x)|^2\,dx-\frac r2\int_{\partial B_r} (\Re h(x))|u(x)|^2
      \, dS \right]} {\int_{\partial B_r} |u|^2 dS}
\end{gather}
in a distributional sense and for a.e. $r\in (0,\overline r)$.
\end{Lemma}

\begin{pf} From Lemmas \ref{welld}, \ref{l:hprime}, and
  \ref{l:dprime}, it follows that ${\mathcal N}\in W^{1,1}_{{\rm
      loc}}(0,\overline r)$. Multiplying both sides of (\ref{u}) by
  $\overline u$, integrating by parts, and taking the real part we
  obtain the identity
$$
\int_{B_r} \left[\left|\nabla u(x)+i
    \frac{A(x/|x|)}{|x|}u(x)\right|^2-\frac{a(x/|x|)}{|x|^2}|u(x)|^2 -(\Re
  h(x))|u(x)|^2 \right] \, dx=\int_{\partial B_r} \Re\left(u\frac{\partial
    \,\overline u}{\partial \nu}\right) dS.
$$
Therefore, by (\ref{D(r)}) and (\ref{H'}) we infer
\begin{equation} \label{ul}
D(r)=\frac{1}{2}r H'(r)
\end{equation}
for every $r\in (0,\bar r)$.
From (\ref{ul}) we have that
$$
{\mathcal N}'(r)=\frac{D'(r)H(r)-D(r)H'(r)}{(H(r))^2}
=\frac{D'(r)H(r)-\frac{1}{2} r (H'(r))^2}{(H(r))^2}
$$
and, using (\ref{H'}) and (\ref{D'F}), the proof of the lemma easily follows.
\end{pf}

\noindent We now prove that ${\mathcal N}(r)$ admits a finite limit as
$r\to 0^+$.

\begin{Lemma} \label{gamma}
Under the same assumptions as in Lemma \ref{mono}, the limit
$$
\gamma:=\lim_{r\rightarrow 0^+} {\mathcal N}(r)
$$
exists and is finite.
\end{Lemma}

\begin{pf}
  We start by proving that ${\mathcal N}(r)$ is bounded from below as
  $r\rightarrow 0^+$.  By Lemma \ref{l:hardyboundary}, proceeding as
  in (\ref{sdo}) we arrive, for some positive constant $c_h>0$ depending
only on $h$, to
\begin{align} \label{47} \int_{B_{r}} & \left|\nabla
    u(x)+i\frac{\A(x/|x|)}{|x|}u(x)\right|^2 dx -\int_{B_{r}}
  \frac{a(x/|x|)}{|x|^2}|u(x)|^2 dx
  -\int_{B_r} (\Re h(x))|u(x)|^2 dx \\
  & \notag \geq -\frac{N-2}{2r} \int_{\partial B_r} |u(x)|^2 dS +
  \left(\mu_1(\A,a)+\bigg(\frac{N-2}{2}\bigg)^{\!\!2}-c_h r^\e \right)
  \int_{B_{r}}\frac{|u(x)|^2}{|x|^2}\,dx \\
  & \notag > -\frac{N-2}{2r} \int_{\partial B_r} |u(x)|^2 dS
\end{align}
for $r>0$ sufficiently small.
This with (\ref{D(r)})-(\ref{H(r)}) yields
 \begin{equation}\label{Nbelow}
   {\mathcal N}(r)>-\frac{N-2}{2}
 \end{equation}
 for any $r>0$ sufficiently close to zero.  Thanks to (\ref{hph}), for
 some $C_1>0$, we estimate
 \begin{multline*}
   \left| \int_{B_r} \Re \big(h(x) {\overline{u(x)}} (x\cdot \nabla
     u(x))\big) \, dx+ \frac{N-2}2 \int_{B_r} (\Re
     h(x))|u(x)|^2\,dx-\frac
     r2\int_{\partial B_r} (\Re h(x))|u(x)|^2 \, dS \right|\\
   \leq C_1 r^\e\left(\int_{B_r} \left|\nabla u+i
       \frac{A(x/|x|)}{|x|}u\right|^2\,dx+ \int_{B_r}
     \frac{|u(x)|^2}{|x|^2} \, dx+r^{N-2}H(r)\right).
 \end{multline*}
 Together with (\ref{47}), this implies that there exist $C_2>0$ and
 $\tilde r>0$ such that, for any $r\in (0,\tilde r)$,
\begin{multline*}
\left| \int_{B_r} \Re\big(h(x)
    \overline{u(x)} (x\cdot \nabla u(x))\big) \, dx+ \frac{N-2}2 \int_{B_r}
    (\Re h(x))|u(x)|^2\,dx-\frac r2\int_{\partial B_r} (\Re
    h(x))|u(x)|^2 \, dS \right| \\[5pt]
  \leq C_2\, r^{\e+N-2}\left[D(r)+H(r)\right].
\end{multline*}
Therefore, for any $r\in (0,\tilde r)$, we have that
\begin{gather}\label{B00}
  \left| \frac{ \int_{B_r} \Re\big(h(x)\overline{u(x)} (x\cdot \nabla
      u(x))\big) \, dx + \frac{N-2}2 \int_{B_r} (\Re
      h(x))|u(x)|^2\,dx-\frac r2\int_{\partial B_r} (\Re
      h(x))|u(x)|^2 \, dS }{\int_{\partial B_r} |u(x)|^2 dS} \right|\\[5pt]
  \notag\leq C_2 \,r^{-1+\e} \frac{D(r)+H(r)}{H(r)} \leq C_2\, r^{-1+\e}
  {\mathcal N}(r)+ C_2 \, r^{-1+\e}.
 \end{gather}
 By Lemma \ref{mono} and Schwarz's inequality, one sees that
$$
{\mathcal N}'(r)\geq 2\,
\frac{
\int_{B_r} \Re(h(x) \overline{u(x)} (x\cdot \nabla u(x))) \, dx
+ \frac{N-2}2 \int_{B_r} (\Re h(x))|u(x)|^2\,dx-\frac
      r2\int_{\partial B_r} (\Re
      h(x))|u(x)|^2 \, dS
}{\int_{\partial B_r} |u(x)|^2 dS}
$$
and hence by (\ref{B00}) we obtain
\begin{equation}\label{eq:40}
{\mathcal N}'(r)\geq
-2\,C_2\, r^{-1+\e} {\mathcal N}(r)-2\,C_2\, r^{-1+\e}
\end{equation}
for any $r\in (0,\tilde r)$.
After integration it follows that, for some $C_3>0$,
\begin{equation} \label{Nabove}
{\mathcal N}(r)\leq {\mathcal N}(\tilde r) e^{\frac{2C_2}{\e}(\tilde r^\e-r^\e)}
+2\,C_2 e^{-\frac{2C_2}{\e} r^{\e}} \int_r^{\tilde r}
s^{\e-1}e^{\frac{2C_2}{\e} s^\e} \, ds \leq C_3
\end{equation}
for any $r\in (0,\tilde r)$.  This shows that the left hand side of
(\ref{B00}) belongs to $L^1(0,\tilde r)$.  In particular by Lemma
\ref{mono} and Schwarz's inequality we see that ${\mathcal N}'$ is the
sum of a nonnegative function and of a $L^1$-function.  Therefore
$$
{\mathcal N}(r)={\mathcal N}(\tilde r)-\int_r^{\tilde r} {\mathcal N}'(s)\, ds
$$
admits a limit as $r\rightarrow 0^+$ which is necessarily finite in view of
(\ref{Nbelow}) and (\ref{Nabove}).
\end{pf}

\noindent A first consequence of the above analysis on the Almgren's
frequency function is the following estimate of $H(r)$.
\begin{Lemma} \label{l:uppb}
  Under the same assumptions as in Lemma \ref{mono}, let
  $\gamma:=\lim_{r\rightarrow 0^+} {\mathcal N}(r)$ be as in Lemma \ref{gamma}.
 Then there exists a constant
$K_1>0$ such that
\begin{equation} \label{1stest}
H(r)\leq K_1 r^{2\gamma}  \quad \text{for all } r\in (0,\bar r).
\end{equation}
On the other hand for any $\sigma>0$ there exists a constant
$K_2(\sigma)>0$ depending on $\sigma$ such that
\begin{equation} \label{2ndest}
H(r)\geq K_2(\sigma)\, r^{2\gamma+\sigma}   \quad \text{for all } r\in (0,\bar r).
\end{equation}
\end{Lemma}

\begin{pf}
  We start by proving (\ref{1stest}). Since, by Lemma \ref{gamma},
  ${\mathcal N}'\in L^1(0,\bar r )$ and ${\mathcal N}$ is bounded, then
  by (\ref{eq:40}), we infer that
  \begin{equation} \label{qsopra}
{\mathcal N}(r)-\gamma=\int_0^r
    {\mathcal N}'(s) \, ds\geq -C_4 r^\e
\end{equation}
for some constant $C_4>0$ and $r\in(0,\tilde r)$ with $0<\tilde r<\bar
r$.  Therefore by (\ref{ul}) and (\ref{qsopra}) we deduce that for
$r\in(0,\tilde r)$
$$
\frac{H'(r)}{H(r)}=\frac{2\,{\mathcal N}(r)}{r}\geq
\frac{2\gamma}{r}-2C_4 r^{-1+\e}.
$$
The proof of (\ref{1stest})  follows
immediately after integration in the previous differential inequality
over the interval $(r,\tilde r)$ and by continuity of $H$ outside $0$.

Let us prove (\ref{2ndest}). Since $\gamma=\lim_{r\rightarrow 0^+}
{\mathcal N}(r)$, for any $\sigma>0$ there exists $r_\sigma>0$ such
that ${\mathcal N}(r)<\gamma+\sigma/2$ for any $r\in (0,r_\sigma)$ and
hence
$$
\frac{H'(r)}{H(r)}=\frac{2\,{\mathcal N}(r)}{r}<\frac{2\gamma+\sigma}{r}
\quad \text{for all } r\in (0,r_\sigma).
$$
Integrating over the interval $(r,r_\sigma)$ and by continuity of $H$
outside $0$, we obtain (\ref{2ndest}) for some constant $K_2(\sigma)$
depending on $\sigma$.
\end{pf}

\section{Proofs of  Theorems \ref{t:asym} and
  \ref{asy_infinity}}\label{sec:proofs-theor-reft}

In this section we use the monotonicity properties established in
section \ref{sec:monot-prop} combined with a blow-up technique to
deduce asymptotics of solutions near the singularity and to prove
Theorems~\ref{t:asym} and \ref{asy_infinity}.

\begin{Lemma}\label{l:blowup}
  Let $\Omega\subset\R^N$, $N\geq 2$, be a bounded open set containing
  $0$, $a, \A$ such that (\ref{eq:reg}), (\ref{eq:transversality}),
  and (\ref{positivity}) hold, and $h$ as in
  (\ref{hph}).  For $u\in H^1_{*}(\Omega,\C)$ weakly solving (\ref{u}),
  $u\not\equiv 0$, let $\gamma:=\lim_{r\rightarrow 0^+} {\mathcal
    N}(r)$ as in Lemma \ref{gamma}. Then
\begin{itemize}
\item[\rm (i)] there exists $k_0\in \N$ such that
  $\gamma=-\frac{N-2}2+\sqrt{\left(\frac{N-2}{2}\right)^2+\mu_{k_0}(\A,a)}$;
\item[\rm (ii)] for every sequence $\lambda_n\to0^+$, there exist a subsequence
$\{\lambda_{n_k}\}_{k\in\N}$ and an eigenfunction $\psi$ of the operator
$L_{\A,a}$ associated to the eigenvalue $\mu_{k_0}(\A,a)$ such that
$\|\psi\|_{L^{2}({\mathbb S}^{N-1},\C)}=1$ and
\[
\frac{u(\lambda_{n_k}x)}{\sqrt{H(\lambda_{n_k})}}\to
|x|^{\gamma}\psi\Big(\frac x{|x|}\Big)
\]
weakly in $H^1(B_1,\C)$, strongly in $H^1(B_r,\C)$ for every $0<r<1$,
and in $C^{1,\tau}_{\rm loc}(B_1\setminus\{0\},\C)$ for any $\tau\in (0,1)$.
\end{itemize}
\end{Lemma}
\begin{pf}
Let us set
\[
w^\lambda(x)=\frac{u(\lambda x)}{\sqrt{H(\lambda)}}.
\]
We notice that $\int_{\partial B_1}|w^{\lambda}|^2dS=1$. Moreover, by
scaling and (\ref{Nabove}),
\begin{gather}\label{eq:8}
  \int_{B_{1}} \left|\nabla w^\lambda(x)+ i\frac{\A(\frac
      x{|x|})}{|x|}w^\lambda(x)\right|^2 dx -\int_{B_{1}}
  \frac{a(\frac x{|x|})}{|x|^2}|w^\lambda(x)|^2 dx
  -\int_{B_1} \lambda^2(\Re h(\lambda x))|w^\lambda(x)|^2 dx \\[5pt]
  \notag ={\mathcal N}(\lambda)\leq {\rm const}.
\end{gather}
Hence, by (\ref{eq:10}) and (\ref{hph}) there exists $c_h>0$ such that
\[
\bigg(\mu_1(\A,a)+\bigg(\frac{N-2}{2}\bigg)^{\!\!2}-c_h\lambda^\e\bigg)
\int_{B_{1}}
  \frac{|w^\lambda(x)|^2}{|x|^2}\, dx\leq \frac{N-2}2+{\mathcal N}(\lambda),
\]
and, consequently, there exist $\bar\lambda>0$ and ${\rm const}>0$ such that
\[
\int_{B_{1}} \frac{|w^\lambda(x)|^2}{|x|^2}\, dx\leq {\rm const} \quad\text{for
  every }0<\lambda<\bar\lambda,
\]
which, in view of (\ref{eq:8}), implies
that $\{w^\lambda\}_{\lambda\in(0,\bar\lambda)}$ is bounded in
$H^1_{*}(B_1,\C)$.

Therefore, for any given sequence $\lambda_n\to 0^+$, there exists a
subsequence $\lambda_{n_k}\to0^+$ such that $w^{\lambda_{n_k}}\weakly
w$ weakly in $H^1_{*}(B_1,\C)$ for some $w\in H^1_{*}(B_1,\C)$. We notice
that $H^1_{*}(B_1,\C)$ is continuously
embedded into $H^1(B_1,\C)$, hence $w^{\lambda_{n_k}}\weakly w$ weakly also
in $H^1(B_1,\C)$.  Due to compactness of the trace imbedding
$H^1(B_1,\C)\hookrightarrow L^2(\partial B_1,\C)$, we obtain that
$\int_{\partial B_1}|w|^2dS=1$. In particular $w\not\equiv
0$. Furthermore, weak convergence allows passing to the weak limit in
the equation
\begin{equation}\label{eq:19}
\mathcal L_{\A,a} w^{\lambda_{n_k}}(x)={\lambda^2_{n_k}}
h({\lambda_{n_k}} x) w^{\lambda_{n_k}}(x)
\end{equation}
which holds in a weak sense in
$B_{{\overline{R}}/{\lambda_{n_k}}}\supset B_1$ (see the beginning of
section \ref{sec:monot-prop} for the definition of $\overline R$),
thus yielding
\begin{equation} \label{eq:w}
\mathcal L_{\A,a} w(x)=0\quad\text{in }B_1.
\end{equation}
A bootstrap argument and classical regularity theory lead to
\[
w^{\lambda_{n_k}}\to w\quad\text{in }C^{1,\tau}_{\rm loc}(B_{1}\setminus \{0\},\C)
\]
for any $\tau\in (0,1)$
and
\begin{equation} \label{strongH1}
w^{\lambda_{n_k}}\to w\quad\text{in } H^1(B_r,\C)
\text{ and in } H^1_{*}(B_r,\C)
\end{equation}
for any $r\in (0,1)$.
Since the functions $w^{\lambda_{n_k}}$ solve
 equation (\ref{eq:19}), then for any $r\in (0,1)$ we may define
the functions
\begin{align*}
  D_k(r)=& \frac{1}{r^{N-2}} \int_{B_r} \left[\left|\nabla
      w^{\lambda_{n_k}}(x)+i \frac{\A(x/|x|)}{|x|}
      w^{\lambda_{n_k}}(x) \right|^2
  \right]\,dx \\
  \notag & -\frac{1}{r^{N-2}} \int_{B_r} \left[
    \frac{a(x/|x|)}{|x|^2}|w^{\lambda_{n_k}}(x)|^2 +\lambda_{n_k}^2
    (\Re h(\lambda_{n_k}x))|w^{\lambda_{n_k}}(x)|^2 \right] \, dx
\end{align*}
and
\begin{equation*}
H_k(r)=\frac{1}{r^{N-1}}\int_{\partial B_r}|w^{\lambda_{n_k}}|^2 \, dS.
\end{equation*}
On the other hand, since $w$ solves (\ref{eq:w}), then we put
\begin{equation} \label{Dw(r)}
  D_w(r)=\frac{1}{r^{N-2}} \int_{B_r}
  \left[\left|\nabla w(x)+i
      \frac{\A(x/|x|)}{|x|}w(x)\right|^2-\frac{a(x/|x|)}{|x|^2}|w(x)|^2
  \right] \, dx \quad \text{for all } r\in (0,1)
\end{equation}
and
\begin{equation} \label{Hw(r)}
H_w(r)=\frac{1}{r^{N-1}}\int_{\partial B_r}|w|^2 \, dS
\quad  \text{for all } r\in (0,1).
\end{equation}
Using a change of variables, one sees that
\begin{equation}\label{NkNw}
{\mathcal
    N}_k(r):=\frac{D_k(r)}{H_k(r)}=\frac{D(\lambda_{n_k}r)}{H(\lambda_{n_k}r)}
  ={\mathcal N}(\lambda_{n_k}r) \quad \text{for all } r\in (0,1).
\end{equation}
By (\ref{hph}) and (\ref{strongH1}), we have
for any fixed $r\in (0,1)$
\begin{equation} \label{convDk}
 D_k(r)\to D_w(r).
\end{equation}
On the other hand, by compactness of the trace imbedding
$H^1(B_r,\C)\hookrightarrow L^2(\partial B_r,\C)$, we also have
\begin{equation} \label{convHk}
H_k(r)\to H_w(r) \quad \text{for any fixed } r\in (0,1).
\end{equation}
From (\ref{eq:10}) it follows that $D_w(r)>-\frac{N-2}2H_w(r)$ for all
$r\in(0,1)$.  Therefore, if, for some $r\in(0,1)$, $H_w(r)=0$ then
$D_w(r)>0$, and passing to the limit in (\ref{NkNw}) should give a
contradiction with Lemma \ref{gamma}. Hence $H_w(r)>0$ for all
$r\in(0,1)$. Thus the function
\[
{\mathcal N}_w(r):=\frac{D_w(r)}{H_w(r)}
\]
is well defined for $r\in (0,1)$.
This, together with (\ref{NkNw}), (\ref{convDk}), (\ref{convHk}), and
Lemma \ref{gamma}, shows that
\begin{equation} \label{Nw(r)} {\mathcal
    N}_w(r)=\lim_{k\to \infty} {\mathcal
    N}(\lambda_{n_k}r)=\gamma
\end{equation}
for all $r\in (0,1)$.
Therefore ${\mathcal N}_w$ is  constant in $(0,1)$ and hence ${\mathcal
  N}_w'(r)=0$ for any $r\in (0,1)$.  By (\ref{eq:w}) and Lemma
\ref{mono} with $h\equiv 0$, we obtain
\begin{equation*}
  \left(\int_{\partial B_r} \left|\frac{\partial
        w}{\partial\nu}\right|^2 dS\right) \cdot \left(\int_{\partial B_r} |w|^2
    dS\right)-\left(\int_{\partial B_r} \Re\left(w\frac{\partial \overline
        w}{\partial \nu}\right) dS\right)^{\!\!2}=0 \quad
  \text{for all } r\in (0,1),
\end{equation*}
i.e.
$$
\left|\int_{\partial B_r} \Re\left(w\,\frac{\partial\overline
      w}{\partial \nu}\right) \, dS\right|^2
=\left\|w\right\|^2_{L^2({\partial B_r},\C)} \cdot \left\|\frac{\partial
    w}{\partial \nu}\right\|^2_{L^2({\partial B_r},\C)}.
$$
This shows that $w$ and $\frac{\partial w}{\partial \nu}$
have the same direction as vectors in $L^2(\partial B_r,\C)$ and hence
there exists a real valued function $\eta=\eta(r)$ such that
$\frac{\partial w}{\partial \nu}(r,\theta)=\eta(r) w(r,\theta)$ for $r\in(0,1)$.
After integration we obtain
\begin{equation} \label{separate}
w(r,\theta)=e^{\int_1^r \eta(s)ds} w(1,\theta)
=\varphi(r) \psi(\theta) \quad  r\in(0,1), \ \theta\in \SN,
\end{equation}
where we put $\varphi(r)=e^{\int_1^r \eta(s)ds}$ and $\psi(\theta)=w(1,\theta)$.
Since
$$
\mathcal L_{\A,a}w=-\frac{\partial^2 w}{\partial
  r^2}-\frac{N-1}{r}\frac{\partial w}{\partial r} +\frac{1}{r^2}
L_{\A,a} w,
$$
then (\ref{separate}) yields
$$
\left(-\varphi''(r)-\frac{N-1}{r} \varphi'(r) \right)\psi(\theta)
+\frac{\varphi(r)}{r^2} L_{\A,a} \psi(\theta)=0.
$$
Taking $r$ fixed we deduce that $\psi$ is an eigenfunction of the
operator $L_{\A,a}$. If $\mu_{k_0}(\A,a)$ is the corresponding
eigenvalue then $\varphi(r)$ solves the equation
$$
-\varphi''(r)-\frac{N-1}{r} \varphi(r)+\frac{\mu_{k_0}(\A,a)}{r^2}\varphi(r)=0
$$
and hence $\varphi(r)$ is of the form
$$
\varphi(r)=c_1 r^{\sigma_{k_0}^+}+c_2 r^{\sigma_{k_0}^-}
$$
for some $c_1,c_2\in\R$, where
\begin{equation*}
  \sigma^+_{k_0}=-\frac{N-2}{2}+\sqrt{\bigg(\frac{N-2}
    {2}\bigg)^{\!\!2}+\mu_{k_0}(\A,a)}\quad\text{and}\quad
  \sigma^-_{k_0}=-\frac{N-2}{2}-\sqrt{\bigg(\frac{N-2}{2}
    \bigg)^{\!\!2}+\mu_{k_0}(\A,a)}.
\end{equation*}
Since the function
$\frac1{|x|}\big(|x|^{\sigma_{k_0}^-}\psi(\frac{x}{|x|})\big)\notin
L^2(B_1,\C)$ and hence $|x|^{\sigma_{k_0}^-}\psi(\frac{x}{|x|})\notin
H^1_{*}(B_1,\C)$, then $c_2=0$ and  $\varphi(r)=c_1
r^{\sigma_{k_0}^+}$. Since $\varphi(1)=1$, we obtain that $c_1=1$ and
then
\begin{equation} \label{expw}
w(r,\theta)=r^{\sigma_{k_0}^+} \psi(\theta),  \quad
\text{for all }r\in (0,1)\text{ and }\theta\in \SN.
\end{equation}
It remains to prove part (i).
Since $w$ solves (\ref{eq:w}), after integration by parts
$$
\int_{B_r} \left[\left|\nabla w(x)+i
\frac{A(x/|x|)}{|x|}w(x)\right|^2-\frac{a(x/|x|)}{|x|^2}|w(x)|^2 \right]
\, dx=\int_{\partial B_r} \frac{\partial w}{\partial \nu}
\,\overline w \, dS.
$$
Therefore, by (\ref{Dw(r)}), (\ref{Hw(r)}), (\ref{Nw(r)}) and
(\ref{expw}), it follows
$$
\gamma={\mathcal N}_w(r)=\frac{D_w(r)}{H_w(r)}=\frac{r \int_{\partial B_r}
\frac{\partial w}{\partial \nu} \overline w \, dS}{\int_{\partial
B_r} |w|^2 dS }
=\sigma_{k_0}^+.
$$
This completes the proof of the lemma.
\end{pf}

\noindent A further step towards a-priori bounds for solutions to (\ref{u})
relies in uniformly estimating the supremum of $|u|$ on $\partial B_r$
with $H(r)$.
\begin{Lemma}\label{l:8}
  Let $\Omega\subset\R^N$, $N\geq 2$, be a bounded open set containing
  $0$, $a, \A$ such that (\ref{eq:reg}), (\ref{eq:transversality})
  and (\ref{positivity}) hold, and $h$ as in
  (\ref{hph}).  Then, for any weak $H^1_{*}(\Omega,\C)$-solution $u$ to
  (\ref{u}) there exist $\bar s>0$ and $C>0$ such that
\[
\sup_{\partial B_{s}}|u|^2\leq \frac{C}{s^{N-1}}\int_{\partial
  B_{s}}|u|^2\,dS \quad\text{for every }0<s<\bar s.
\]
\end{Lemma}
\begin{pf}
  Let $\gamma=\lim_{r\rightarrow 0^+} {\mathcal N}(r)$ as in Lemma \ref{gamma}
  and $k_0\in \N$ such that
  \[
  \gamma=-\frac{N-2}2+\sqrt{\bigg(\frac{N-2}{2}\bigg)^{\!\!2}+\mu_{k_0}(\A,a)},
  \]
  see Lemma \ref{l:blowup}. Denote as ${\mathcal A}_0$ the eigenspace
  of the operator $L_{\A,a}$ associated to the eigenvalue
  $\mu_{k_0}(\A,a)$. Since $\mathop{\rm dim}{\mathcal A}_0$ is finite,
  it is easy to verify that
\[
\Lambda=\sup_{v\in{\mathcal A}_0\setminus \{0\}}\frac{\sup_{{\mathbb
      S}^{N-1}}|v|^2}{\int_{{\mathbb S}^{N-1}}|v|^2\,dS}<+\infty.
\]
Let $\widetilde C>2^{N-1}\Lambda$. We claim that there exists
$\bar\lambda$ such that
\begin{equation}\label{eq:claim}
  \sup_{\partial B_{1/2}}|w^{\lambda}|^2\leq
   \widetilde C\int_{\partial B_{1/2}}|w^{\lambda}|^2dS\quad
  \text{for every} \ \ \lambda\in(0,\bar\lambda).
\end{equation}
To prove (\ref{eq:claim}), assume by contradiction that there exists a
sequence $\{\lambda_{n}\}_{n\in\N}$ such that $\lambda_n\to0^+$ and
\begin{equation}\label{eq:18}
\sup_{\partial B_{1/2}}|w^{\lambda_{n}}|^2>\widetilde C\int_{\partial
  B_{1/2}}|w^{\lambda_{n}}|^2dS.
\end{equation}
Lemma \ref{l:blowup} implies that
there exist a subsequence $\{\lambda_{n_j}\}_{j\in\N}$ and an
eigenfunction $\psi\in{\mathcal A}_0$  such that $\|\psi\|_{L^{2}({\mathbb
    S}^{N-1},\C)}^2=1$ and $w^{\lambda_{n_j}}\to
|x|^{\gamma}\psi\Big(\frac x{|x|}\Big)$ weakly in
$H^1(B_1,\C)$ and in $C^{1,\tau}_{\rm loc}(B_1\setminus\{0\},\C)$
 for any $\tau\in(0,1)$.
Passing to limit in (\ref{eq:18}), this should
imply that
\[
\sup_{{\mathbb S}^{N-1}}|\psi|^2\geq \frac{\widetilde C}{2^{N-1}}
\int_{{\mathbb S}^{N-1}}|\psi|^2dS
>\Lambda \int_{{\mathbb S}^{N-1}}|\psi|^2dS
\]
giving rise to a contradiction with the definition of $\Lambda$. Claim
(\ref{eq:claim}) is thereby proved.

Estimate (\ref{eq:claim}) can be written as
\begin{equation*}
  \sup_{\partial B_{\lambda/2}}|u|^2\leq
  \frac{\widetilde C}{\lambda^{N-1}}\int_{\partial B_{\lambda/2}}|u|^2dS\quad
  \text{for every }\lambda\in(0,\bar\lambda).
\end{equation*}
Choosing $\bar s=\frac12\bar \lambda$ and $C=2^{1-N}\widetilde C$, the conclusion
follows.
\end{pf}

\noindent From Lemmas  \ref{l:uppb} and \ref{l:8} we deduce the following
pointwise estimate for solutions to (\ref{u}).
\begin{Corollary}\label{c:upbound}
  Let $\Omega\subset\R^N$, $N\geq 2$, be a bounded open set containing
  $0$, $a, \A$ such that (\ref{eq:reg}), (\ref{eq:transversality})
  and (\ref{positivity}) hold, and $h$ as in
  (\ref{hph}). Then, for any weak $H^1_{*}(\Omega,\C)$-solution $u$ to
  (\ref{u}) there exist $\bar s>0$ and $\bar C>0$ such that
\[
|u(x)|\leq \bar C\,|x|^{\gamma} \quad\text{for every }x\in B_{\bar s},
\]
where $\gamma=\lim_{r\rightarrow 0^+} {\mathcal N}(r)$ as in Lemma \ref{gamma}.
\end{Corollary}
\begin{pf}
It follows from (\ref{1stest}) and Lemma \ref{l:8}.
\end{pf}

\noindent Let us now describe the behavior of $H(r)$ as $r\to 0^+$.
\begin{Lemma} \label{l:limite}
Under the same assumptions as in Lemma
  \ref{mono} and letting $\gamma:=\lim_{r\rightarrow 0^+} {\mathcal
    N}(r)\in \R$ as in Lemma \ref{gamma}, the limit
\[
\lim_{r\to0^+}r^{-2\gamma}H(r)
\]
exists and it is finite.
\end{Lemma}
\begin{pf}
In view of (\ref{1stest}) it is sufficient to prove that the limit exists.
By (\ref{H(r)}), (\ref{ul}), and Lemma~\ref{gamma} we have
$$
\frac{d}{dr} \frac{H(r)}{r^{2\gamma}}
=-2\gamma r^{-2\gamma-1} H(r)+r^{-2\gamma} H'(r)
=2r^{-2\gamma-1} (D(r)-\gamma H(r))=2r^{-2\gamma-1} H(r) \int_0^r {\mathcal N}'(s) ds.
$$
Denote by $M_1(r)$ and $M_2(r)$ respectively the first and the second
term in the right hand side of (\ref{formulona}) in order to obtain,
after integration over $(r,\tilde r)$,
\begin{equation} \label{inte} \frac{H(\tilde r)}{\tilde
    r^{2\gamma}}-\frac{H(r)}{r^{2\gamma}}=\int_r^{\tilde r} 2s^{-2\gamma-1}
  H(s) \left( \int_0^s M_1(t) dt \right) ds +\int_r^{\tilde r} 2s^{-2\gamma-1}
  H(s) \left( \int_0^s M_2(t) dt \right) ds.
\end{equation}
By Schwarz's inequality we have that $M_1(t)\geq 0$ and hence
$$
\lim_{r\to 0^+} \int_r^{\tilde r} 2s^{-2\gamma-1} H(s) \left( \int_0^s
  M_1(t) dt \right) ds
$$
exists.  On the other hand, by (\ref{B00}) and (\ref{1stest}) we
deduce that $|M_2(r)|=O(r^{-1+\e})$ and $H(r)=O(r^{2\gamma})$ as $r\to
0^+$. Therefore, if $\tilde r$ is sufficiently small, for some
${\rm const\,}>0$ there holds
$$
\left|
s^{-2\gamma-1} H(s) \left( \int_0^s M_2(t) dt \right)
\right|\leq \frac{\rm const}{\e} s^{-1+\e}
$$
for all $r\in(0,\tilde r)$, which proves that $s^{-2\gamma-1} H(s)
\left( \int_0^s M_2(t) dt \right)\in L^1(0,\tilde r)$.  We may
conclude that both terms in the right hand side of (\ref{inte}) admit
a limit as $r\to 0^+$ thus completing the proof of the lemma.
\end{pf}

\noindent The limit $\lim_{r\to0^+}r^{-2\gamma}H(r)$ is indeed
strictly positive, as we prove in the following lemma.

\begin{Lemma} \label{l:limitepositivo}
Under the same assumptions as in Lemma
  \ref{mono} and letting $\gamma:=\lim_{r\rightarrow 0^+} {\mathcal
    N}(r)\in \R$ as in Lemma \ref{gamma},  there holds
\[
\lim_{r\to0^+}r^{-2\gamma}H(r)>0.
\]
\end{Lemma}
\begin{pf}
  Let us fix $R>0$ such that $\overline{B_{R}}\subset\Omega$.
 For any $k\in\N\setminus\{0\}$,
let $\psi_k$   be a $L^2$-normalized
eigenfunction of the operator $L_{\A,a}$ on the sphere associated to
the $k$-th eigenvalue $\mu_{k}(\A,a)$, i.e. satisfying
\begin{equation}\label{eq:2rad}
\begin{cases}
L_{\A,a}\psi_k(\theta)
=\mu_k(\A,a)\,\psi_k(\theta),&\text{in }{\mathbb S}^{N-1},\\[3pt]
\int_{{\mathbb S}^{N-1}}|\psi_k(\theta)|^2\,dS(\theta)=1.
\end{cases}
\end{equation}
We can choose the functions $\psi_k$ in such a way that they form an
orthonormal basis of $L^2({\mathbb S}^{N-1},\C)$, hence $u$ and $hu$
 can be expanded as
\begin{equation}\label{eq:21}
u(x)=u(\lambda\,\theta)=\sum_{k=1}^\infty\varphi_k(\lambda)\psi_k(\theta),
\quad
h(x)u(x)=h(\lambda\,\theta)u(\lambda\,\theta)=
\sum_{k=1}^\infty\zeta_k(\lambda)\psi_k(\theta),
\end{equation}
where $\lambda=|x|\in(0,R]$, $\theta=x/|x|\in{{\mathbb S}^{N-1}}$, and
\begin{equation}\label{eq:22}
  \varphi_k(\lambda)=\int_{{\mathbb S}^{N-1}}u(\lambda\,\theta)
  \overline{\psi_k(\theta)}\,dS(\theta),
  \quad
  \zeta_k(\lambda)=\int_{{\mathbb S}^{N-1}}
  h(\lambda\,\theta)u(\lambda\,\theta)\overline{\psi_k(\theta)}\,dS(\theta).
\end{equation}
Equations (\ref{u}) and (\ref{eq:2rad}) imply that, for every $k$,
\begin{equation*}
  -\varphi_k''(\lambda)-\frac{N-1}{\lambda}\varphi_k^\prime(\lambda)+
  \frac{\mu_k(\A,a)}{\lambda^2}\varphi_k(\lambda)=
  \zeta_k(\lambda),\quad\text{in }(0,R).
\end{equation*}
A direct calculation shows that, for some $c_1^k,c_2^k\in\R$,
\begin{equation}\label{eq:42}
\varphi_k(\lambda)=\lambda^{\sigma^+_k}
\bigg(c_1^k+\int_\lambda^R\frac{s^{-\sigma^+_k+1}}{\sigma^+_k-\sigma^-_k}
\zeta_k(s)\,ds\bigg)+\lambda^{\sigma^-_k}
\bigg(c_2^k+\int_\lambda^R\frac{s^{-\sigma^-_k+1}}{\sigma^-_k-\sigma^+_k}
\zeta_k(s)\,ds\bigg),
\end{equation}
where
\begin{equation*}
  \sigma^+_k=-\frac{N-2}{2}+\sqrt{\bigg(\frac{N-2}
    {2}\bigg)^{\!\!2}+\mu_k(\A,a)}\quad\text{and}\quad
  \sigma^-_k=-\frac{N-2}{2}-\sqrt{\bigg(\frac{N-2}{2}\bigg)^{\!\!2}+\mu_k(\A,a)}.
\end{equation*}
In view of Lemma \ref{l:blowup}, there exist $j_0,m\in\N$, $j_0,m\geq 1$
such that $m$ is the
multiplicity of the eigenvalue
$\mu_{j_0}(\A,a)=\mu_{j_0+1}(\A,a)=\cdots=\mu_{j_0+m-1}(\A,a)$ and
\begin{equation}\label{eq:15}
  \gamma=\lim_{r\rightarrow 0^+} {\mathcal N}(r)=\sigma_{i}^+,
  \quad i=j_0,\dots,j_0+m-1.
\end{equation}
The Parseval identity yields
\begin{equation}\label{eq:17}
H(\lambda)=\int_{{\mathbb
    S}^{N-1}}|u(\lambda\,\theta)|^2\,dS(\theta)=
\sum_{k=1}^{\infty}|\varphi_k(\lambda)|^2,\quad\text{for all }0<\lambda\leq R.
\end{equation}
Let us assume by contradiction that
$\lim_{\lambda\to0^+}\lambda^{-2\gamma}H(\lambda)=0$ and fix
$i\in\{j_0,\dots,j_0+m-1\}$. Then, (\ref{eq:15}) and (\ref{eq:17}) imply that
\begin{equation}\label{eq:11}
\lim_{\lambda\to0^+}\lambda^{-\sigma_{i}^+}\varphi_{i}(\lambda)=0.
\end{equation}
From (\ref{hph}) and Corollary \ref{c:upbound}, we obtain that
\begin{equation}\label{eq:zeta}
  \zeta_{i}(\lambda)=O(\lambda^{-2+\e+\sigma_{i}^+})\quad\text{as }\lambda\to 0^+,
\end{equation}
and, consequently, the functions
\[
s\mapsto \frac{s^{-\sigma^+_{i}+1}}{\sigma^+_{i}-\sigma^-_{i}}
\zeta_{i}(s)\quad\text{and}\quad s\mapsto
\frac{s^{-\sigma^-_{i}+1}}{\sigma^-_{i}-\sigma^+_{i}} \zeta_{i}(s)
\]
belong to $L^1((0,R),\C)$. Hence
\[
\lambda^{\sigma^+_{i}}
\bigg(c_1^{i}+\int_\lambda^R\frac{s^{-\sigma^+_{i}+1}}{\sigma^+_{i}-\sigma^-_{i}}
\zeta_{i}(s)\,ds\bigg)=o(\lambda^{\sigma^-_{i}})\quad\text{as }\lambda\to0^+,
\]
and then, since $\frac{u}{|x|}\in L^2(B_R,\C)$ and
$\frac{|x|^{\sigma^-_{i}}}{|x|}\not
\in L^2(B_R,\C)$, we conclude that there must be
\begin{equation*}
c_2^{i}=-\int_0^R\frac{s^{-\sigma^-_{i}+1}}{\sigma^-_{i}-\sigma^+_{i}}
\,\zeta_{i}(s)\,ds.
\end{equation*}
Using (\ref{eq:zeta}), we then deduce that
\begin{align}\label{eq:12}
  \lambda^{\sigma^-_{i}}
  \bigg(c_2^{i}+\int_\lambda^R\frac{s^{-\sigma^-_{i}+1}}{\sigma^-_{i}-\sigma^+_{i}}
  \zeta_{i}(s)\,ds\bigg)&=\lambda^{\sigma^-_{i}}
  \bigg(\int_0^\lambda
  \frac{s^{-\sigma^-_{i}+1}}{\sigma^+_{i}-\sigma^-_{i}}
  \zeta_{i}(s)\,ds\bigg)=O(\lambda^{\sigma^+_{i}+\e})
\end{align}
as $\lambda\to0^+$.  From (\ref{eq:42}), (\ref{eq:11}), and
(\ref{eq:12}), we obtain that
\[
c_1^{i}+\int_0^R\frac{s^{-\sigma^+_{i}+1}}{\sigma^+_{i}-\sigma^-_{i}}
\zeta_{i}(s)\,ds=0,
\]
thus implying, together with (\ref{eq:zeta}),
\begin{align}\label{eq:13}
  \lambda^{\sigma^+_{i}}
  \bigg(c_1^{i}+\int_\lambda^R\frac{s^{-\sigma^+_{i}+1}}{\sigma^+_{i}-\sigma^-_{i}}
  \zeta_{i}(s)\,ds\bigg)= \lambda^{\sigma^+_{i}}
  \int_0^\lambda\frac{s^{-\sigma^+_{i}+1}}{\sigma^-_{i}-\sigma^+_{i}}
  \zeta_{i}(s)\,ds=O(\lambda^{\sigma^+_{i}+\e})
\end{align}
as $\lambda\to0^+$. Collecting (\ref{eq:42}), (\ref{eq:12}), and
(\ref{eq:13}), we conclude that
\[
\varphi_{i}(\lambda)=O(\lambda^{\sigma^+_{i}+\e})\quad\text{as
}\lambda\to0^+\quad\text{for every }i\in\{1,\dots,m\},
\]
namely, setting $u^\lambda(\theta)=u(\lambda \theta)$,
\[
(u^\lambda,\psi)_{L^2({\mathbb S}^{N-1},\C)}=
O(\lambda^{\gamma+\e})\quad\text{as
}\lambda\to0^+
\]
for every $\psi\in {\mathcal A}_0$, where ${\mathcal A}_0$ is the eigenspace
of the operator $L_{\A,a}$ associated to the eigenvalue
$\mu_{j_0}(\A,a)=\mu_{j_0+1}(\A,a)=\cdots=\mu_{j_0+m-1}(\A,a)$.
  Let
$w^\lambda(\theta)=(H(\lambda))^{-1/2}u(\lambda \theta)$.  From (\ref{2ndest}),
there exists $C(\e)>0$ such that $\sqrt{H(\lambda)}\geq
C(\e)\lambda^{\gamma+\frac\e2}$ for $\lambda$ small, and therefore
\begin{equation}\label{eq:14}
(w^\lambda,\psi)_{L^2({\mathbb S}^{N-1},\C)}=
O(\lambda^{\e/2})=o(1)\quad\text{as
}\lambda\to0^+
\end{equation}
for every $\psi\in{\mathcal A}_0$.  From Lemma \ref{l:blowup}, for
every sequence $\lambda_n\to0^+$, there exist a subsequence
$\{\lambda_{n_j}\}_{j\in\N}$ and an eigenfunction $\widetilde
\psi\in{\mathcal A}_0$ such that
\begin{equation}\label{eq:16}
\int_{{\mathbb S}^{N-1}}|\widetilde\psi(\theta)|^2dS=1\quad\text{and} \quad
w^{\lambda_{n_j}}\to \widetilde \psi\quad\text{in } L^2({\mathbb
  S}^{N-1},\C).
\end{equation}
From (\ref{eq:14}) and (\ref{eq:16}), we infer that
\[
0=\lim_{j\to+\infty}(w^{\lambda_{n_j}},\widetilde\psi)_{L^2({\mathbb S}^{N-1},\C)}
=\|\widetilde\psi\|_{L^2({\mathbb S}^{N-1},\C)}^2=1,
\]
thus reaching a contradiction.
\end{pf}

\noindent The analysis carried out in this section leads to a complete
description of the behavior of solutions to (\ref{u}) near the
singularity and hence to the proof of Theorem \ref{t:asym}.

\medskip\noindent \begin{pfn}{Theorem \ref{t:asym}} Identity
  (\ref{eq:20}) follows from part (i) of Lemma \ref{l:blowup}, thus
  there exists $k_0\in \N$, $k_0\geq 1$, such that $\lim_{r\to
    0^+}{\mathcal N}_{u,h}(r)=-\frac{N-2}2+\sqrt{\big(\frac{N-2}{2}
    \big)^{\!2}+\mu_{k_0}(\A,a)}$.  Let us denote as $m$ the
  multiplicity of $\mu_{k_0}(\A,a)$, so that, for some $j_0\in\N$,
  $j_0\geq 1$, $j_0\leq k_0\leq j_0+m-1$,
  $\mu_{j_0}(\A,a)=\mu_{j_0+1}(\A,a)=\cdots=\mu_{j_0+m-1}(\A,a)$ and
  let $\{\psi_i:\,j_0\leq i\leq j_0+m-1\}$ be an $L^2({\mathbb
    S}^{N-1},\C)$-orthonormal basis for the eigenspace of $L_{\A,a}$
  associated to $\mu_{k_0}(\A,a)$.  Set
  \[
  \gamma=-\frac{N-2}2+\sqrt{\bigg(\frac{N-2}{2}\bigg)^{\!\!2}+\mu_{k_0}(\A,a)}
  \]
and let $\{\lambda_n\}_{n\in\N}\subset (0,+\infty)$ such that
$\lim_{n\to+\infty}\lambda_n=0$. Then, from part (ii) of Lemma
\ref{l:blowup} and Lemmas \ref{l:limite} and \ref{l:limitepositivo},
there exist a subsequence $\{\lambda_{n_k}\}_{k\in\N}$ and $m$ real
numbers $\beta_{j_0},\dots,\beta_{j_0+m-1}\in\R$ such that
$(\beta_{j_0},\beta_{j_0+1},\dots,\beta_{j_0+m-1})\neq(0,0,\dots,0)$
and
\begin{equation}\label{eq:23}
\lambda_{n_k}^{-\gamma}u(\lambda_{n_k}\theta)\to
\sum_{i=j_0}^{j_0+m-1} \beta_i\psi_{i}(\theta)\quad \text{in }
C^{1,\tau}({\mathbb S}^{N-1},\C) \quad \text{as }k\to+\infty
\end{equation}
and
\begin{equation} \label{eq:23grad} \lambda_{n_k}^{1-\gamma}\nabla
  u(\lambda_{n_k}\theta)\to \sum_{i=j_0}^{j_0+m-1}
  \beta_i(\gamma\psi_{i}(\theta)\theta+\nabla_{\SN} \psi_i(\theta))
  \quad \text{in } C^{0,\tau}({\mathbb S}^{N-1},\C^N) \quad \text{as
  }k\to+\infty
\end{equation}
for any $\tau\in(0,1)$. We now prove
that the $\beta_i$'s depend neither on the sequence
$\{\lambda_n\}_{n\in\N}$ nor on its subsequence
$\{\lambda_{n_k}\}_{k\in\N}$.

Let us fix $R>0$ such that
$\overline{B_{R}}\subset\Omega$.  Defining $\varphi_i$ and $\zeta_i$
as in (\ref{eq:22}) and expanding $u$ as in (\ref{eq:21}), from
(\ref{eq:23}) it follows that, for any $i=j_0,\dots, j_0+m-1$,
\begin{equation}\label{eq:25}
\lambda_{n_k}^{-\gamma}\varphi_i(\lambda_{n_k}) =\int_{{\mathbb
    S}^{N-1}}\frac{u(\lambda_{n_k}\theta)}{\lambda_{n_k}^{\gamma}}
\overline{\psi_i(\theta)}\,dS(\theta)
\to\sum_{j=j_0}^{j_0+m-1} \beta_j\int_{{\mathbb
    S}^{N-1}}\psi_{j}(\theta)\overline{\psi_i(\theta)}\,dS(\theta)=\beta_i
\end{equation}
as $k\to+\infty$.  As deduced in the proof of Lemma
\ref{l:limitepositivo}, for any $i=j_0,\dots, j_0+m-1$
and $\lambda\in(0,R]$ there holds
\begin{align}\label{eq:24}
\varphi_i(\lambda)&=\lambda^{\sigma^+_i}
\bigg(c_1^i+\int_\lambda^R\frac{s^{-\sigma^+_i+1}}{\sigma^+_i-\sigma^-_i}
\zeta_i(s)\,ds\bigg)+\lambda^{\sigma^-_{i}}
  \bigg(\int_0^\lambda
  \frac{s^{-\sigma^-_i+1}}{\sigma^+_i-\sigma^-_i}
  \zeta_i(s)\,ds\bigg)\\
&\notag=\lambda^{\sigma^+_i}
\bigg(c_1^i+\int_\lambda^R\frac{s^{-\sigma^+_i+1}}{\sigma^+_i-\sigma^-_i}
\zeta_i(s)\,ds\bigg)+O(\lambda^{\sigma^+_i+\e})\quad\text{as }\lambda\to0^+,
\end{align}
for some $c_1^i\in\R$, where
\begin{equation*}
  \sigma^+_i=\gamma=-\frac{N-2}{2}+\sqrt{\bigg(\frac{N-2}
    {2}\bigg)^{\!\!2}+\mu_{k_0}(\A,a)},\quad
  \sigma^-_i=-\frac{N-2}{2}-\sqrt{\bigg(\frac{N-2}{2}
    \bigg)^{\!\!2}+\mu_{k_0}(\A,a)}.
\end{equation*}
Choosing $\lambda=R$ in the first line of (\ref{eq:24}), we obtain
\[
c_1^i=R^{-\sigma^+_i}\varphi_i(R)-R^{\sigma^-_i-\sigma^+_i}\int_0^R
  \frac{s^{-\sigma^-_i+1}}{\sigma^+_i-\sigma^-_i}
  \zeta_i(s)\,ds.
\]
Hence (\ref{eq:24}) yields
\[
\lambda^{-\gamma}\varphi_i(\lambda)\to
R^{-\sigma^+_i}\varphi_i(R)-R^{\sigma^-_i-\sigma^+_i}\int_0^R
  \frac{s^{-\sigma^-_i+1}}{\sigma^+_i-\sigma^-_i}
  \zeta_i(s)\,ds+\int_0^R\frac{s^{-\sigma^+_i+1}}{\sigma^+_i-\sigma^-_i}
\zeta_i(s)\,ds\quad\text{as }\lambda\to0^+,
\]
and therefore, from (\ref{eq:25}) we deduce that
\begin{align*}
  \beta_i&= R^{-\gamma}\int_{{\mathbb S}^{N-1}}u(R\theta)
  \overline{\psi_{i}(\theta)}\,dS(\theta)
  \\
  &\quad-R^{-2\gamma-N+2}\int_{0}^R\frac{s^{\gamma+N-1}}{2\gamma+N-2}\bigg(
  \int_{{\mathbb S}^{N-1}}
  h(s\,\eta)u(s\,\eta)\overline{\psi_{i}(\eta)}\,dS(\eta)\bigg) ds\\
  &\quad +\int_{0}^R\frac{s^{1-\gamma}}{2\gamma+N-2}\bigg( \int_{{\mathbb
      S}^{N-1}}
  h(s\,\eta)u(s\,\eta)\overline{\psi_{i}(\eta)}\,dS(\eta)\bigg) ds .
\end{align*}
In particular the $\beta_i$'s depend neither on the sequence
$\{\lambda_n\}_{n\in\N}$ nor on its subsequence
$\{\lambda_{n_k}\}_{k\in\N}$, thus implying that the convergences in
(\ref{eq:23}) and \eqref{eq:23grad} actually hold as $\lambda\to 0^+$
and proving the theorem.~\end{pfn}

\medskip\noindent \begin{pfn}{Corollary \ref{cor:holdercont}}
  Statement (i) follows directly from (\ref{estu}). Statement (iii) is
  an immediate consequence of (\ref{estu}) and (\ref{estgradu}).  To
  prove (ii), we notice that classical elliptic regularity theory
  yields H\"older continuity away from $0$, so it remains to prove
  that $u$ is H\"older continuous in every $\overline{B_{r}}\subset
  \Omega$. To this aim, we argue by contradiction and assume that
  there exist sequences $\{x_n\}_{n\in\N}$,
  $\{y_n\}_{n\in\N}\subset\overline{B_{r}}$ such that
  \begin{equation}\label{eq:45}
    \lim_{n\to+\infty}\frac{|u(x_n)-u(y_n)|}{|x_n-y_n|^\gamma}=+\infty.
  \end{equation}
H\"older continuity away from $0$ implies that either $|x_n|\to 0$ or
$|y_n|\to 0$ along a subsequence. Hence we can assume without loss of
generality that $|y_n|\to 0$ and $|x_n|\geq|y_n|$.  Two cases can occur.
\begin{description}
\item[Case 1] there exists a positive constant $c>1$ such that
  $\frac{|x_n|}{|y_n|}\leq c$. In this case, $|x_n|\to 0$ and, letting
  $\lambda_n=2c|x_n|$ and observing that $\frac{x_n}{\lambda_n},
  \frac{y_n}{\lambda_n}\in \overline B_{1/(2c)}\setminus {B_{1/(2c^2)}}
\Subset
  B_1\setminus\{0\}$, from part (ii) of Lemma \ref{l:blowup}
and Lemmas \ref{l:limite} and \ref{l:limitepositivo}
it follows
\[
\lim_{n\to+\infty}\frac{\left|\lambda_n^{-\gamma}
    u\big(\lambda_n\frac{x_n}{\lambda_n}\big)
    -(2c)^{-\gamma}\psi\big(\frac{x_n}{|x_n|}\big)
    -\lambda_n^{-\gamma}u\big(\lambda_n\frac{y_n}{\lambda_n}\big)+
    \frac{|y_n|^\gamma}{\lambda_n^\gamma}
    \psi\big(\frac{y_n}{|y_n|}\big)\right|}
{\big|\frac{x_n}{\lambda_n}-\frac{y_n}{\lambda_n}\big|^\gamma}=0
\]
for some eigenfunction $\psi$ of the operator $L_{\A,a}$. Since the
function $|x|^{\gamma}\psi\big(\frac x{|x|}\big)$ is H\"older
continuous away from $0$, from above we conclude that
\[
\frac{|u(x_n)-u(y_n)|}{|x_n-y_n|^\gamma}=
\frac{\left|\lambda_n^{-\gamma}u\big(\lambda_n\frac{x_n}{\lambda_n}\big)
-\lambda_n^{-\gamma}u\big(\lambda_n\frac{y_n}{\lambda_n}\big)\right|}
{\big|\frac{x_n}{\lambda_n}-\frac{y_n}{\lambda_n}\big|^\gamma}
\]
is bounded uniformly in $n$, thus giving rise to a contradiction.
\item[Case 2] There exists subsequences $\{x_{n_k}\}_{k\in\N}$ and
  $\{y_{n_k}\}_{k\in\N}$ such that
  $\frac{|x_{n_k}|}{|y_{n_k}|}\to+\infty$. In particular
$|y_{n_k}|=o(|x_{n_k}|)$ as $k\to+\infty$. From (\ref{eq:45}) we
  deduce that $|x_{n_k}|\to 0$ as $k\to+\infty$ and by Corollary
  \ref{c:upbound}
\begin{align*}
\frac{|u(x_{n_k})-u(y_{n_k})|}{|x_{n_k}-y_{n_k}|^\gamma}
&=|x_{n_k}|^{-\gamma}\frac{\big|u(x_{n_k})-u(y_{n_k})\big|}
{\big|\frac{x_{n_k}}{|x_{n_k}|}-\frac{y_{n_k}}{|x_{n_k}|}\big|^\gamma}\\[5pt]
&\leq {\rm const}\,
|x_{n_k}|^{-\gamma}\frac{|x_{n_k}|^\gamma+|y_{n_k}|^\gamma}
{\big|\frac{x_{n_k}}{|x_{n_k}|}-\frac{y_{n_k}}{|x_{n_k}|}\big|^\gamma}
\leq{\rm const}
\end{align*}
thus giving rise to a contradiction with (\ref{eq:45}).
\end{description}

\end{pfn}

\noindent Invariance by Kelvin's transform allows rewriting equations
in exterior domains as equations in bounded neighborhoods of $0$, thus
reducing the problem of asymptotics at infinity to the problem of
asymptotics at $0$. Hence we can deduce Theorem \ref{asy_infinity} from
Theorem \ref{t:asym}.

\medskip\noindent
\begin{pfn}{Theorem \ref{asy_infinity}}
Let $u$ be a weak solution of (\ref{u}) where $\Omega$ is
an external domain as in the statement of the theorem.
Let $v$ be the Kelvin transform of $u$, i.e.
\begin{equation} \label{Kelvin}
v(x)=|x|^{2-N} u\left(\frac{x}{|x|^2}\right),
\quad  x\in \widetilde \Omega=\big\{x\in\R^N: x/{|x|^2}\in\Omega\big\}.
\end{equation}
If we put $y=\frac{x}{|x|^2}$, then
we have
\begin{equation} \label{DELTA} \Delta u(x)=|y|^{N+2} \Delta v(y) \quad
  \text{for all }y\in\widetilde\Omega,
\end{equation}
and
\begin{align} \label{DIV}
  &\frac{a(x/|x|)-|\A(x/|x|)|^2+i{\rm div}_{\SN} \A(x/|x|)}{|x|^2}u(x) \\
  \notag &=|y|^{N+2} \frac{a(y/|y|)-|\A(y/|y|)|^2+i{\rm div}_{\SN}
    \A(y/|y|)}{|y|^2}v(y) \quad \text{for all }y\in\widetilde\Omega.
\end{align}
Moreover, by the transversality assumption (\ref{eq:transversality})
we also have
\begin{equation} \label{NABLA}
\frac{\A(x/|x|)}{|x|}\cdot \nabla u(x)
=|y|^{N+2}\frac{\A(y/|y|)}{|y|}\cdot \nabla v(y)\quad
  \text{for all }y\in\widetilde\Omega.
\end{equation}
Therefore, by (\ref{Kelvin}--\ref{NABLA}) we obtain
\begin{equation} \label{EQ_TRAS}
\mathcal L_{\A,a}v(y)= |y|^{-4} h\left(\frac{y}{|y|^2}\right) v(y)
\quad \text{in }\widetilde \Omega\setminus \{0\}.
\end{equation}
From a direct computation we infer that $\nabla v\in L^2(\widetilde
\Omega,\C^N)$, $\frac{v}{|x|}\in L^{2}(\widetilde \Omega,\C)$, and hence $v\in
H^1_{*}(\widetilde \Omega,\C)$.  This is sufficient for proving that $v$ is a
$H^1_{*}$-weak solution of equation (\ref{EQ_TRAS}) in $\widetilde \Omega$.

 On the other hand, by (\ref{hph_00})
$$
\left| |y|^{-4}h\left(\frac{y}{|y|^2}\right)\right|=O(|y|^{-2+\e}),
\qquad {\rm as} \ |y|\rightarrow 0^+
$$
and hence $v$ satisfies all the assumptions of Theorem \ref{t:asym}.
The proof of \eqref{eq:20_00} and the asymptotic estimate for $u$ then
follows by Theorem \ref{t:asym}, (\ref{Kelvin}), and the fact that
\begin{equation} \label{NN}
{\mathcal N}_{v,|y|^{-4}h(y/|y|^2)}(r)
=  \widetilde{\mathcal N}_{u,h}\big({\textstyle{\frac1r}}\big)-N+2
\end{equation}
with $\widetilde{\mathcal N}_{u,h}$ as in (\ref{eq:37ext}).  For
proving the estimate on the gradient one may proceed as follows.  Let
$\tilde\gamma$ be as in the statement of the theorem and let
$\gamma=\lim_{r\rightarrow 0^+}\mathcal N_{v,|y|^{-4}h(y/|y|^2)}(r)$.
From \eqref{NN} it follows that $\gamma=\tilde\gamma-N+2$, hence by
\eqref{Kelvin} we have
\begin{equation} \label{***} \lambda^{1-\gamma}\nabla v(\lambda\theta)
  =(2-N)\lambda^{-\tilde\gamma}u\left(\tfrac{\theta}{\lambda}\right)\theta
  +\lambda^{-\tilde\gamma-1}\nabla
  u\left(\tfrac{\theta}{\lambda}\right)
  -2\lambda^{-\tilde\gamma-1}\left(\nabla
    u\left(\tfrac{\theta}{\lambda}\right)\cdot \theta\right)\theta
\end{equation}
for any $\lambda$ such that $B_{\lambda}\subset \widetilde \Omega$ and
for any $\theta\in \SN$.  Applying Theorem \ref{t:asym} to the
function $v$, from the previous identity we infer
$$
(2-N)\lambda^{-\tilde\gamma}u\left(\tfrac{\theta}{\lambda}\right)
-\lambda^{-\tilde\gamma-1}\left(\nabla
  u\left(\tfrac{\theta}{\lambda}\right)\cdot \theta\right) \rightarrow
\gamma \sum_{i=j_0}^{j_0+m-1} \widetilde\beta_i \psi_i(\theta)
$$
in $C^{0,\tau}(\SN,\C)$ for any $\tau\in(0,1)$ as $\lambda\rightarrow 0^+$.
From the first part of the theorem we also have that
\begin{equation} \label{firstpart}
\lambda^{-\tilde\gamma}u\left(\tfrac{\theta}{\lambda}\right)\rightarrow
\sum_{i=j_0}^{j_0+m-1} \widetilde\beta_i \psi_i(\theta)
\end{equation}
from which we obtain
\begin{equation}\label{sopra**}
  \lambda^{-\tilde\gamma-1}\left(\nabla
    u\left(\tfrac{\theta}{\lambda}\right)\cdot
    \theta\right)\rightarrow -\tilde\gamma\sum_{i=j_0}^{j_0+m-1}
  \widetilde\beta_i \psi_i(\theta)
\end{equation}
in $C^{0,\tau}(\SN,\C)$ for any $\tau\in(0,1)$ as $\lambda\rightarrow
0^+$.  Letting $\lambda\rightarrow 0^+$ in \eqref{***}, applying again
Theorem \ref{t:asym} to the function $v$ and using
\eqref{firstpart}-\eqref{sopra**} we deduce that
$$
\lambda^{-\tilde\gamma-1}\nabla u\left(\tfrac{\theta}{\lambda}\right)
\rightarrow \sum_{i=j_0}^{j_0+m-1} \widetilde\beta_i(-\tilde \gamma
\psi_i(\theta)\theta+\nabla_{\SN} \psi_i(\theta))
$$
in $C^{0,\tau}(\SN,\C^N)$ for any $\tau\in(0,1)$ as $\lambda\rightarrow 0^+$.
By replacing $\lambda$ with $1/\lambda$ we obtain the desired estimate.
\end{pfn}

\section{An example: Aharonov-Bohm magnetic potentials in dimension
  $2$}\label{sec:ahar-bohm-magn}

In this section we discuss an application of Theorem \ref{t:asym} to
Schr\"odinger equations with Aharonov-Bohm vector potentials
(\ref{eq:31}), i.e. we let $N=2$, $\A(\cos t,\sin t)=\alpha(-\sin
t,\cos t)$, $a(\cos t,\sin t)=a_0$ for some
$a_0\in\R$, and consider the corresponding equation
\[
\left(-i\,\nabla+\alpha\left(-\frac{x_2}{|x|^2},\frac{x_1}{|x|^2}
\right)\right)^{\!\!2}u-\dfrac{a_0}{|x|^2}u=h\,u,
\]
with $x=(x_1,x_2)$ in a bounded domain of $\R^2$ containing $0$ and
$h$ verifying (\ref{hph}). In this case, an explicit calculation yields
\[
\{\mu_k(\A,a):k\in\N\setminus\{0\}\}=\{(\alpha-j)^2-a_0:j\in\Z\}
\]
hence, in particular,
\[
\mu_1(\A,a)=\big(\mathop{\rm dist}(\alpha,\Z)\big)^2-a_0.
\]
If $\mathop{\rm dist}(\alpha,\Z)\neq\frac12$, then all eigenvalues are
simple and the eigenspace associated to the eigenvalue
$(\alpha-j)^2-a_0$ is generated by $\psi(\cos t,\sin t)=e^{-ijt}$.  If
$\mathop{\rm dist}(\alpha,\Z)=\frac12$, then all eigenvalues have
multiplicity $2$.
Theorem \ref{t:asym} hence yields:
\begin{itemize}
\item[i)] if $a_0<\big(\mathop{\rm
    dist}(\alpha,\Z)\big)^2$ and $\mathop{\rm
    dist}(\alpha,\Z)\neq\frac12$, then there exists $j_0\in\Z$ and
  $\beta\in \C$ such that
\[
\lambda^{-\sqrt{(\alpha-j_0)^2-a_0}}u(\lambda\cos t,\lambda\sin t)\to
\beta e^{-ij_0t}\quad\text{as }\lambda\to 0^+,
\]
in $C^{1,\tau}(0,2\pi,\C)$ for all $\tau\in(0,1)$;

\medskip\noindent
\item[ii)] if $a_0<\big(\mathop{\rm
    dist}(\alpha,\Z)\big)^2$ and $\mathop{\rm
    dist}(\alpha,\Z)=\frac12$, then there exists $j_0\in\Z$ and
  $\beta_1,\beta_2\in \C$ such that $2\alpha-j_0\in\Z$ and
\[
\lambda^{-\sqrt{(\alpha-j_0)^2-a_0}}u(\lambda\cos t,\lambda\sin t)\to
\beta_1 e^{-ij_0t}+\beta_2 e^{-i(2\alpha-j_0)t}\quad\text{as }\lambda\to 0^+,
\]
in $C^{1,\tau}(0,2\pi,\C)$ for all $\tau\in(0,1)$.
\end{itemize}
The constants $\beta,\beta_1,\beta_2$ can be computed as in
(\ref{eq:38}).  Furthermore, in view of Corollary \ref{cor:holdercont},
if $(\mathop{\rm
  dist}(\alpha,\Z))^2<1+a_0$ then  $u\in C^{0,\gamma}_{\rm loc}(\Omega,\C)$
with $\gamma=\sqrt{(\mathop{\rm
  dist}(\alpha,\Z))^2-a_0}$, whereas
 $u$ is locally Lipschitz continuous in
$\Omega$ if $(\mathop{\rm
    dist}(\alpha,\Z))^2\geq 1+a_0$.

\section{Magnetic Hardy-Sobolev type
  inequalities}\label{sec:magn-hardy-sobol}

This section is devoted to the proof of a weighted electromagnetic
Hardy-Sobolev inequality in dimension $N\geq 3$. We start by observing
that, from Lemma \ref{l:pos} and classical Sobolev's inequality, the
following electromagnetic Hardy-Sobolev inequality holds.
\begin{Proposition}\label{p:magn-hardy-sobol}
 Let $N\geq 3$ and  $a,\A$ satisfying (\ref{eq:reg}), (\ref{eq:transversality}), and (\ref{positivity}).
Then
\[
S(\A,a):=\inf_{u\in \Di\setminus\{0\}}\frac{Q_{\A,a}(u)}
{\left(\int_{\R^N}|u(x)|^{2^*}\,dx\right)^{2/2^*}}>0.
\]
\end{Proposition}
\begin{pf}
  The proof follows from Lemma \ref{l:pos}, part (i) of Lemma
  \ref{l:spaces}, and Sobolev's inequality.
\end{pf}

\noindent We assume $N\geq 3$ and \eqref{pos-mu10a}
so that the number
\begin{equation}\label{eq:sigma}
  \sigma=\sigma(a,N):=-\frac{N-2}2+
  \sqrt{\bigg(\frac{N-2}2\bigg)^{\!\!2}+\mu_1(0,a)}
\end{equation}
is well defined. Let $\phi\in H^1(\SN,\R)$, $\|\phi\|_{L^2(\SN,\R)}=1$,
be the first positive eigenfunction of the eigenvalue problem
\begin{equation*}
-\Delta_{\SN}
\phi(\theta)-a(\theta)\phi(\theta)=\mu_1(0,a)\phi(\theta) \qquad
{\rm in } \ \SN.
\end{equation*}
We recall from \cite[Lemma 2.1]{FMT2} that $\mu_1(0,a)$ is simple and
$\min_{\SN} \phi>0$.
Let
\begin{equation} \label{PHI}
w(x)=|x|^{\sigma} \phi\bigg(\frac{x}{|x|}\bigg) \qquad {\rm for \
all \ } x\in \R^N\setminus \{0\}
\end{equation}
and introduce the weighted space ${\mathcal
  D}^{1,2}_{w}(\R^N,\C)$ as the closure of $C^\infty_{\rm
  c}(\R^N,\C)$ with respect to the norm
\begin{displaymath}
  \|v\|_{{\mathcal D}^{1,2}_{w}(\R^N,\C)}:
  =\bigg(\int_{\R^N} w^2(x)\big|\nabla v(x)\big|^2dx\bigg)^{\!\!1/2}.
\end{displaymath}
We notice that, by the
Caffarelli-Kohn-Nirenberg inequality (see \cite{CKN} and
\cite{CatrinaWang}), $v\in{\mathcal
  D}^{1,2}_{w}(\R^N,\C)$ if and only if  $w v\in\Di$ and
there exists  ${\mathcal C}_w>0$ such that
\begin{equation*} %\label{CAFF}
{\mathcal C}_w \int_{\R^N}w^2(x) \,\frac{|v(x)|^2}{|x|^2}\,dx\leq
\int_{\R^N}w^2(x) |\nabla v(x)|^2dx
\end{equation*}
for every $v\in{\mathcal D}^{1,2}_{w}(\R^N,\C)$.
\begin{Proposition}\label{p:ckn}
  Let $N\geq 3$ and $a,\A$ satisfying (\ref{eq:reg}),
  (\ref{eq:transversality}), (\ref{pos-mu10a}), and let $w$ be the
  function defined in (\ref{PHI}). Then
  \begin{align}
    \label{eq:ckn}
    \int_{\R^N}&{\textstyle{w^2(x)\big|\nabla v(x)+i\,
          \frac{\A({x}/{|x|})}
          {|x|}\,v(x)\big|^2\,dx}}
    \geq S(\A,a)\left(
      \int_{\R^N}w^{2^*}(x)|v(x)|^{2^*}\,dx\right)^{\!\!\frac2{2^*}}
  \end{align}
  for all $v\in{\mathcal D}^{1,2}_{w}(\R^N,\C)$.
\end{Proposition}

\begin{pf}
First of all, one can check by explicit computation that the function
$w$ solves the equation
\begin{equation} \label{phi}
-\Delta w(x)-\frac{a(x/|x|)}{|x|^2}w(x)=0
%(\mu_1(0,a)-\mu_1(\A,a)) \frac{w}{|x|^2}
\quad  \mbox{in} \ \R^N\setminus \{0\}.
\end{equation}
Let $v\in C^{\infty}_c(\R^N\setminus\{0\},\C)\subset {\mathcal
  D}^{1,2}_{w}(\R^N,\C)$ so that
$u(x):=w(x)v(x)\in C^{\infty}_c(\R^N\setminus\{0\},\C)\subset
{\mathcal D}^{1,2}(\R^N,\C)$.  By (\ref{phi}) and integration by parts
we have
\begin{align} \label{BOO}
\int_{\R^N}  & \nabla w(x) \nabla(w(x) |v(x)|^2) \, dx
-\int_{\R^N} \frac{a(x/|x|)}{|x|^2} w^2(x) |v(x)|^2 \, dx=0. 
\end{align}
By a direct computation we infer
\begin{equation}\label{eq:43}
 \nabla w \nabla(w |v|^2) =
  |\nabla w|^2 |v|^2+w\nabla w(\overline v\nabla
  v+v\nabla\overline v)
\end{equation}
and
\begin{align}\label{eq:50}
  \left|\nabla u+i\frac{\A(x/|x|)}{|x|}u\right|^2 &= |\nabla
  w|^2 |v|^2+w\nabla w(\overline v\nabla
  v+v\nabla\overline v)
  +w^2 |\nabla v|^2 \\
  &\notag\quad-2\Im\left(\frac{\A}{|x|}w^2 v\nabla \overline v\right)
  +\frac{|\A|^2}{|x|^2} w^2 |v|^2.
\end{align}
From (\ref{BOO}), (\ref{eq:43}), and (\ref{eq:50}), we obtain that
\begin{align*}
 \int_{\R^N}  \bigg|\nabla
    u(x)+&i\frac{\A(x/|x|)}{|x|}u(x)\bigg|^2 dx
  -\int_{\R^N} \frac{a(x/|x|)}{|x|^2} |u(x)|^2 \, dx \\
  & \notag = \int_{\R^N} w^2(x) |\nabla v(x)|^2 \, dx
  -\int_{\R^N} 2\Im\left(\frac{\A(x/|x|)}{|x|}w^2(x) v(x)\nabla
    \overline v(x)\right) \, dx
  \\
  & \notag \quad  +\int_{\R^N} \frac{|\A(x/|x|)|^2}{|x|^2}\,
  w^2(x) |v(x)|^2 \, dx
  \\
  & \notag = \int_{\R^N} w^2(x) \left|\nabla
    v(x)+i\frac{\A(x/|x|)}{|x|}v(x)\right|^2 \, dx.
\end{align*}
By the above identity and Proposition \ref{p:magn-hardy-sobol}, we
obtain (\ref{eq:ckn}) for any $v\in
C^\infty_c(\R^N\setminus\{0\},\C)$.  By a density argument (see
\cite[Lemma 2.1]{CatrinaWang}), we deduce that (\ref{eq:ckn})
holds for any $v\in{\mathcal D}^{1,2}_{w}(\R^N,\C)$.
\end{pf}

\section{A Brezis-Kato type lemma for $N\geq 3$}\label{sec:bk}

This section is devoted to the proof of a Brezis-Kato type result
in dimension $N\geq 3$.  Let $w$ be the function defined in
(\ref{PHI}).  We define the weighted space $H^1_{w}(\Omega,\C)$ as
the closure of $H^1(\Omega,\C)\cap C^\infty(\Omega,\C)$ with
respect to the norm
\begin{equation} \label{H1As}
  \|v\|_{H^1_{w}(\Omega,\C)}:
  =\bigg(\int_{\Omega} w^2(x) \bigg[|\nabla v(x)|^2
+|v(x)|^2\bigg]\,dx\bigg)^{\!\!1/2},
\end{equation}
and the space ${\mathcal D}^{1,2}_{w}(\Omega,\C)$
 as the closure of $C^\infty_{\rm c}(\Omega,\C)$ with respect to
\begin{align*}
  \|v\|_{{\mathcal D}^{1,2}_{w}(\Omega,\C)} :=
\bigg(\int_{\Omega} w^2(x) |\nabla v(x)|^2\,dx\bigg)^{\!\!1/2}.
\end{align*}
It is easy to verify that $v\in H^1_{w}(\Omega,\C)$ if and only if
$w v\in H^1(\Omega,\C)$. For $N\geq 3$ and any $q\geq 1$, we also
denote as $L^q(w^{2^*},\Omega,\C)$ the weighted $L^q$-space
endowed with the norm
$$
\|v\|_{L^q(w^{2^*},\Omega,\C)}:= \bigg(\int_{\Omega}
w^{2^*}(x)|v(x)|^q\,dx\bigg)^{\!\!1/q},
$$
where $2^*=\frac{2N}{N-2}$ is the critical Sobolev exponent.

\begin{Lemma}\label{l:bk1}
  Let $\Omega\subset\R^N$, $N\geq 3$, be a bounded open set containing
  $0$, (\ref{eq:reg}), (\ref{eq:transversality}) and
  (\ref{pos-mu10a}) hold, and $v\in H^1_{w}(\Omega,\C)\cap
  L^q(w^{2^*},\Omega,\C)$, $q>2$, be a weak solution to
  \begin{align}\label{eq:6}
    -\mathop{\rm div}&(w^2(x)\nabla
    v(x))
    - \frac{2i \frac{\A(x/|x|)}{\phi(x/|x|)}\nabla_{\SN} 
   \phi(\frac{x}{|x|})-|\A(\frac{x}{|x|})|^2
      +i\dive\nolimits_{{\mathbb
          S}^{N-1}}\A(\frac{x}{|x|})}{|x|^2} \, w^2(x)v(x)\\
    &\notag \hskip5cm
    -2i \, w^2(x) \,\frac{\A\big(\frac{x}{|x|}\big)}{|x|}\cdot \nabla
    v(x)= w^{2^*}(x)V(x)v(x),\quad\text{in }\Omega,
  \end{align}
  where $(\Re(V))_+\in L^s(w^{2^*},\Omega,\C)$ for some $s>
  N/2$. Then, for any $\Omega' \Subset \Omega$ such that $0\in\Omega'$, $v\in
  L^{\frac{2^*q}2}(w^{2^*},\Omega',\C)$ and
\begin{multline}\label{eq:66}
  \|v\|_{L^{\frac{2^*q}2}(w^{2^*},\Omega',\C)}\leq
  S(\A,a)^{-\frac1q}  \|v\|_{L^{q}(w^{2^*},\Omega,\C)}
  \bigg(\frac{32}{C(q)} \frac{M^{2-2^*}\big(
\widetilde C(\Omega,\Omega')\big)^{\sigma(2-2^*)}}{(\mathop{\rm dist}
    (\Omega',\partial\Omega))^2}+\frac{2\ell_q}{C(q)} 
\bigg)^{\!\!\frac1q},
\end{multline}
where $C(q):=\min\big\{\frac14, \frac4{q+4}\big\}$,
$M=\min_{\SN} \phi>0$,
\begin{align*}
&  \widetilde C(\Omega,\Omega') =
\begin{cases}
\mathop{\rm
      diam}\Omega&\text{if }\mu_1(0,a)\leq0,\\
\mathop{\rm dist}
    (0,\R^N\setminus\Omega')&\text{if }\mu_1(0,a)>0,
  \end{cases}
\quad\text{and}\\
&  \ell_q=\bigg[\max\bigg\{\frac{8}{S(\A,a)}\,
  \|(\Re(V))_+\|_{L^s(w^{2^*},\Omega,\C)} ^{2s/N},
  \frac{q+4}{2\,S(\A,a)} \|(\Re(V))_+\|_{L^s(w^{2^*},\Omega,\C)}
  ^{2s/N}\bigg\}\bigg]^{\frac{N}{2s-N}}.
\end{align*}
\end{Lemma}

\begin{pf}
  H\"older's inequality and (\ref{eq:ckn}) yield for any $u\in
  {\mathcal D}^{1,2}_{w}(\Omega,\C)$
\begin{align}
  \label{eq:4bk}
  &\qquad\int_{\Omega}w^{2^*}(x)(\Re(V(x)))_+|u(x)|^2\,dx\\
  &\notag \leq\ell_q\!\!\!\!\!\!
  \int\limits_{(\Re(V(x)))_+\leq\ell_q}\!\!\!\!\!\!w^{2^*}(x)
  |u(x)|^2\,dx +\!\!\!\!\!\!  \int\limits_{(\Re(V(x)))_+\geq
    \ell_q}\!\!\!\!\!\!
  w^{(2^*-2)}(x)(\Re(V(x)))_+w^{2}(x)|u(x)|^2\,dx\\
  \notag &\leq \ell_q\int_{\Omega}w^{2^*}(x) |u(x)|^2\,dx
  +\bigg(\int_{\Omega}w^{2^*}(x)|u(x)|^{2^*}dx\bigg)^{\!\!\frac2{2^*}}
  \bigg( \int\limits_{(\Re(V(x)))_+\geq \ell_q}
  w^{2^*}(x)(\Re(V(x)))_+^{\frac N2}\,dx\bigg)^{\!\!\frac2N}\\
  \notag &\leq \frac1{S(\A,a)} \bigg(\int_{\Omega}
  {\textstyle{w^{2}(x) \big|\nabla u(x)+i\,
        \frac{\A({x}/{|x|})}
        {|x|}\,u(x)\big|^2
 \,dx}}\bigg)\times\\
  &\notag\hskip3cm\times \bigg( \int\limits_{(\Re(V(x)))_+\geq \ell_q}
  w^{2^*}(x)(\Re(V(x)))_+^{\frac N2}\,dx\bigg)^{\!\!\frac2N}+
  \ell_q\int_{\Omega}w^{2^*}(x) |u(x)|^2\,dx.
\end{align}
By H\"older's inequality and by the choice of
$\ell_q$ it follows that
\begin{gather}\label{eq:70}
  \int\limits_{(\Re(V(x)))_+\geq \ell_q}
  \!\!\!\!\! w^{2^*}(x)(\Re(V(x)))_+^{\frac N2}\,dx \leq
  \bigg(\int_{\Omega} w^{2^*}(x)(\Re(V(x)))_+^s\,dx
  \bigg)^{\!\!\frac{N}{2s}}\bigg( \int\limits_{(\Re(V(x)))_+\geq
    \ell_q} \!\!\!\!\! w^{2^*}(x)\,dx
  \bigg)^{\!\!\frac{2s-N}{2s}}\\
  \notag\leq
  \bigg(\int_{\Omega}w^{2^*}(x)(\Re(V(x)))_+^s\,dx\bigg)^{\!\!\frac{N}{2s}}\bigg(
  \int\limits_{(\Re(V(x)))_+\geq
    \ell_q}\bigg(\frac{(\Re(V(x)))_+}{\ell_q} \bigg)^{\!\!s}
  w^{2^*}(x)\,dx\bigg)^{\!\!\frac{2s-N}{2s}}\\
  \notag\leq
  \|(\Re(V))_+\|_{L^s(w^{2^*},\Omega,\C)}^{s}\,\ell_q^{-s+\frac
    N2}\leq\min\bigg\{\frac{S(\A,a)}{8}, \frac{2\,
    S(\A,a)}{q+4}\bigg\}^{\!\!\frac N2},
\end{gather}
and hence from (\ref{eq:4bk}) we obtain that for any $u\in
  {\mathcal D}^{1,2}_{w}(\Omega,\C)$
\begin{align}\label{eq:64}
  &\int_{\Omega} w^{2^*}(x)(\Re(V(x)))_+|u(x)|^2\,dx \leq
  \ell_q\int_{\Omega} w^{2^*}(x) |u(x)|^2\,dx\\
  &\notag \quad+ \min\bigg\{\frac{1}{8}, \frac{2}{q+4}\bigg\}
  \bigg(\int_{\Omega} {\textstyle{w^{2}(x)\big|\nabla
        u(x)+i\, \frac{\A({x}/{|x|})}
        {|x|}\,u(x)\big|^2
\,dx}}\bigg).
\end{align}
Let $\eta\in C^\infty_{\rm c}(\Omega,\R)$ be a nonnegative cut-off
function such that
$$
\mathop{\rm supp}(\eta) \Subset \Omega, \quad \eta \equiv 1 \text{ on
} \Omega',\text{ and } |\nabla \eta(x)|\leq \frac{2}{\mathop{\rm dist}
  (\Omega',\partial\Omega)}.
$$
Set $v^n := \min(n,|v|) \in H^1_{w}(\Omega,\C)$.  Let us test
(\ref{eq:6}) with $\eta^2(v^n)^{q-2}\bar v \in{\mathcal
  D}^{1,2}_{w}(\Omega,\C)$ and take the real part.  Observing
that $\Re(\bar v\nabla v)=|v|\nabla|v|$ and using the elementary
inequality $2ab \le 1/2 a^2 + 2 b^2$ and the diamagnetic inequality 
(see Lemma \ref{l:diam}), we thus obtain
\begin{align*}
  &(q-2)\int_{\Omega} w^{2}(x)\eta^2(x)(v^n(x))^{q-2}
  \alchi_{\{y\in\Omega:|v(y)|<n\}}(x)|\nabla|v|(x)|^2\,dx\\
  &\quad\quad+ \int_{\Omega}w^{2}(x) \eta^2(x) (v^n(x))^{q-2}|\nabla
  v(x)|^2\,dx + \int_{\Omega}\frac{|\A(\frac{x}{|x|})|^2}
  {|x|^2}\, w^2(x) \eta^2(x) (v^n(x))^{q-2}|v(x)|^2\,dx\\
  &\quad\quad+2\int_{\Omega} w^2(x) \eta^2(x)
  (v^n(x))^{q-2}\frac{\A(\frac{x}{|x|})}{|x|}
  \cdot\Im(\bar v(x)\nabla v(x))\,dx\\
  &=\!\int_{\Omega}\! w^{2^*}(x)\Re(V(x))\eta^2(x)
  |v(x)|^2(v^n(x))^{q-2}\,dx -2\!\int_{\Omega}\! w^{2}(x)
  \eta(x)(v^n(x))^{q-2}
  |v(x)|\nabla |v|(x)\cdot    \nabla \eta(x)\,dx\\
  &\leq \int_{\Omega}\! w^{2^*}(x)\Re(V(x))\eta^2(x)
  |v(x)|^2(v^n(x))^{q-2}\,dx +2\int_{\Omega}\! w^{2}(x)
  |\nabla\eta(x)|^2(v^n(x))^{q-2}
  |v(x)|^2\,dx\\
  & \quad\quad+\frac12 \int_{\Omega}\! w^{2}(x) \eta^2(x)(v^n(x))^{q-2}
  \left|\nabla v(x)+i\frac{\A(\frac{x}{|x|})}{|x|} v(x)\right|^2\,dx
 \end{align*}
and hence
\begin{gather}\label{eq:19bk}
  (q-2)\int_{\Omega}w^{2}(x)\eta^2(x)(v^n(x))^{q-2}
  \alchi_{\{y\in\Omega:|v(y)|<n\}}(x)|\nabla|v|(x)|^2\,dx\\
  \notag\ +\frac12 \int_{\Omega} w^{2}(x)\eta^2(x) (
  v^n(x))^{q-2}\left|\nabla v(x)+i\frac{\A(\frac{x}{|x|})}{|x|}
    v(x)\right|^2\,dx \\ \notag\leq \int_{\Omega}\!
  w^{2^*}(x)\Re(V(x))\eta^2(x) |v(x)|^2(v^n(x))^{q-2}\,dx
  +2\int_{\Omega}\! w^{2}(x) |\nabla\eta(x)|^2(v^n(x))^{q-2}
  |v(x)|^2\,dx.
 \end{gather}
 Furthermore, by diamagnetic inequality (see Lemma \ref{l:diam}) we
 have that 
\begin{align}\label{eq:7}
&   \left|\nabla \big((v^n)^{\frac q2-1}v\eta\big)+i\,
    \frac{\A(\frac{x}{|x|})} {|x|}\,(v^n)^{\frac
      q2-1}v\eta\right|^2\\
  &\quad\notag= \left|\nabla \big((v^n)^{\frac q2-1}v\eta\big)\right|^2 +2
  \frac{\A(\frac{x}{|x|})} {|x|}\,\eta^2(v^n)^{q-2} \Im(\bar v\nabla
  v) +\frac{|\A(\frac{x}{|x|})|^2}{|x|^2}
  (v^n)^{q-2}\eta^2|v|^2\\
  &\quad\notag\leq\frac{(q+4)(q-2)}4(v^n)^{q-2}\eta^2|\nabla v^n|^2
+2(v^n)^{q-2}\eta^2 \left|\nabla
v+i\frac{\A(\frac{x}{|x|})}{|x|}v \right|^2\\
  &\quad\quad\notag
+\frac{q+2}2(v^n)^{q-2} |v|^2|\nabla \eta|^2.
\end{align}
Letting $C(q):=\min\big\{\frac14, \frac4{q+4}\big\}$, from
(\ref{eq:19bk}) and (\ref{eq:7}) we obtain
\begin{gather}\label{eq:22bk}
  C(q)\int_{\Omega}w^{2}(x)  \bigg|\nabla
      \big((v^n)^{\frac q2-1}v\eta\big)(x)+i\, \frac{\A(\frac{x}{|x|})}
      {|x|}\,(v^n(x))^{\frac
        q2-1}v(x)\eta(x)\bigg|^2 \,dx\\
\notag\leq \int_{\Omega}\!w^{2^*}(x)\Re(V(x))\eta^2(x)
   |v(x)|^2(v^n(x))^{q-2}\,dx\\
\notag +2\int_{\Omega}\! w^{2}(x)
(v^n(x))^{q-2} |v(x)|^2|\nabla \eta(x)|^2\,dx
+C(q)\frac{q+2}2\int_{\Omega}\!w^{2}(x)
(v^n(x))^{q-2} |v(x)|^2|\nabla \eta(x)|^2\,dx.
\end{gather}
Estimate (\ref{eq:64}) applied to $\eta(v^n)^{\frac q2-1}v$ gives
\begin{gather}
\label{eq:5bk}
\int_{\Omega}w^{2^*}(x)(\Re(V(x)))_+|\eta(x)(v^n(x))^{\frac
  q2-1}v(x)|^2\,dx \leq
\ell_q\int_{\Omega}w^{2^*}(x) |\eta(x)(v^n(x))^{\frac q2-1}v(x)|^2\,dx\\
\notag + \min\bigg\{\frac{1}{8}, \frac{2}{q+4}\bigg\}
\bigg(\int_{\Omega} w^{2}(x) \bigg|\nabla
      (\eta(v^n)^{\frac q2-1}v)(x)+i\, \frac{\A({x}/{|x|})}
      {|x|}\,\eta(x)(v^n(x))^{\frac
        q2-1}v(x)\bigg|^2
 \,dx\bigg).
\end{gather}
Using  (\ref{eq:5bk}) to estimate the term with $V$ in (\ref{eq:22bk}),
 (\ref{eq:ckn}) yields
\begin{gather*}
  \bigg(\int_{\Omega} w^{2^*}(x) |v^n(x)|^{(\frac
    q2-1)2^*}|v(x)|^{2^*}\eta^{2^*}(x)\,dx\bigg)^{\frac 2{2^*}} \leq
  \frac{2\ell_q}{C(q)S(\A,a)}\int_{\Omega} w^{2^*}(x)
  \eta^2(x)|v^n(x)|^{q-2}|v(x)|^2\,dx\\
  +\frac{4+C(q)(q+2)}{C(q)S(\A,a)} \int_{\Omega}
  w^{2}(x)|v^n(x)|^{q-2}|v(x)|^2|\nabla
  \eta(x)|^2\,dx\\
\leq \frac{2\ell_q}{C(q)S(\A,a)}\int_{\Omega} w^{2^*}(x)
  \eta^2(x)|v^n(x)|^{q-2}|v(x)|^2\,dx\\
  \qquad+\frac{8}{C(q)S(\A,a)} \int_{\Omega}
  w^{2}(x)|v^n(x)|^{q-2}|v(x)|^2|\nabla
  \eta(x)|^2\,dx.
\end{gather*}
Letting $n\to\infty$ in the above inequality, (\ref{eq:66}) follows.
\end{pf}

\begin{remark} \label{remAL}
It is possible to extend the result of Lemma \ref{l:bk1} also to the case
$$(\Re(V))_+\in L^{N/2}(w^{2^*},\Omega,\C)$$ and obtain
estimate (\ref{eq:66}). Indeed, by the previous summability assumption
on $(\Re(V))_+$, it is possible to find $\ell_q$
such that
$$
\int\limits_{(\Re(V(x)))_+\geq \ell_q}
w^{2^*}(x)(\Re(V(x)))_+^{\frac N2}\,dx
\leq
\min\bigg\{\frac{S(\A,a)}{8},
  \frac{2\, S(\A,a)}{q+4}\bigg\}^{\!\!\frac N2}.
$$
But we have not a control on the constant $\ell_q$ in terms of $q$ as
in Lemma \ref{l:bk1} since it is not possible to apply H\"older's
inequality in (\ref{eq:70}) when $s=N/2$.  The rest of the proof in
the case $s=N/2$ coincide with the proof of Lemma \ref{l:bk1}.
\end{remark}

The previous lemma allows starting a 
Brezis-Kato type iteration.
\begin{Theorem}\label{t:bk}
 Let $\Omega\subset\R^N$, $N\geq 3$, be a bounded open set containing $0$,
(\ref{eq:reg}), (\ref{eq:transversality}) and (\ref{pos-mu10a}) hold.

\begin{itemize}
\item[i)] If $V$ is such such that $(\Re(V))_+\in
  L^s(w^{2^*},\Omega,\C)$ for some $s> N/2$, then, for any
  $\Omega' \Subset \Omega$, there exists a positive constant
  $$
  C_\infty=C_\infty\big(N,\A,a,\|(\Re(V))_+\|_{L^s(w^{2^*},
\Omega,\C)},\mathop{\rm
    dist} (\Omega',\partial\Omega),\widetilde C(\Omega,\Omega') \big)
  $$
  depending only on $N$, $\A$, $a$,
  $\|(\Re(V))_+\|_{L^s(w^{2^*},\Omega,\C)}$, $\mathop{\rm
    dist} (\Omega',\partial\Omega)$, and $\widetilde
  C(\Omega,\Omega')$, such that for any weak
  $H^1(\Omega,\C)$-solution $u$ to
\begin{equation}\label{eq:67}
  {\mathcal L}_{\A,a}u(x)=w^{2^*-2}(x)V(x)u(x),
\quad \text{in }\Omega,
\end{equation}
there holds
$|x|^{-\sigma}{u}\in L^{\infty}(\Omega',\C)$ and
\begin{equation*}
  \big\||x|^{-\sigma}{u}\big\|_{L^{\infty}(\Omega',\C)}\leq
C_\infty\,\|u\|_{L^{2^*}(\Omega,\C)}.
\end{equation*}
\item[ii)] If $V$ is such
  that $(\Re(V))_+\in L^{N/2}(w^{2^*},\Omega,\C)$, then, for
  any $\Omega' \Subset \Omega$ and for any $s\geq 1$, there
  exists a positive constant
  $$
C_s=C_s\big(N,\A,a,(\Re(V))_+,s,\mathop{\rm
    dist} (\Omega',\partial\Omega),\widetilde C(\Omega,\Omega') \big)
$$
  depending only on $N$, $\A$, $a$,
  $(\Re(V))_+$, $s$, $\mathop{\rm dist}
  (\Omega',\partial\Omega)$, and $\widetilde C(\Omega,\Omega')$, such that
for any weak
  $H^1(\Omega,\C)$-solution  $u$ to (\ref{eq:67}) in $\Omega$
there holds
$|x|^{-\sigma}{u}\in L^{s}(w^{2^*},\Omega',\C)$ and
\begin{equation*}
  \big\||x|^{-\sigma}{u}\big\|_{L^{s}(w^{2^*},\Omega',\C)}\leq
C_s\,\|u\|_{L^{2^*}(\Omega,\C)}.
\end{equation*}
\end{itemize}
\end{Theorem}

\begin{pf}
i) Let $u$ be a weak
  $H^1(\Omega,\C)$-solution to (\ref{eq:67}).  It is easy to verify that
   $v:=w^{-1}u$ belongs to $H^1_{w}(\Omega,\C)$ and is a
weak solution to (\ref{eq:6}). Let $R>0$ be such that
$$
\Omega'\Subset
\Omega'+B(0,2R)\Subset \Omega.
$$
Using Lemma \ref{l:bk1} in $\Omega_1:=\Omega'+B(0,R(2-r_1))\Subset
\Omega'+B(0,2R)$, $r_1=1$, with $q=q_1=2^*$, we infer that $v\in
L^{\frac{(2^*)^2}2}(w^{2^*},\Omega_1,\C)$ and the following
estimate holds
\begin{multline*}
  \|v\|_{L^{\frac{(2^*)^2}2}(w^{2^*},\Omega_1,\C)}\leq
S(\A,a)^{-\frac1{q_1}}  \|v\|_{L^{2^*}(w^{2^*},\Omega,\C)}
\bigg(
\frac{32}{C(q_1)} \frac{M^{2-2^*}\big(
\widetilde C(\Omega,\Omega')\big)^{\sigma(2-2^*)}}{(Rr_1)^2}
  +\frac{2\ell_{q_1}}{C(q_1)} 
\bigg)^{\!\!\frac1{q_1}}.
\end{multline*}
Using again Lemma \ref{l:bk1} in $\Omega_2:=\Omega'+B(0,R(2-r_1-r_2))\Subset
\Omega_1$, $r_2=\frac14$, with $q=q_2=(2^*)^2/2$, we infer that
$v\in
L^{\frac{(2^*)^3}4}(w^{2^*},\Omega_2,\C)$ and
\begin{align*}
  &\|v\|_{L^{\frac{(2^*)^3}4}(w^{2^*},\Omega_2,\C)}\leq
  S(\A,a)^{-\frac1{q_2}}
  \bigg(
\frac{32}{C(q_2)} \frac{M^{2-2^*}\big(
\widetilde C(\Omega,\Omega')\big)^{\sigma(2-2^*)}}{(Rr_2)^2}
  +\frac{2\ell_{q_2}}{C(q_2)} 
\bigg)^{\!\!\frac1{q_2}}
\|v\|_{L^{q_2}(w^{2^*},\Omega_1,\C)}
  \\
  & \leq
  S(\A,a)^{-(\frac1{q_1}+\frac1{q_2})}
\times
 \bigg(
\frac{32}{C(q_1)} \frac{M^{2-2^*}\big(
\widetilde C(\Omega,\Omega')\big)^{\sigma(2-2^*)}}{(Rr_1)^2}
  +\frac{2\ell_{q_1}}{C(q_1)} 
\bigg)^{\!\!\frac1{q_1}}
\times \\
& \hskip4cm\times
 \bigg(
\frac{32}{C(q_2)} \frac{M^{2-2^*}\big(
\widetilde C(\Omega,\Omega')\big)^{\sigma(2-2^*)}}{(Rr_2)^2}
  +\frac{2\ell_{q_2}}{C(q_2)} 
\bigg)^{\!\!\frac1{q_2}}
\|v\|_{L^{2^*}(w^{2^*},\Omega,\C)}.
\end{align*}
Setting, for any  $n\in\N$, $n\geq 1$,
$$
q_n=2\bigg(\frac{2^*}2\bigg)^{\!n},\quad
\Omega_n:=\Omega'+B\bigg(0,R\bigg(2-\sum_{k=1}^n r_k\bigg)\bigg),\quad\text{and}
\quad
r_n=\frac1{n^2},
$$
and using iteratively Lemma \ref{l:bk1}, we obtain that, for any
$n\in\N$, $n\geq 1$,
\begin{multline}\label{eq:68}
  \|v\|_{L^{q_{n+1}}(w^{2^*},\Omega',\C)}\leq
 \|v\|_{L^{q_{n+1}}(w^{2^*},\Omega_n,\C)}
\leq  \|v\|_{L^{2^*}(w^{2^*},\Omega,\C)}
   (S(\A,a))^{-\sum\limits_{k=1}^n\frac1{q_k}}\times\\
\times \prod_{k=1}^n
 \bigg(
\frac{32}{C(q_k)} \frac{M^{2-2^*}\big(
\widetilde C(\Omega,\Omega')\big)^{\sigma(2-2^*)}}{(Rr_k)^2}
  +\frac{2\ell_{q_k}}{C(q_k)} 
\bigg)
^{\!\!\frac1{q_k}}.
\end{multline}
We notice that
\begin{align*}
\prod_{k=1}^n 
 \bigg(
\frac{32}{C(q_k)} \frac{M^{2-2^*}\big(
\widetilde C(\Omega,\Omega')\big)^{\sigma(2-2^*)}}{(Rr_k)^2}
  +\frac{2\ell_{q_k}}{C(q_k)} 
\bigg)
^{\!\!\frac1{q_k}}=\exp\bigg[\sum_{k=1}^n b_k\bigg]
\end{align*}
where
$$
b_k= \frac1{q_k}\log  \bigg(
\frac{32}{C(q_k)} \frac{M^{2-2^*}\big(
\widetilde C(\Omega,\Omega')\big)^{\sigma(2-2^*)}}{(Rr_k)^2}
  +\frac{2\ell_{q_k}}{C(q_k)} 
\bigg)
,
$$
and, for some constant
$C=C\big(N,\A,a,\|(\Re(V))_+\|_{L^s(w^{2^*},\Omega,\C)},\mathop{\rm
  dist} (\Omega',\partial\Omega),\widetilde C(\Omega,\Omega') \big)$,
$$
b_k\sim
\frac12\bigg(\frac{2}{2^*}\bigg)^{\!\!k}\log\left[C
\bigg(2\bigg(\frac{2^*}2\bigg)^{\!\!k}\bigg)^{\!\frac{2s}{2s-N}}\right]
\quad\text{as }k\to+\infty.
$$
Hence $\sum_{n=1}^\infty b_n$ converges to some positive sum depending
only on
$N$, $\A$, $a$, $\|(\Re(V))_+\|_{L^s(w^{2^*},\Omega,\C)}$,
$\mathop{\rm
  dist} (\Omega',\partial\Omega)$, $\widetilde C(\Omega,\Omega')$,
 hence
$$
\lim_{n\to+\infty}
(S(\A,a))^{-\sum\limits_{k=1}^n\frac1{q_k}}\prod_{k=1}^n
\bigg(
\frac{32}{C(q_k)} \frac{M^{2-2^*}\big(
\widetilde C(\Omega,\Omega')\big)^{\sigma(2-2^*)}}{(Rr_k)^2}
  +\frac{2\ell_{q_k}}{C(q_k)} 
\bigg)^{\!\!\frac1{q_k}}
$$
is finite and depends only on $N$, $\A$, $a$,
$\|(\Re(V))_+\|_{L^s(w^{2^*},\Omega,\C)}$, $\mathop{\rm
  dist} (\Omega',\partial\Omega)$, $\widetilde C(\Omega,\Omega')$.
Hence, from (\ref{eq:68}), we deduce that there exists a positive
constant $C$ depending only on
$\|(\Re(V))_+\|_{L^s(w^{2^*},\Omega,\C)}$, $N$, $\A$, $a$,
$\mathop{\rm dist} (\Omega',\partial\Omega)$, $\widetilde
C(\Omega,\Omega')$, such that
\begin{equation*}
  \|v\|_{L^{q_{n+1}}(w^{2^*},\Omega',\C)}\leq
C\,\|v\|_{L^{2^*}(w^{2^*},\Omega,\C)}\quad\text{for all }n\in\N.
\end{equation*}
Letting $n\to+\infty$ we deduce that $|v|$ is essentially bounded in
$\Omega'$ with respect to the measure $w^{2^*}\!dx$ and
\begin{equation*}
 \|v\|_{L^{\infty}(w^{2^*},\Omega',\C)}\leq
  C\,\|v\|_{L^{2^*}(w^{2^*},\Omega,\C)}=C\,
\|u\|_{L^{2^*}(\Omega,\C)},
\end{equation*}
where $\|v\|_{L^{\infty}(w^{2^*},\Omega',\C)}$ denotes the
essential supremum of $v$ with respect to the measure
$w^{2^*}\!dx$. Since $w^{2^*}\!dx$ is absolutely
continuous with respect to the Lebesgue measure and vice versa, there holds
$ \|v\|_{L^{\infty}(w^{2^*},\Omega',\C)}=\|v\|_{L^{\infty}(\Omega',\C)}$, hence
 $v\in L^{\infty}(\Omega',\C)$ and
\begin{equation*}
  \|v\|_{L^{\infty}(\Omega',\C)}\leq
  C\,
\|u\|_{L^{2^*}(\Omega,\C)},
\end{equation*}
thus completing the proof of part i). We recall that for any $x\in \Omega\setminus\{0\}$ we have 
$$
|x|^{-\sigma} u(x)=w^{-1}(x) \phi(x/|x|) u(x)=\phi(x/|x|)v(x)\leq (\max_{\SN} \phi) v(x).
$$ 

ii) Since $u\in H^1(\Omega,\C)$ is a weak solution to (\ref{eq:67})
then $v:=w^{-1}u\in H^1_{w}(\Omega,\C)$ is a weak solution
of (\ref{eq:6}). Using Remark \ref{remAL} and the iterative scheme
used to prove part i), for any $1\leq s<\infty$, after a finite number
of iterations we arrive to $v\in L^s(w^{2^*},\Omega',\C)$ and
$$
\|v\|_{L^s(w^{2^*},\Omega',\C)}
\leq C_s \|v\|_{L^{2^*}(w^{2^*},\Omega,\C)}.
$$
This completes the proof.
\end{pf}

\noindent Applying Theorem \ref{t:bk}
to the nonlinear equation (\ref{nonlin0}), we can obtain a pointwise
estimate for solutions to (\ref{nonlin0}).

\begin{Theorem} \label{nonlinear}

 Let $\Omega\subset\R^N$, $N\geq 3$, be a bounded open set containing $0$,
(\ref{eq:reg}), (\ref{eq:transversality}) and (\ref{pos-mu10a}) hold.
Let $u$ be a weak $H^1(\Omega,\C)$-solution of (\ref{nonlin0}) with
$f(x,u)$ satisfying (\ref{subcrit_0}).
Then for any $\Omega'\Subset \Omega$,
\begin{equation}\label{eq:41}
|x|^{-\sigma}u\in L^\infty(\Omega',{\mathbb C}).
\end{equation}

\end{Theorem}

\begin{pf} If we put
$$
V(x):=\left\{
\begin{array}{ll}
w^{2-2^*} \frac{f(x,u(x))}{u(x)}, & \quad {\rm if} \ u(x)\neq 0, \\
0, & \quad {\rm if} \ u(x)=0,
\end{array}
\right.
$$
then, by (\ref{subcrit_0}) and the Sobolev embedding
$H^1(\Omega,\C)\subset L^{2^*}(\Omega,\C)$,
we have that $V\in L^{N/2}(w^{2^*},\Omega,\C)$ and  $u$ weakly solves
$$
\mathcal L_{\A,a} u(x)=w^{2^*-2}V(x) u(x) \quad \text{in }  \Omega.
$$
From part ii) of Theorem \ref{t:bk}, it follows that
$|x|^{-\sigma}u\in L^s(w^{2^*},\Omega',\C)$ for any
$\Omega'\Subset \Omega$ and for any $s\geq 1$.  Fix now $s_0=N/2+\e_0$
with $0<\e_0<\frac{N(N-2)}{4|\sigma|}$.  By (\ref{subcrit_0}) we
easily deduce that $V\in L^{s_0}(w^{2^*},\Omega',\C)$.  The
proof of the theorem follows now by part i) of Theorem \ref{t:bk}.
\end{pf}

The a-priori estimate of solutions to the nonlinear problem obtained
above, allows deducing  Theorem \ref{t:asym-nonlin} from
 Theorem \ref{t:asym}.

\medskip

\begin{pfn}{Theorem \ref{t:asym-nonlin} for $N\geq 3$}
Note that all the assumptions of Theorem \ref{nonlinear}
are verified and hence
\begin{equation} \label{sopra}
|u(x)|=O(|x|^\sigma) \quad \text{as } |x|\to 0,
\end{equation}
where $\sigma>-\frac{N-2}{2}$ is defined by (\ref{eq:sigma}).
Therefore, by (\ref{subcrit_0}) and (\ref{sopra}),
\begin{equation*}
\left|
\frac{f(x,u)}{u} \right|
\leq
{\rm const\,} \left(1+|x|^{-2+\frac4{N-2}
 \sqrt{\left(\frac{N-2}2\right)^{2}+\mu_1(0,a)}}\right)
\end{equation*}
for some constant ${\rm const\,}>0$.
Hence, the function
\begin{equation*}
h(x):=\left\{
\begin{array}{ll}
\frac{f(x,u(x))}{u(x)} & \quad {\rm if} \ u(x)\neq 0 \\
0 & \quad {\rm if} \ u(x)=0
\end{array}
\right.
\end{equation*}
satisfies $h(x)=O(|x|^{-2+\e})$ as $|x|\rightarrow 0^+$ for some
$\e>0$.  On the other hand, by Remark \ref{remA0} we also have $u\in
L^{\infty}_{\rm loc}(\Omega\setminus \{0\})$ and in turn by
(\ref{subcrit_0}), $h\in L^{\infty}_{\rm loc}(\Omega\setminus \{0\})$.
This shows that all the assumptions of Theorem \ref{t:asym} are
satisfied and the proof of Theorem \ref{t:asym-nonlin} follows in the
case $N\geq 3$. The proof of Theorem \ref{t:asym-nonlin} in the case
$N=2$ is postponed to section \ref{sec:brezis-kato-type}.
\end{pfn}

\begin{pfn}{Theorem \ref{asy_infinitynonlin} for $N\geq 3$}
It follows from Theorems \ref{asy_infinity} and \ref{t:asym-nonlin}
by the use of the Kelvin transform.
\end{pfn}

Since the proof of the pointwise a-priori estimate (\ref{eq:41}) (and
then of Theorems \ref{t:asym-nonlin} and \ref{asy_infinitynonlin}) in
dimension $N=2$ originates from a different inequality than
(\ref{eq:ckn}) and requires a little bit different notation, we devote
the next section to a sketched description of the modifications to be
made in the above argument to treat the case $N=2$.

\section{A Brezis-Kato type lemma in dimension
  $N=2$}\label{sec:brezis-kato-type}

Similarly to section \ref{sec:bk}, for $N=2$ we define the spaces $\mathcal
D^{1,2}_{*}(\Omega,\C)$ and $\mathcal D^{1,2}_{*,w}(\Omega,\C)$ as
the completion of $C^\infty_c(\Omega\setminus \{0\},\C)$ respectively
with the norms
\begin{align*}
  \|u\|_{{\mathcal D}^{1,2}_{*}(\Omega,\C)} :=
\bigg(\int_{\Omega} \left(
|\nabla u(x)|^2+\frac{|u(x)|^2}{|x|^2} \right)dx \bigg)^{\!\!1/2}
\end{align*}
and
\begin{align*}
  \|v\|_{{\mathcal D}^{1,2}_{*,w}(\Omega,\C)} :=
\bigg(\int_{\Omega} w^{2}\left(
|\nabla u(x)|^2+\frac{|u(x)|^2}{|x|^2} \right) dx \bigg)^{\!\!1/2}
\end{align*}
where $\Omega\subset \R^2$ is a bounded domain containing the origin
and $w$ is defined by (\ref{PHI}). We observe that
the space $\mathcal D^{1,2}_*(\Omega,\C)$ is smaller than
$H^1_0(\Omega,\C)$.  Moreover, it easy to verify that
$v\in
  \mathcal D^{1,2}_{*,w}(\Omega,\C)$ if and only if $w
  v\in \mathcal D^{1,2}_*(\Omega,\C)$.
Similarly, we define the space
$H^1_{*,w}(\Omega,\C)$ as
the completion of $\{v\in
H^1(\Omega,\C)\cap C^\infty(\Omega,\C):v\text{ vanishes in a
  neighborhood of }0\}$ with respect to the norm
$$
\|v\|_{H^1_{*,w}(\Omega,\C)}
:=
\bigg(\int_{\Omega} w^{2}\bigg[|\nabla v(x)|^2
+\frac{|v(x)|^2}{|x|^2}+|v(x)|^2\bigg]\,dx\bigg)^{\!\!1/2}.
$$
The following weighted Poincar\'e-Sobolev
inequality holds.

\begin{Proposition}\label{p:weighted_poincare_sobolev}
  Let $N=2$ and $a,\A$ satisfying (\ref{eq:reg}),
  (\ref{eq:transversality}) and (\ref{pos-mu10a}).  Then, for any
  $1\leq p<\infty$,
  \begin{equation} \label{firstp} S(\A,a,p,\Omega)= \inf_{u\in
      \mathcal D^{1,2}_*(\Omega,\C)\setminus \{0\}} \frac{
      {\displaystyle{\int_\Omega}}
      \left[\left|\left(\nabla+i\frac{\A}{|x|}\right)u(x)\right|^2
        -\frac{a(x/|x|)}{|x|^2}|u(x)|^2 \right] \, dx}
    {\left({\displaystyle{\int_\Omega}} |u(x)|^p \, dx\right)^{\!\!2/p}}>0.
\end{equation}
Moreover
\begin{align}\label{secondp}
  \int_\Omega w^{2} \left|\nabla
      v(x)+i\tfrac{\A(x/|x|)}{|x|}v(x)\right|^2
  \, dx 
  \geq S(\A,a,p,\Omega) \left(\int_\Omega
    w^{p}|v(x)|^p \, dx\right)^{2/p}
\end{align}
for all $v\in \mathcal D^{1,2}_{*,\sigma}(\Omega,\C)$.
\end{Proposition}

\begin{pf}
  Inequality (\ref{firstp}) follows by Lemma \ref{l:pos} and classical Poincar\'e-Sobolev inequality.
To obtain  the second part of the statement, by density it is sufficient to
  prove inequality (\ref{secondp}) for functions $v\in
  C^\infty_c(\Omega\setminus\{0\},\C)$ as one can easily do by following the
  same procedure developed in the proof of Proposition \ref{p:ckn}.
\end{pf}

\begin{remark}\label{r:costomega}
  We notice that the constant in (\ref{secondp}) depends on the domain
  $\Omega$, unlike the constant appearing in (\ref{eq:ckn}) in the
  case $N=3$ and $p=2^*$.  Moreover $S(\A,a,p,\Omega)$ is decreasing
  with respect to $\Omega$, i.e. if $\Omega_1\subset\Omega_2$
  then $S(\A,a,p,\Omega_1)\geq S(\A,a,p,\Omega_2)$.
\end{remark}

We are now ready to prove the following
2-dimensional version of Lemma \ref{l:bk1}.

\begin{Lemma}\label{l:bk12d}
  Let $\Omega\subset\R^2$ be a bounded open set containing $0$,
  (\ref{eq:reg}), (\ref{eq:transversality}), (\ref{pos-mu10a})
  hold, and, for some $p>2$ and $q>2$, let $v\in
  H^1_{*,w}(\Omega,\C)\cap L^q(w^{p},\Omega,\C)$ be a
  weak solution to
  \begin{align*}
    -\mathop{\rm div}&(w^{2}(x)\nabla
    v(x))
- \frac{2i \frac{\A(x/|x|)}{\phi(x/|x|)}\nabla_{\SN} 
   \phi(\frac{x}{|x|})-|\A(\frac{x}{|x|})|^2
      +i\dive\nolimits_{{\mathbb
          S}^{N-1}}\A(\frac{x}{|x|})}{|x|^2} \, w^2(x)v(x)\\
    &\notag \hskip5cm
    -2i \, w^2(x) \,\frac{\A\big(\frac{x}{|x|}\big)}{|x|}\cdot \nabla
    v(x)=w^{p}(x)V(x)v(x),\quad\text{in }\Omega,
  \end{align*}
where $(\Re(V))_+\in L^s(w^{p},\Omega,\C)$ for some $s>\frac{p}{p-2}$.
Then, for any $\Omega' \Subset \Omega$ such that $0\in\Omega'$, $v\in
  L^{\frac{p\,q}2}(w^{p},\Omega',\C)$ and
\begin{multline*}
  \|v\|_{L^{\frac{p\, q}2}(w^{p},\Omega',\C)}\leq
  S(\A,a,p,\Omega)^{-\frac1q}  \|v\|_{L^{q}(w^{p},\Omega,\C)}
  \times 
\bigg(\frac{32}{C(q)} \frac{M^{2-p}\big(
\widetilde C(\Omega,\Omega')\big)^{\sigma(2-p)}}{(\mathop{\rm dist}
    (\Omega',\partial\Omega))^2}+\frac{2\ell_q}{C(q)} 
\bigg)^{\!\!\frac1q},
\end{multline*}
where $C(q):=\min\big\{\frac14, \frac4{q+4}\big\}$, 
$  \widetilde C(\Omega,\Omega') =
\mathop{\rm dist}
    (0,\R^N\setminus\Omega')$, $M=\min_{\SN} \phi>0$ and
\[
\ell_q=\bigg[\max\bigg\{\frac{8\,\|(\Re(V))_+\|_{L^s(w^{p},\Omega,\C)}
  ^{s(p-2)/p}}{S(\A,a,p,\Omega)}, \frac{q+4}{2\,S(\A,a,p,\Omega)}
\|(\Re(V))_+\|_{L^s(w^{p},\Omega,\C)}
^{s(p-2)/p}\bigg\}\bigg]^{\frac{p}{s(p-2)-p}}.
\]
\end{Lemma}
\begin{pf}
  The proof  may be obtained proceeding as in the proof
  of Lemma \ref{l:bk1} and using (\ref{secondp}) in place of
  (\ref{eq:ckn}).
\end{pf}

The counterpart in dimension $N=2$ of Theorem \ref{t:bk} is the
following Brezis-Kato type result.

\begin{Theorem}\label{t:bk2d}
 Let $\Omega\subset\R^2$ be a bounded open set containing $0$,
(\ref{eq:reg}), (\ref{eq:transversality}), (\ref{pos-mu10a}) hold,
and let $p>2$.
\begin{itemize}
\item[i)]
If  $V$ is such that $(\Re(V))_+\in
  L^s(w^{p},\Omega,\C)$ for some $s>\frac{p}{p-2}$, then,
  for any $\Omega' \Subset \Omega$, there exists a positive constant
  $$
  C_{\infty,2}=C_{\infty,2}\big(\Omega,p,
  \A,a,\|(\Re(V))_+\|_{L^s(w^{p},\Omega,\C)},\mathop{\rm
    dist} (\Omega',\partial\Omega),\widetilde C(\Omega,\Omega') \big)
$$
depending only on $\Omega$, $p$, $\A$, $a$,
$\|(\Re(V))_+\|_{L^s(w^{p},\Omega,\C)}$, $\mathop{\rm dist}
(\Omega',\partial\Omega)$, and $\widetilde C(\Omega,\Omega')$, such
that for any weak $H^1_{*}(\Omega,\C)$-solution $u$ to
\begin{equation}\label{eq:672d}
  {\mathcal L}_{\A,a}u(x)=w^{p-2}V(x)u(x),
\quad \text{in }\Omega,
\end{equation}
there holds
$|x|^{-\sigma}{u}\in L^{\infty}(\Omega',\C)$ and
\begin{equation*}
  \big\||x|^{-\sigma}{u}\big\|_{L^{\infty}(\Omega',\C)}\leq
C_{\infty,2}\,\|u\|_{L^{p}(\Omega,\C)}.
\end{equation*}
\item[ii)]  If $V$ is such
  that $(\Re(V))_+\in L^{\frac p{p-2}}(w^{p},\Omega,\C)$, then, for
  any $\Omega' \Subset \Omega$ and for any $1\leq s<\infty$ there
  exists a positive constant
  $$
C_{s,2}=C_{s,2}\big(\Omega,p,
  \A,a,(\Re(V))_+,\mathop{\rm
    dist} (\Omega',\partial\Omega),\widetilde C(\Omega,\Omega') \big)
$$
depending only on $\Omega$, $p$, $\A$, $a$,
$(\Re(V))_+$, $\mathop{\rm dist}
(\Omega',\partial\Omega)$, and $\widetilde C(\Omega,\Omega')$,  such that
for any weak
  $H^1_{*}(\Omega,\C)$-solution  $u$ to (\ref{eq:672d}) in $\Omega$
there holds
$|x|^{-\sigma}{u}\in L^{s}(w^{p},\Omega',\C)$ and
\begin{equation*}
  \big\||x|^{-\sigma}{u}\big\|_{L^{s}(w^{p},\Omega',\C)}\leq
C_{s,2}\,\|u\|_{L^{p}(\Omega,\C)}.
\end{equation*}
\end{itemize}
\end{Theorem}

\begin{pf}
  This theorem can be proved by iterating the estimate proved in Lemma
  \ref{l:bk12d} and following the same scheme as  in the proof of
  Theorem \ref{t:bk}. We notice that the constants $S(\A,a,p,\Omega_i)$
appearing at each step (at a negative power) can be uniformly controlled
with  $S(\A,a,p,\Omega)$ in  view of Remark \ref{r:costomega}.
\end{pf}

\noindent From the above analysis, the proofs of Theorems
\ref{t:asym-nonlin} and \ref{asy_infinitynonlin} in dimension $N=2$
follow.

\medskip\noindent
\begin{pfn}{Theorem \ref{t:asym-nonlin} for $N=2$}
  Arguing as in the proof of Theorem \ref{nonlinear}, from Theorem
  \ref{t:bk2d} we deduce that $|u(x)|=O(|x|^\sigma)$ as $|x|\to 0$.
  In particular, from (\ref{subcrit_0}), the function
  $\frac{f(x,u(x))}{u(x)}\alchi_{\{x:u(x)\neq 0\}}$ is bounded. The
  conclusion then follows from Theorem \ref{t:asym}.
\end{pfn}

\begin{pfn}{Theorem \ref{asy_infinitynonlin} for $N=2$}
  As in dimension $N\geq 3$, it follows from Theorems
  \ref{asy_infinity} and \ref{t:asym-nonlin} by the use of the Kelvin
  transform.
\end{pfn}

\appendix
 \section*{Appendix}
 \setcounter{section}{1}
 \setcounter{Theorem}{0}

We recall the following well known result proved in \cite{LL}.
\begin{Lemma}\label{l:diam} {\bf (Diamagnetic inequality)}
  Let $N\geq 2$.  If  $u\in \Distar$ then
$$
|\nabla |u|(x)|\leq
\left|\nabla u(x)+i\frac{\A(x/|x|)}{|x|}u(x)\right| \qquad {\rm for \
  a.e.} \ x\in \R^N.
$$
\end{Lemma}

\begin{pf}
We only give an idea of the proof. We have 
\begin{align} \label{dia-formula}
|\nabla|u|(x)|&=\left |\Re\left(\frac{\overline u(x)}{|u(x)|}\nabla u(x)\right)
\right|  \\
& \notag
\leq  
\left|\Re\left( \left(\nabla u(x)+i\frac{\A(x/|x|)}{|x|} u(x) \right)
\frac{\overline u(x)}{|u(x)|} \right) \right|
\leq \left|\nabla u(x)+i \frac{\A(x/|x|)}{|x|} u(x) \right|   
\end{align}
for a.e. $x\in\R^N$. 
\end{pf}

\noindent An analogous result can be easily shown also for
$H^1_{*}(\Omega,\C)$-functions. The following lemma allows comparing
assumptions (\ref{positivity}) and (\ref{pos-mu10a}).

\begin{Lemma} \label{irrotational} Let $N\geq 2$ and assume that
  (\ref{eq:reg}) and (\ref{eq:transversality}) hold. Then
  $\mu_1(\A,a)\geq \mu_1(0,a)$ with equality holding if and only if
  ${\rm curl}\frac{\A}{|x|}=0$ in a distributional sense.
\end{Lemma}

\begin{pf}
The fact that $\mu_1(\A,a)\geq \mu_1(0,a)$ follows by \eqref{firsteig}
and the diamagnetic inequality on the sphere
\begin{equation}  \label{dia-sfera}
|\nabla_{\SN} |\psi|(\theta)|\leq  |\nabla_{\SN} \psi(\theta)+iA(\theta)\psi(\theta)|
\qquad {\rm for \ a.e.} \ \theta\in \SN
\end{equation}
which holds for any function $\psi\in H^1(\SN)$.
Indeed if 
$\psi_1\in H^1(\SN)$
is a nontrivial eigenfunction of $\mu_1(\A,a)$ then 
\begin{align} \label{mu1-mu1} \mu_1(\A,a)&= \frac{\int_{\SN}
    |\nabla_{\SN}\psi_1(\theta)+i\A(\theta)\psi_1(\theta)|^2 \, dS
    -\int_{\SN} a(\theta)|\psi_1(\theta)|^2 \, dS}
  {\int_{\SN} |\psi_1(\theta)|^2 \, dS} \\
  & \notag \geq \frac{\int_{\SN} |\nabla_{\SN}|\psi_1|(\theta)|^2 \,
    dS -\int_{\SN} a(\theta)|\psi_1(\theta)|^2 \, dS}{\int_{\SN}
    |\psi_1(\theta)|^2 \, dS} \geq \mu_1(0,a).
\end{align}

\noindent
We start by assuming that $\mu_1(\A,a)=\mu_1(0,a)$. Let $\psi_1$ be as in 
\eqref{mu1-mu1} so that by \eqref{dia-sfera} we infer
\begin{equation} \label{IDSN}
|\nabla_{\SN}\psi_1(\theta)+i\A(\theta)\psi_1(\theta)|
=
|\nabla_{\SN}|\psi_1|(\theta)|  
\qquad {\rm for \ a.e.} \ \theta\in \SN.
\end{equation}
Similarly to \eqref{dia-formula} we have 
\begin{align} \label{dia-formula-SN}
|\nabla_{\SN}|\psi_1|(\theta)|
\leq  
\left|\Re\left( \frac{\overline \psi_1(\theta)}{|\psi_1(\theta)|}\left(\nabla_{\SN} \psi_1(\theta)+i\A(\theta) \psi_1(\theta) \right)
 \right) \right|
\leq \left|\nabla_{\SN} \psi_1(\theta)+i \A(\theta) \psi_1(\theta) \right|   
\end{align}
which together with \eqref{IDSN} gives 
$$
\Im (\overline \psi_1(\theta) (\nabla_{\SN} \psi_1(\theta)+i\A(\theta)
\psi_1(\theta)) )=0 \qquad {\rm for \ a.e.} \ \theta\in \SN
$$
and in turn
$$
\A(\theta)=-\Im\left(\frac{\nabla_{\SN} \psi_1(\theta)}{\psi_1(\theta)} \right)
\qquad {\rm for \ a.e.} \ \theta\in \SN.
$$
This implies
$$
\frac{\A(x/|x|)}{|x|}=-\Im \left(\frac{\nabla
    (\psi_1(x/|x|))}{\psi_1(x/|x|)} \right) \qquad {\rm for \ a.e.} \
x\in \R^N.
$$
By direct computation this gives ${\rm curl} \frac{\A}{|x|}=0$ in a
distributional sense.
 
Suppose now that ${\rm curl} \frac{\A}{|x|}=0$ in a distributional
sense and let us prove that $\mu_1(\A,a)=\mu_1(0,a)$.  By \cite{Lein}
we have that there exists $\phi \in L^1_{{\rm loc}} (\R^N)$ such that
$\nabla \phi=\frac{\A}{|x|}$ in a distributional sense.  From
(\ref{eq:transversality}) it follows that
$\phi(x)=\phi(\frac{x}{|x|})$ and $\nabla_{{\mathbb S}^{N-1}}\phi=\A$.
 Let $\Psi$ be
a nontrivial eigenfunction of $\mu_1(0,a)$ and define the angular
function $\psi(\theta)$ by
$$
\psi(\theta)=e^{-i\phi(\theta)} \Psi (\theta).
$$ 
Then 
\begin{align*}
  \mu_1(\A,a)&\leq \frac{\int_{\SN}
    |\nabla_{\SN}\psi(\theta)+i\A(\theta)\psi(\theta)|^2 \, dS
    -\int_{\SN} a(\theta)|\psi(\theta)|^2 \, dS}
{\int_{\SN} |\psi(\theta)|^2 \, dS} \\
  &= \frac{\int_{\SN} |\nabla_{\SN}\Psi(\theta)|^2 \, dS -\int_{\SN}
    a(\theta)|\Psi(\theta)|^2 \, dS}{\int_{\SN} |\Psi(\theta)|^2 \,
    dS} =\mu_1(0,a).
\end{align*}
Since the reverse inequality is always
verified the proof is complete. 
\end{pf}

The following Hardy type inequality with boundary terms is due to Wang
and Zhu \cite{wz}.
\begin{Lemma}[\bf Wang and Zhu]\label{l:wz}
  For every $r>0$ and $u\in H^1(B_r,\C)$ there holds
\begin{align}\label{eq:32}
  \int_{B_r} |\nabla u(x)|^2\,dx+ \frac{N-2}{2r}\int_{\partial
    B_r}|u(x)|^2\,dS \geq
  \bigg(\frac{N-2}{2}\bigg)^{\!\!2}
  \int_{B_r}\frac{|u(x)|^2}{|x|^2}\,dx.
  \end{align}
  \end{Lemma}
\begin{pf}
See \cite[Theorem 1.1]{wz}.
\end{pf}

The following lemma establishes the relation between the classical
$H^1$-space on the sphere and its magnetic counterpart,

\begin{Lemma}\label{l:h1sph}
  If $N\geq 2$ and $\A\in L^{\infty}({\mathbb S}^{N-1},\R^N)$, then the space
  $H^1_{\A}(\mathbb S^{N-1})$ defined in
  (\ref{eq:1magn}--\ref{eq:2magn}) coincides with the Sobolev space
\begin{equation*}
  H^1(\mathbb S^{N-1},\C):=\big\{\psi\in L^2(\mathbb S^{N-1},\C):\,
  \nabla_{\mathbb S^{N-1}}\psi\in L^2(\mathbb S^{N-1},\C^N)\big\}.
\end{equation*}
Moreover the norms $\|\cdot\|_{H^1_{\A}(\mathbb S^{N-1})}$ and
\begin{equation*}
  \|\cdot\|_{H^1(\mathbb S^{N-1},\C)}:= \bigg(
  \|\nabla_{\mathbb S^{N-1}}\cdot\|_{L^2(\mathbb S^{N-1},\C^N)}^2
+ \|\cdot\|_{L^2(\mathbb S^{N-1},\C)}^2\bigg)^{\!\!1/2},
\end{equation*}
are equivalent.
\end{Lemma}

\begin{pf}
  It follows easily from boundedness of the function
  $\theta\mapsto|\bf A(\theta)|$.
\end{pf}

We finally describe the spectrum of angular operator $L_{\A,a}$.

\begin{Lemma}\label{l:spe}
  Let $a\in L^{\infty}({\mathbb S}^{N-1},\R)$ and $\A\in C^1({\mathbb
    S}^{N-1},\R^N)$.  Then the spectrum of the operator $L_{\A,a}$ on $\mathbb
  S^{N-1}$ consists in  a diverging sequence of real eigenvalues
with finite multiplicity
  $\mu_1(\A,a)\leq\mu_2(\A,a)\leq\cdots\leq\mu_k(\A,a)\leq\cdots$ the
  first of which admits the variational  characterization (\ref{firsteig}).
\end{Lemma}
\begin{pf} For $\lambda=1+\|a\|_{L^{\infty}({\mathbb S}^{N-1},\R)}$, the operator
$T:L^{2}({\mathbb S}^{N-1},\R)\to L^{2}({\mathbb S}^{N-1},\R)$ defined as
\[
Tf=u\quad\text{if and only if}\quad
\big(-i\,\nabla_{\mathbb S^{N-1}}+\A\big)^2u-au+\lambda u=f
\]
is well-defined, symmetric, and compact. The lemma follows then from
classical spectral theory.~\end{pf}


\begin{thebibliography}{99}

\bibitem{almgren} F. J. Jr. Almgren, {\it $Q$ valued functions
    minimizing Dirichlet's integral and the regularity of area
    minimizing rectifiable currents up to codimension two},
  Bull. Amer. Math. Soc., 8 (1983), no. 2, 327--328.

\bibitem{as} G. Arioli, A. Szulkin, {\it A semilinear Schr\"odinger
    equation in the presence of a magnetic field},
  Arch. Ration. Mech. Anal., 170 (2003), no. 4, 277--295.

\bibitem{balinsky} A. A. Balinsky, {\it Hardy type inequalities for
    Aharonov-Bohm magnetic potentials with multiple singularities},
  Math. Res. Lett., 10 (2003), no. 2-3, 169--176.

\bibitem{BrezisKato} H.~Br{\'e}zis, T.~Kato, {\it Remarks on the
    {S}chr\"odinger operator with singular complex potentials}, J.
  Math. Pures Appl. (9), 58 (1979), no.~2, 137--151.

\bibitem{CKN} L.~Caffarelli, R.~Kohn, L.~Nirenberg, {\it First order
    interpolation inequalities with weights}, Compositio Math., 53
  (1984), no.~3, 259--275.

\bibitem{CatrinaWang} F.~Catrina, Z.-Q. Wang, {\it On the
    {C}affarelli-{K}ohn-{N}irenberg inequalities: sharp constants,
    existence (and nonexistence), and symmetry of extremal functions},
  Comm. Pure Appl. Math., 54 (2001), no.~2, 229--258.

\bibitem{cs} J. Chabrowski, A. Szulkin, {\it On the Schr\"odinger
    equation involving a critical Sobolev exponent and magnetic
    field}, Topol. Methods Nonlinear Anal., 25 (2005), no. 1, 3--21.

\bibitem{cingolani} S. Cingolani, {\it Semiclassical stationary states
    of nonlinear Schr\"dinger equations with an external magnetic
    field}, J. Differential Equations, 188 (2003), no. 1, 52--79.

\bibitem{cingolanisecchi} S. Cingolani, S. Secchi, {\it Semiclassical
    limit for nonlinear Schr\"odinger equations with electromagnetic
    fields}, J. Math. Anal. Appl., 275 (2002), no. 1, 108--130.

\bibitem{el} M. J. Esteban, P.-L. Lions, {\it Stationary solutions of
    nonlinear Schr\"odinger equations with an external magnetic
    field}, Partial differential equations and the calculus of
  variations, Vol. I, 401--449, Progr. Nonlinear Differential
  Equations Appl., 1, Birkh\"auser Boston, Boston, MA, 1989.

\bibitem{FMT1} V. Felli, E.M. Marchini, S. Terracini, {\it On
    {S}chr\"odinger operators with multipolar inverse-square
    potentials}, Journal of Functional Analysis, 250 (2007), 265--316.

\bibitem{FMT2} V. Felli, E.M. Marchini, S. Terracini, {\it On the
    behavior of solutions to {S}chr\"odinger equations with
    dipole-type potentials near the singularity}, Discrete
  Contin. Dynam. Systems, 21 (2008), 91--119.

\bibitem{FMT3} V. Felli, E.M. Marchini, S. Terracini, {\it On
    {S}chr\"odinger operators with multisingular inverse-square
    anisotropic potentials}, Indiana Univ. Math. Journal, to appear.

\bibitem{GL} N. Garofalo, F.-H.  Lin, {\it Monotonicity properties of
    variational integrals, $A\sb p$ weights and unique continuation},
  Indiana Univ. Math. J.  35 (1986), no. 2, 245--268.

\bibitem{hhlt} M. Hoffmann-Ostenhof, T. Hoffmann-Ostenhof, A. Laptev,
  J.  Tidblom, {\it Many-particle Hardy inequalities},
  J. Lond. Math. Soc., (2) 77 (2008), no. 1, 99--114.

\bibitem{jackson} J. D. Jackson, {\it From Lorenz to Coulomb and other
    explicit gauge transformations}, Am. J. Phys., 70 (2002),
  917--928.

\bibitem{Kurata} K. ~Kurata, {\it A unique continuation theorem for
    the Schr\"odinger equation with singular magnetic field},
  Proc. Amer. Math. Soc., 125 (1997), no. ~3, 853--860

\bibitem{Kurata2} K. ~Kurata, {\it Local boundedness and continuity
    for weak solutions of $-(\nabla-ib)\sp 2u+Vu=0$},  Math. Z.,  224
    (1997), no. 4, 641--653.

\bibitem{lw} A. Laptev, T. Weidl, {\it Hardy inequalities for
      magnetic Dirichlet forms}, Mathematical results in quantum
    mechanics (Prague, 1998), 299--305, Oper. Theory Adv. Appl., 108,
    Birkh\"auser, Basel, 1999.

\bibitem{Lein} H. Leinfelder, {\it Gauge invariance of Schr\"odinger
      operators and related spectral properties}, J. Operator Theory 9
    (1983), 163--179

\bibitem{leblond} J. M. L\'evy-Leblond, {\it Electron capture by
      polar molecules,} Phys. Rev., 153 (1967), no. 1, 1--4.

\bibitem{LL} E. H. Lieb, M. Loss, {\it Analysis,} Graduate Studies in
  Mathematics 14, AMS (1997)

\bibitem{mor} M. Melgaard, E.-M. Ouhabaz, G. Rozenblum, {\it Negative
    discrete spectrum of perturbed multivortex Aharonov-Bohm
    Hamiltonians},  Ann. Henri Poincar\'e, 5 (2004), no. 5, 979--1012.

\bibitem{pinchover94} Y. Pinchover, {\it On positive Liouville
    theorems and asymptotic behavior of solutions of Fuchsian type
    elliptic operators,} Ann. Inst. H. Poincar\'e Anal. Non
  Lin\'eaire, 11 (1994), no. 3, 313--341.

\bibitem{Ter} S. ~Terracini, {\it On positive entire solutions to a
    class of equations with a singular coefficient and a critical
    exponent}, Adv. Diff. Eq., 1 (1996), no. ~2, 241--264

\bibitem{thaller} B. Thaller, {\it The Dirac equation}, Texts and
  Monographs in Physics, Springer-Verlag, Berlin, 1992.

\bibitem{wz} Z.-Q. Wang, M. Zhu, {\it {H}ardy inequalities with
    boundary terms}, Electron. J. Differential Equations 2003, No. 43.

\bibitem{wolff} T. H. Wolff, {\it A property of measures in $\R^ N$
    and an application to unique continuation},  Geom. Funct. Anal.,  2
    (1992), no. 2, 225--284.

\end{thebibliography}
\end{document}